\documentclass[10pt]{amsart}
\usepackage{amsmath, amscd, amsthm} 
\usepackage{amssymb, amsfonts} 
\usepackage{color}
\usepackage[all]{xy} 
\newcommand{\C}{\mathbb{C}} 
\newcommand{\A}{\mathbb{A}}
\renewcommand{\P}{\mathbb{P}}
\newcommand{\Z}{\mathbb{Z}}
\newcommand{\R}{\mathbb{R}}
\newcommand{\G}{\mathbb{G}}

\newcommand{\setm}{-}

\newcommand{\iso}{\cong}

\newcommand{\wkeq}{\simeq}

\newcommand{\id}{\mathrm{id}}

\newcommand{\mcal}[1]{\mathcal{#1}}
\newcommand{\sdot}{\smash\cdot}

\renewcommand{\div}{\mathrm{div}}

\DeclareMathOperator*{\colim}{\mathrm{colim}}

\DeclareMathOperator*{\tot}{\mathrm{Tot}}

\DeclareMathOperator{\spec}{\mathrm{Spec}}
\DeclareMathOperator{\proj}{\mathrm{Proj}}
\DeclareMathOperator{\stab}{\mathrm{Stab}}
\DeclareMathOperator{\cone}{\mathrm{cone}}
\DeclareMathOperator{\supp}{Supp}
\newcommand{\Sym}{\mathrm{Sym}}
\newcommand{\Ind}{\mathrm{Ind}}
\newcommand{\Div}{\mathrm{Div}}
\newcommand{\cyc}{\mathrm{cyc}}
\newcommand{\Pic}{\mathrm{Pic}}

\DeclareMathOperator{\coker}{\mathrm{coker}}
\DeclareMathOperator{\codim}{\mathrm{codim}}

\DeclareMathOperator{\Hom}{Hom}
\DeclareMathOperator{\Ext}{Ext}

\DeclareMathOperator{\Cor}{Cor}
\DeclareMathOperator{\Aut}{Aut}
\newcommand{\Sch}{\mathrm{Sch}}
\newcommand{\Sm}{\mathrm{Sm}}
\newcommand{\Ab}{\mathrm{Ab}}
\newcommand{\Shv}{\mathrm{Shv}}
\newcommand{\Pre}{\mathrm{Pre}}
\DeclareMathOperator{\ch}{char}
\DeclareMathOperator{\rank}{rank}

\numberwithin{equation}{section} 

\theoremstyle{plain}
\newtheorem{theorem}[equation]{Theorem}
\newtheorem*{theorem*}{Theorem}
\newtheorem{proposition}[equation]{Proposition}
\newtheorem{lemma}[equation]{Lemma}
\newtheorem{corollary}[equation]{Corollary}

\theoremstyle{definition}
\newtheorem{definition}[equation]{Definition}
\newtheorem{example}[equation]{Example}
\newtheorem{notation}[equation]{Notation}

\newtheorem{remark}[equation]{Remark}
\newtheorem{condition}[equation]{Condition}

\begin{document}
\title[Equivariant Cancellation]{Equivariant cycles and cancellation for motivic cohomology}
\author{J. Heller}
\email{heller@math.uni-wuppertal.de}
\address{Bergische Universit\"at Wuppertal, Gau{\ss}str. 20, D-42119 Wuppertal, Germany}
\thanks{The first author received partial support from DFG grant HE6740/1-1}
\author{M. Voineagu}
\email{m.voineagu@unsw.edu.au}
\address{UNSW Sydney, NSW 2052 Australia}
\thanks{The second author received partial support from JSPS Grant in Aid (B), No. 23740006}
\author{P. A. {\O}stv{\ae}r}
\email{paularne@math.uio.no}
\address{Department of Mathematics, University of Oslo, Norway.}
\thanks{The third author received partial support from the Leiv Eriksson 
mobility programme and RCN ES479962}

\begin{abstract}
We introduce a Bredon motivic cohomology theory for smooth schemes defined over a field and equipped with an action by a finite group. 
These cohomology groups are defined for finite dimensional representations as the hypercohomology of complexes of equivariant correspondences in the equivariant Nisnevich topology. 
We generalize the theory of presheaves with transfers to the equivariant setting and prove a Cancellation Theorem.
\end{abstract}
\maketitle
\tableofcontents

\section{Introduction}

The theory of motivic cohomology for smooth schemes over a base field is a well established one. It is a powerful computational tool with ramifications for many branches of algebra, algebraic  and arithmetic geometry: quadratic forms, algebraic $K$-theory, special values of $L$-functions,  to name a few. The success of this theory is best exemplified by its fundamental role in Voevodsky's resolution of Milnor's conjecture on Galois cohomology \cite{Voev:miln}.

The purpose of this article is to generalize Suslin and Voevodsky's  construction of motivic cohomology \cite{OrangeBook}, especially Voevodsky's machinery of presheaves with transfers \cite{V:pth}, to the equivariant setting of  smooth schemes over a field $k$ equipped with an action of a  finite group $G$ (and $|G|$ coprime to $\ch(k)$).

Using Totaro's construction \cite{Tot:BG} of an algebro-geometric version of the classifying space of an algebraic group, 
Edidin-Graham \cite{EG:int} have constructed an equivariant version of Bloch's higher Chow groups. 
This theory has proved to be interesting, amongst other reasons, for its connection to equivariant algebraic $K$-theory and the equivariant Riemann-Roch theorem \cite{EG:RR, EG:compl}.
These equivariant higher Chow groups are an algebro-geometric version of topological Borel cohomology. Our construction follows a different route altogether and results in a more refined \textit{Bredon style} cohomology theory. These Bredon motivic cohomology groups
form a new set of invariants for smooth schemes with $G$-action. As a first indication that the Bredon motivic cohomology theory we construct is a refinement of equivariant higher Chow groups, we note that the  former is equipped with a grading by the representations of $G$ while the latter is graded by integers. In fact in Section \ref{sub:examples} we construct a comparison map from our Bredon motivic cohomology to the equivariant higher Chow groups. 

Topological Bredon cohomology has recently experienced a surge of interest in part because of its appearence in the work of Hill-Hopkins-Ravenel \cite{HHR}, specifically through the equivariant slice spectral sequence for certain spectra. The
$\Z/2$-equivariant case of this spectral sequence was first constructed by Dugger in \cite{dug:kr} where he constructed a spectral sequence relating Bredon cohomology groups (with coefficients in the constant Mackey functor $\underline{\Z}$) and Atiyah's $KR$-theory. The motivic analog of Atiyah's $KR$-theory is Hermitian $K$-theory, constructed as a $\Z/2$-equivariant motivic spectrum $K\R^{alg}$ by Hu-Kriz-Ormsby \cite{HKO}. Our construction of Bredon motivic cohomology is motivated in part by  a program to construct and use a $\Z/2$-equivariant motivic generalization of Dugger's spectral sequence as a tool for studying these Hermitian $K$-theory groups.

We define our theory via hypercohomology in the equivariant Nisnevich topology, introduced in \cite{Deligne:V} and \cite{HKO}.
A surjective, equivariant \'etale map $f:Y\to X$ is an equivariant Nisnevich cover if for each point $x\in X$ there is a point $y\in Y$ such that $f$ induces an isomorphism $k(x)\iso k(y)$ on residue fields and an isomorphism $G_{y}\iso G_{x}$ on \textit{set-theoretic} stabilizers. A different generalization of the Nisnevich topology was introduced in \cite{Herrmann:EMHT}, the fixed point Nisnevich topology. A cover in this topology is as above but with the requirement that $f$ induce an isomorphism $I_{y}\iso I_{x}$ on scheme-theoretic stabilizers rather than the set-theoretic ones. (An equivalent formulation is that for each subgroup $H$, the map on fixed points $f^{H}:Y^{H}\to X^{H}$ is a Nisnevich cover.) 
We focus on the equivariant Nisnevich topology in the present work for two reasons. One is that 
equivariant algebraic $K$-theory fulfills descent in the equivariant Nisnevich topology but not in the fixed point Nisnevich topology. The second is that the fixed point Nisnevich topology does not behave well with respect to transfers, see Example \ref{ex:fpfail}.

A presheaf with equivariant transfers is a presheaf $F$ of abelian groups on smooth $G$-schemes which are equipped with functorial maps $\mcal{Z}^*:F(Y)\to F(X)$ for finite equivariant correspondences $\mcal{Z}$ from $X$ to $Y$. The presheaf with equivariant transfers freely generated by $X$ is denoted by $\Z_{tr,G}(X)$. 
We define the equivariant motivic complexes $\Z_{G}(n)$ by forming the $\A^{1}$-singular chain complex on $\Z_{tr,G}(\P(k[G]^{n}\oplus 1))/\Z_{tr,G}(\P(k[G]^{n})$, where $k[G]$ is the regular representation.
Bredon motivic cohomology is defined as equivariant Nisnevich hypercohomology with coefficients in the complex $\Z_{G}(n)$.  In light of the equivariant Dold-Thom theorem \cite{DS:equiDT}, our construction is analogous to the topological Bredon cohomology with coefficients in the constant Mackey functor $\underline{\Z}$. Our setup is however flexible enough  to allow for more sophisticated coefficient systems. The benefit of allowing Mackey functor coefficients is well understood in topology. As explained in Section \ref{sub:coefficients}, we have an embedding of cohomological Mackey functors into our category.

Over the complex numbers there is a topological realization functor which relates our Bredon motivic cohomology groups and topological Bredon cohomology groups. This comparison map is the subject of an equivariant Beilinson-Lichtenbaum type conjecture, predicting that a range of these groups should be isomorphic with torsion coefficients. The ordinary Beilinson-Lichtenbaum conjecture, relating motivic cohomology and \'etale cohomology (or singular cohomology), is equivalent to the Milnor and Bloch-Kato conjectures which have been proved in work of Voevodsky and Rost. 
In a sequel paper \cite{HVO:BL} we show that this equivariant conjecture is true for $G=\Z/2$ and any torsion coefficients.    
The key new ingredient in that work is a $\Z/2$-equivariant generalization of Voevodsky's Cancellation Theorem,   which is the main result of this present paper. The Cancellation Theorem is an algebro-geometric version of the familiar suspension isomorphism in singular cohomology of topological spaces. It asserts that in motivic cohomology, we can cancel suspension factors of the algebro-geometric sphere $\G_{m}$. Besides the usual $\G_{m}$, our equivariant version also allows for $\G_{m}$ equipped with the canonical involution. 
 
\begin{theorem}[Equivariant Cancellation]
Let $X$ be a smooth $\Z/2$-scheme over a perfect field $k$. Then
$$
H^{n}_{GNis}(X,C_*\Z_{tr,G}(Y)) = H^{n}_{GNis}(X\wedge \G, C _*\Z_{tr,G}(Y\wedge \mathbb{G}))
$$
where $\mathbb{G}$ denotes either $\G_{m}$ with the trivial action or $\G_{m}$ equipped with the involution $x\mapsto x^{-1}$.
\end{theorem}

We establish the equivariant Cancellation Theorem as Theorem \ref{thm:can} below. Its proof uses equivariant modifications of Voevodsky's arguments in \cite{Voev:can} and relies on equivariant versions of the main results of Voevodsky's techniques for analyzing the cohomological behaviour of presheaves with transfers. Most of the paper is focused on the generalization of this machinery. 

%

As mentioned, our motivation is to study $\Z/2$-equivariant phenomena, but where possible we establish our results in a greater generality. Everything works  
best under the assumption that the irreducible representations of $G$ (which are defined over the base field) are all one-dimensional,
but we anticipate further generalizations are possible. 
A main source of this condition on the group $G$  arises from the question of existence of equivariant triples. 
We refer to Section \ref{sec:etrip} for a precise definition of an equivariant triple, but point out here that it is in particular a smooth equivariant relative curve $X\to S$ between smooth $G$-schemes. 
A typical step in several key arguments used in the course of establishing the main homotopy invariance result (see Theorem \ref{thm:hi}) 
is that in order to establish an isomorphism of sheaves with equivariant transfers, it suffices to show the isomorphism on the generic points of smooth $G$-schemes. To establish this reduction step we need a good supply of equivariant triples. More precisely, if $x\in X$ is a point, there should exist an invariant open neighborhood $U\subseteq X$ of $x$ and a smooth $G$-scheme $S$ such that $U\to S$ is a smooth equivariant curve. 
As a simple illustrative example, consider a  representation $V$  viewed as a smooth $G$-scheme. Suppose there is a smooth equivariant curve $f:U\to C$, where $U$ is an invariant neighborhood of the origin in $V$.  Then there is an equivariant surjection 
$df:V\iso T_{0}U \to T_{f(0)}C$. This surjection splits (as $|G|$ and $\ch(k)$ are coprime) and so $V$ contains a one-dimensional summand. It could however happen that $V$ admits no one-dimensional summand and thus there is no neighborhood of the origin in $V$ fitting into a smooth equivariant curve. We establish the existence of triples around an arbitrary point of a smooth quasi-projective $G$-schemes under the assumption that all irreducible representations of $G$ 
are one-dimensional. 

%
%

%

Other main results are as follows.
Theorem \ref{thm:complex} provides a Mayer-Vietoris exact sequence for certain special equivariant Nisnevich covers. This has important consequences. It in particular allows for the computation of the equivariant Nisnevich cohomology of open invariant subsets of $G$-line bundles over smooth zero dimensional $G$-schemes, see Theorem \ref{thm:complex}, which is the precursor to the homotopy invariance theorem.  

In Theorem \ref{thm:hi} we establish our homotopy invariance result.
\begin{theorem}
Suppose that all irreducible $k[G]$-modules are one-dimensional. Let $F$ be a homotopy invariant presheaf with equivariant transfers on $G\Sm/k$. Then 
$H^{n}_{GNis}(- ,F_{GNis})$ is also a homotopy invariant presheaf with equivariant transfers.
\end{theorem}

Lastly we mention that several other  constructions of equivariant cohomology theories related to algebraic cycles exist. 
There is for example work of Lawson, Lima-Filho, and Michelson \cite{LLFM:equi}, Joshua \cite{Joshua:bredon}, and Levine-Serpe \cite{LS:ss}. 
It would be interesting to see how these different constructions relate to the one carried out here.

An outline of this paper is as follows. In Section \ref{sec:Gsch} we record results about $G$-schemes, equivariant divisors, bundles, and cohomology that we need in later sections. In Section \ref{sec:enis} we recall the equivariant Nisnevich topology and establish some of its basic properties not already appearing in the literature. In Section \ref{sec:pretran} we introduce equivariant finite correspondences and presheaves with equivariant transfers. 
We formally define Bredon motivic cohomology in Section \ref{sec:bred} and using the machinery developed in later sections we establish properties and some computations. 
We relate equivariant transfers and equivariant divisors in Section \ref{sec:relcart}, in particular in Theorem \ref{thm:shcurve} we compute the equivariant Suslin homology of equivariant affine curves. In Section \ref{sec:etrip} we study equivariant triples and establish a Mayer-Vietoris sequence in Theorem \ref{thm:complex}. The homotopy invariance of cohomology is established in Section \ref{sec:hi}.  Finally in Section \ref{sec:cancellation}  we establish a $\Z/2$-equivariant generalization of Voevodsky's Cancellation Theorem.

\textbf{Notations and conventions.} Throughout $k$ is a field, which is assumed to be perfect starting in Section \ref{sec:etrip} and $G$ is a finite group whose order is coprime to $\ch(k)$. The finite group $G$ is viewed as a group scheme over $k$ via $G_{k}:=\coprod_{G}\spec(k)$. Usually we simply write again $G$ for this group scheme. 
 Write $G\Sch/k$ (resp.~$G\Sm/k$) for the category whose objects are seperated schemes of finite type (resp.~smooth schemes) over $k$ equipped with a left action by $G$ and morphisms are equivariant morphisms.

We write $\A(V) = \spec (\Sym (V^{\vee}))$  for the affine scheme associated to a vector space over $k$ and $\P(V) = \proj(\Sym (V^{\vee}))$ for the associated projective scheme.
Sometimes we write $V$ for both the vector space as well as its associated scheme $\A(V)$. 

It is important to distinguish between two types of stabilizer groups of $x\in X$,
\begin{enumerate}
 \item the \textit{set-theoretic stabilizer} $G_{x}$ of $x$ is $G_{x}=\{g\in G \,|\, gx =x\}$ \item the  \textit{inertia group} of $x$ is $I_{x} = \ker(G_{x}\to \Aut(k(x)/k))$.
\end{enumerate}
Given a subset $Z\subseteq X$ we write $G\smash\cdot Z$ or $GZ$  for the \textit{orbit} $\cup_{g\in G} gZ$ of $Z$. For a nonclosed point $x\in X$, $G\smash\cdot x$ is given a scheme structure via 
$G\smash\cdot x = (G\times \{x\})/G_{x}$.

If $F$ is a presheaf on $G\Sm/k$ and $S = \lim_{i}S_{i}$ is an inverse limit of smooth $G$-schemes over $k$, with equivariant transition maps, then we set $F(S) = \colim_iF(S_i)$.

\section{Schemes with $G$-action}\label{sec:Gsch}
In this section we collect several useful facts used throughout this paper about schemes with an action of a finite group.

\subsection{Quotients by group actions}

We first recall some basic facts about quotients of schemes by finite groups, for full details see e.g., \cite{SGA1} or \cite{GIT}. 
By a quotient $\pi:X\to X/G$ we simply mean a categorical quotient. In particular quotients are unique when they exist. If $X = \spec(A)$ then $X/G = \spec(A^{G})$. More generally if $X$ is quasi-projective then a quotient $\pi:X\to X/G$ exists. 
The categorical quotient of a scheme by a finite group satisfies the following additional properties:
\begin{enumerate}
\item $\pi$ is finite and surjective,
\item The fibers of $\pi$ are the $G$-orbits of the $G$-action on $X$,
\item $\mcal{O}_{X/G} = \pi_{*}(\mcal{O}_{X})^{G}$,
\item if $Y\to X/G$ is flat, then $X\times_{X/G} Y\to Y$ is a quotient, and
\item if $|G|$ and $\ch(k)$ are coprime and $W\subseteq X$ is a closed  and invariant, then $W\to \pi(W)$ is a quotient.   
\end{enumerate}

\begin{definition}
 Say that $X$ is \textit{equivariantly irreducible} or \textit{$G$-irreducible} provided there is an irreducible component $X_{0}$ of $X$ such that $G\smash\cdot X_{0} = X$.
\end{definition}

Let $H\subseteq G$ be a subgroup and $X$ an $H$-scheme. The scheme $G\times X$ becomes an $H$-scheme under the action $h(g,x) = (gh^{-1}, hx)$ and we define  
$$
G\times^{H} X = (G\times X)/H.
$$
The scheme $G\times^{H} X$ has a left $G$-action through the action of $G$ on itself. 
Concretely $G\times^{H}X$ has the following description. Let $g_i$ be a complete set of left coset representatives. Then  $G\times^{H}X= \coprod_{g_{i}}X_i$, each $X_i$ is a copy of $X$ and $g\in G$ acts as $k:X_{i}\to X_{j}$ where $k\in H$ satisfies $gg_i=g_jk$.
The functor
$G\times^{H} - : H\Sch/k \to G\Sch/k$ is left adjoint to the restriction functor $G\Sch/k\to H\Sch/k$.

The $H$-action on $G\times X$ is free and so $\pi:G\times X \to G\times^{H} X$ is a principle $H$-bundle. In particular, $\pi$ is \'etale and surjective. It follows that if $X$ is smooth, then so is $G\times^{H} X$. This defines a left adjoint to the restriction functor $G\Sm/k\to H\Sm/k$,
$$
G\times^{H} - :H\Sm/k \to G\Sm/k.
$$

\subsection{{$G$}-sheaves and cohomology} \label{sub:Gshv}
Let $X$ be a $G$-scheme. Write $\sigma:G\times X \to X$ for the action map. 
We write $\tau$ for any  one of the Zariski, Nisnevich, or \'etale Grothendieck topologies on $X$. Write $pr_{2}:G\times X\to X$ for the projection and $\sigma:G\times X\to X$ for the action map.

\begin{definition}
Let $\mcal{F}$ be a sheaf of abelian groups on $X$. 
\begin{enumerate}
\item A \textit{$G$-linearization} of $\mcal{F}$ is an isomorphism
$
\phi:\sigma^{*}\mcal{F}\to pr_{2}^{*}\mcal{F}$
 of sheaves on $G\times X$
which satisfies the cocyle condition  
$$
[pr_{23}^{*}(\phi)]\circ[(1\times\sigma)^{*}(\phi)] = (m\times 1)^{*}(\phi)
$$ 
on $G\times G \times X$. Here $m:G\times G\to G$ is multiplication and $pr_{23}:G\times G\times X \to G\times X$ is the projection to second and third factors. 

\item A \textit{{$G$}-sheaf} (in the $\tau$-topology) on $X$ is a pair $(\mcal{F},\phi)$ consisting of a sheaf $\mcal{F}$ on $X$ and a $G$-linearization $\phi$ of $\mcal{F}$.

\item A \textit{{$G$}-module} on $X$ is a $G$-sheaf $(\mcal{M}, \phi)$ where $\mcal{M}$ is an $\mcal{O}_{X}$-module and the $G$-linearization $\phi:\sigma^{*}\mcal{M}\iso pr_{2}^{*}\mcal{M}$ is an isomorphism of $\mcal{O}_{G\times X}$-modules. Similarly a \textit{$G$-vector bundle}  on $X$ consists of a $G$-module $(\mcal{V},\phi)$ such that $\mcal{V}$ is a locally free $\mcal{O}_{X}$-module on $X$. 

\end{enumerate}
\end{definition}

We usually   write $\mcal{F}$ rather than $(\mcal{F},\phi)$ for a $G$-sheaf or module, leaving the $G$-linearization implicit.

\begin{remark}
The previous definition works for any algebraic group $G$. Our groups are always finite, in which case a $G$-linearization of $\mcal{F}$ is equivalent to the data of isomorphisms 
$\phi_{g}:\mcal{F} \xrightarrow{\iso} g_{*}\mcal{F} \,\,\,\textrm{for each $g\in G$}$ 
which are subject to the conditions
that  $\phi_{e} = id$ and $\phi_{gh} = h_{*}(\phi_{g})\circ\phi_{h}$ for all $g,h\in G$.
\end{remark}

An equivariant morphism $f:(\mcal{E},\phi_{\mcal{E}}) \to (\mcal{F},\phi_{\mcal{F}})$ of $G$-sheaves is a morphism $f$ of sheaves which is compatible with the $G$-linearizations in the sense that $\phi_{\mcal{F}}\circ \sigma^{*}f = pr_{2}^{*}f\circ \phi_{\mcal{E}}$. Write $\Ab_{\tau}(G,X)$ for the category whose objects are $G$-sheaves on $X$ and morphisms are the equivariant morphisms. The category of $G$-sheaves on $X$ has enough injectives. We have similarly the category  $\mathrm{Mod}_{G}(X)$ of $G$-modules on $X$ and $\mathrm{Vec}_{G}(X)$ of $G$-vector bundles on $X$.

\begin{remark}\label{rem:skew}
Recall that if $G$ acts on the ring $R$, the skew-group ring $R^{\#}[G]$ is defined as follows. As a (left) $R$-module it is free with basis $\{[g]\,|\, g\in G\}$. Multiplication is defined by setting
$(r_{g}[g])(r_h[h]) = r_g(g\smash\cdot r_h)[gh]$ and extending linearly. If $G$ acts trivially on $R$, then $R^{\#}[G]$ is the usual group ring $R[G]$.

If $X=\spec(R)$, then the category of $G$-modules on $X$ is equivalent to the category of modules over the skew-group ring $R^{\#}[G]$. 
\end{remark}

Given a $G$-sheaf $\mcal{M}$, the group $G$ acts on the global sections $\Gamma(X,\mcal{M})$.
Define the \textit{invariant global sections} functor $\Gamma^{G}_{X}:\Ab_{\tau}(G,X) \to \Ab$  by $\Gamma^{G}_{X}(\mcal{M}) = \Gamma(X,\mcal{M})^{G}$. 
%
The $\tau$-$G$-cohomology groups $H^{p}_{\tau}(G;X,\mcal{M})$ are defined as the right derived functors
$$
H^{p}_{\tau}(G;X, \mcal{M}) := R^{p}\Gamma^{G}_{X}(\mcal{M}).
$$

The functor $\Gamma^{G}_{X}$ can be expressed as a composition $\Gamma^{G}_{X} = (-)^{G}\circ \Gamma(X,-)$. The functor $\Gamma(X,-)$ sends injective $G$-sheaves 
to injective $G$-modules and so the Grothendieck spectral sequence for this composition 
yields the convergent spectral sequence
\begin{equation}\label{eqn:ss}
E_{2}^{p,q} = H^{p}(G, H^{q}_{\tau}(X,\mcal{M})) \Rightarrow H^{p+q}_{\tau}(G;X,\mcal{M}),
\end{equation}
where $H^{*}(G, - )$ is group cohomology. 

\begin{definition}[{\cite{GIT}}]
The equivariant Picard group $\Pic^{G}(X)$ is the group of $G$-line bundles on $X$, with group operation given by tensor product. 
\end{definition}

The equivariant Picard group has the following well known cohomological interpretation.

\begin{theorem}[Equivariant Hilbert 90]\label{thm:H90}
Let $X$ be a $G$-scheme.
\begin{enumerate}
\item There is a natural isomorphism
$$
\Pic^{G}(X) \xrightarrow{\iso} H^{1}_{Zar}(G;X,\mcal{O}^{*}_{X}).
$$
\item There are natural isomorphisms
$$
H^{1}_{Zar}(G;X,\mcal{O}^{*}_{X}) \iso H^{1}_{Nis}(G;X,\mcal{O}^{*}_{X}) \iso H^{1}_{et}(G;X,\mcal{O}^{*}_{X}).
$$
\end{enumerate}
\end{theorem}
\begin{proof}
 The second item is an immediate consequence of the first together with the spectral sequence (\ref{eqn:ss}). We sketch a proof of the first item. Consider the classifying topos $BG$, consisting of $G$-objects in the topos $(\Sch/k)^{\sim}_{Zar}$ of Zariski sheaves. The category $\Ab(BG/X)$ of abelian objects in $BG/X$ is equivalent to $\Ab_{Zar}(G,X)$ and the equivariant cohomology defined above agrees with the cohomology defined in the topos $BG$. Now $(BG/X, \mcal{O}_{X})$ is a ringed topos and the category of line bundles for $(BG/X ,\mcal{O}_{X})$ is equivalent to the category of $G$-line bundles on $X$. In general if $(E,A)$ is a ringed topos, with final object $\mathrm{pt}$ then $H^{1}_{E}(\mathrm{pt}, A^*)$ is isomorphic to equivalence classes of $(E,A)$-line bundles. Specialized to the case $(E,A)=(BG/X, \mcal{O}_{X})$ is the statement of the theorem.
\end{proof}

Theorem \ref{thm:H90} and the spectral sequence (\ref{eqn:ss}) yield the exact sequence
\begin{equation}\label{eqn:picss}
 0\to H^{1}(G,H^{0}(X,\mcal{O}^{*}_{X})) \to \Pic^{G}(X) \to (\Pic(X))^{G}\to H^{2}(G,H^{0}(X,\mcal{O}_{X}^{*})).
\end{equation}

\begin{lemma}\label{lem:picGA1}
 If $X$ is a reduced $G$-scheme then $p^{*}:\Pic^{G}(X) \to \Pic^{G}(X\times\A^{1})$ is injective, where $p:X\times\A^{1}\to X$ is the projection. If $X$ is normal then $p^{*}$ is an isomorphism.
\end{lemma}
\begin{proof}
If $X$ is reduced then $H^{0}(X,\mcal{O}^{*}_{X}) = H^{0}(X\times\A^{1}, \mcal{O}^{*}_{X\times\A^{1}})$ and $\Pic^{G}(X) \to \Pic^{G}(X\times\A^{1})$ is injective. If $X$ is normal then it is an isomorphism. The lemma then follows from (\ref{eqn:picss}) together with the five lemma.  
\end{proof}

\subsection{Divisors on {$G$}-schemes}\label{sub:div}
The notion of Cartier divisor and rational equivalence of Cartier divisors admits a straightforward equivariant generalization.
\begin{definition}\label{def:ecart}
\begin{enumerate}
\item An \textit{equivariant Cartier divisor} on $X$ is an element of $\Gamma(X,\mcal{K}_{X}^{*}/\mcal{O}_{X}^{*})^{G}$. 
Write $\Div^{G}(X)$ for this group with the group law written additively.

\item A \textit{principal equivariant Cartier divisor} is defined to be an invariant rational function on $X$; 
that is, an element in the image of $\Gamma(X,\mcal{K}^{*})^{G}$.

\item Two equivariant Cartier divisors are \textit{equivariantly rationally equivalent}, written $D\sim D'$, if $D-D'$ is principal.  
Write $\Div^{G}_{rat}(X)$ for the group of equivariant Cartier divisors modulo rational equivalence. 
\end{enumerate}
\end{definition}

A global section of $\Gamma(X,\mcal{K}_{X}^{*}/\mcal{O}_{X}^{*})$ is specified by giving an open covering $U_{i}$ and 
$f_{i}\in \Gamma(U_{i},\mcal{K}^{*})$ such that $f_{i}/f_{j} \in \Gamma(U_{i}\cap U_{j},\mcal{O}_{X}^{*})$ for all $i,j$. 
This section is $G$-invariant if $\{(U_{i},f_{i})\}$ and  $\{(gU_{i}, gf_{i})\}$ determine the same global section for all $g\in G$, 
where $gf_{i}$ is the rational function $gf_{i}(x) = f_{i}(g^{-1}x)$). 
This means  that  $\{(U_{i},f_{i})\}$ is $G$-invariant if and only if    
 $gf_{i}/f_{j} \in \Gamma(gU_{i}\cap U_{j},\mcal{O}_{X}^{*})$ for all $i,j$ and $g\in G$.

Write $\mcal{Z}^{1}(X)$ for the group of codimension one cycles on $X$. The homomorphism 
$\cyc:\Div(X) \to \mcal{Z}^{1}(X)$ is defined by $\cyc(D) = \sum_{Z\in X^{(1)}}\mathrm{ord}_{Z}(D)Z$ where $X^{(1)}$ is the set of closed integral codimension one subschemes. 
\begin{lemma}\label{lem:divagr}
 Let $X$ be a smooth $G$-scheme. Then 
$$
\cyc:\Div(X) \to \mcal{Z}^{1}(X)
$$ 
is an equivariant isomorphism. 
\end{lemma}
\begin{proof}
 Note that if $D\in \Div(X)$ then $\mathrm{ord}_{Z}(gD) = \mathrm{ord}_{g^{-1}Z}(D)$. It follows that $\cyc$ is equivariant. Since $X$ is smooth, $\cyc$ is an isomorphism. 
\end{proof}

Given an equivariant Cartier divisor the usual construction yields a $G$-line bundle. Recall that if  $D= \{(U_{i},f_{i})\}$ is a Cartier divisor, then the associated line bundle $\mcal{O}(D)$ is defined by $\mcal{O}(D)|_{U_{i}} = \mcal{O}_{U_{i}}f_{i}^{-1}$. It is straightforward to check that when $D$ is an equivariant Cartier divisor then the associated the line bundle $\mcal{O}(D)$ has a canonical $G$-linearization. We write again $\mcal{O}(D)$ for the $G$-line bundle defined by this choice of linearization.

\begin{proposition}\label{prop:rat}
Let $X$ be a regular $G$-scheme. 
\begin{enumerate}
\item The association $D\mapsto \mcal{O}(D)$ induces an injective homomorphism
$$
\Div^{G}_{rat}(X) \hookrightarrow \Pic^{G}(X),
$$ 
whose image consists of $G$-line bundles $\mcal{L}$ which admit an equivariant injection $\mcal{L}\hookrightarrow\mcal{K}_{X}$. 
\item If $G$ acts faithfully on $X$ then every $G$-line bundle admits such an injection into $\mcal{K}_{X}$. In particular $\Div^{G}_{rat}(X) = \Pic^{G}(X)$.
\end{enumerate}
\end{proposition}
\begin{proof}
 The first part is straightforward from the definitions. When $X$ has faithful action, then $H^{1}_{Zar}(G;X,\mcal{K}^{*}_{X}) = 0$ which implies that $\Div^{G}(X)\to \Pic^{G}(X)$ is surjective.
\end{proof}

When the action isn't faithful, $\Div^{G}_{rat}(X) \subseteq \Pic^{G}(X)$ can be a proper subgroup. For example, $\Pic^{G}(k)$ is isomorphic to the character group of $G$ over $k$ while $\Div^{G}_{rat}(k) = 0$.

\begin{proposition}
If $X$ is regular then $\Div^{G}_{rat}(X\times \A^{1}) = \Div^{G}_{rat}(X)$.
\end{proposition}
\begin{proof}

Let $K = \ker(G\to \Aut(X))$. Then $G/K$ acts faithfully on $X$ and $\Div^{G}_{rat}(X) = \Div^{G/K}_{rat}(X)$. Since $\Pic^{G/K}(X\times\A^{1}) = \Pic^{G/K}(X)$ the proposition follows from Proposition \ref{prop:rat}.
\end{proof}

\section{Equivariant Nisnevich topology}\label{sec:enis}
In this section we introduce the equivariant Nisnevich topology and list some of its properties. The equivariant Nisnevich topology on quasiprojective $G$-schemes was defined by Voevodsky
\cite{Deligne:V} in order to extend the functor of taking quotients by group actions to motivic spaces. More recently, the equivariant Nisnevich topology (on not necessarily quasiprojective smooth $G$-schemes) has been used by Hu-Kriz-Ormsby \cite{HKO} and Krishna-{\O}stv{\ae}r \cite{KO:K}. 
A related topology, the fixed point Nisnevich topology, was defined and studied by Herrmann in \cite{Herrmann:EMHT}. The fixed point Nisnevich topology has pleasant homotopical properties but unfortunately does not seem well suited for our constructions involving presheaves with equivariant transfers.

\subsection{Basic properties}
A \textit{cd-structure} on a category $\mcal{C}$ is a collection $\mcal{P}$ of commutative squares of the form 

$$
\xymatrix{
B\ar[r]\ar[d] & Y \ar[d]^{p} \\
A \ar[r]^{e} & X
}
$$
which are closed under isomorphism. The Grothendieck topology associated to $\mcal{P}$ is the Grothendieck topology generated by declaring all pairs $(A\to X, Y\to X)$ to be coverings.

\begin{definition}
\begin{enumerate}
\item A Cartesian square  in $G\Sch/k$
\begin{equation*}
\xymatrix{
 B\ar[r]\ar[d] & Y \ar[d]^{p} \\
 A \ar[r]^{e} & X
}
\end{equation*}
is said to be an \textit{equivariant distinguished square} if $p$ is \'etale, $e:A\subseteq X$ is an invariant open embedding and $p$ induces an isomorphism $(Y-B)_{red} \xrightarrow{\iso} (X-A)_{red}$.

\item The \textit{equivariant Nisnevich topology} on  $G\Sm/k$ (resp.~$G\Sch/k$) is the Grothendieck topology associated to the $cd$-structure defined by the equivariant distinguished squares and we write  $(G\Sm/k)_{GNis}$ (resp.~$(G\Sch/k)_{GNis}$) for the associated site.
\end{enumerate}
\end{definition}

\begin{lemma}
 A presheaf of sets $F$ is a sheaf in the equivariant Nisnevich topology if and only if $F(\emptyset) = \ast$ and for any distinguished square $Q$ as above the square
$$
\xymatrix{
 F(X)\ar[r]\ar[d] & F(Y) \ar[d] \\
 F(A) \ar[r] & F(B)
}
$$ 
is a pull back square.
\end{lemma}
\begin{proof}
 This follows from \cite[Lemma 2.9]{VV:cd}
\end{proof}

\begin{example}\label{ex:P1}
Let $V$ be a representation. Consider the equivariant open covering
of $\P(V\oplus 1)$ given by $\P(V\oplus 1)- \P(V) = \A(V) $
and  $\P(V\oplus 1)- \P(1)$. The intersection of these opens is identified with $\A(V)- 0$. We thus have an equivariant distinguished square
$$
\xymatrix{
\A(V)-0\ar[r] \ar[d] & \A(V) \ar[d] \\
\P(V\oplus 1)- \P(1) \ar[r] & \P(V\oplus 1).
}
$$
\end{example}

The standard characterizations of a Nisnevich cover in the nonequivariant setting admit an equivariant generalization.

\begin{definition}
Let $f:Y\to X$ be an equivariant morphism. An \textit{equivariant splitting sequence} for $f:Y\to X$ is a sequence of invariant closed subvarieties
$$
\emptyset = Z_{m+1}\subseteq Z_{m}\subseteq \cdots \subseteq Z_{1}\subseteq Z_{0}=X
$$
such that $f|_{Z_{i}-Z_{i+1}}:f^{-1}(Z_{i}-Z_{i+1})\to Z_{i}-Z_{i+1}$ has an equivariant section. The integer $m$ is called the length of this splitting sequence.
\end{definition}

\begin{proposition}\label{prop:nischar}
Let $f:Y\to X$ be an equivariant \'etale map between $G$-schemes. The following are equivalent.
\begin{enumerate}
 \item The map $f$ is an equivariant Nisnevich cover.
\item The map $f$ has an equivariant splitting sequence.
\item For every point $x\in X$, there is a point $y\in Y$ such that $f$ induces  an isomorphism $k(x){\iso} k(y)$ of residue fields and an isomorphism $G_{y}{\iso} G_{x}$ of set-theoretic stabilizers. 
\end{enumerate}
\end{proposition}
\begin{proof}
The proof follows along the lines of the nonequivariant arguments in \cite[Lemma 3.1.5]{MV:A1} and \cite[Proposition 2.17]{Voev:unstable}.
\begin{enumerate}
\item[(1)$\Leftrightarrow$(2)] 
Suppose that $\{V_{i}\to X\}$ is an equivariant Nisnevich cover.  Note that there is a dense invariant open subscheme $U\subseteq X$ on which $f:\coprod V_{i}\to X$ has a splitting. Indeed, this is true by definition for covers coming from distinguished squares and this property is preserved by pullbacks and by compositions. Restricting to the complement of this open and repeating the argument we construct an equivariant splitting sequence, which must stop at a finite stage because $X$ is Noetherian.

For the converse, we proceed by induction on the length of a splitting sequence. The case $m=0$ is immediate. Suppose that we have an equivariant splitting sequence of length $m$. The restriction of $f$ to $Z_{m}\times_{X}Y\to Z_{m}$ has an equivariant section $s$. Since $s$ is equivariant and \'etale, $s(Z_{m})\subseteq Z_{m}\times_{X} Y$ is an invariant open. Let $D$ be its closed complement, equipped with the induced reduced structure. 
Consider the map $\widetilde{Y}:=Y\setm D\to X$. Then $\{\widetilde{Y}\to X,\,X-Z_{m}\}$ forms an equivariant distinguished covering of $X$. 
The pullback of $f:Y\to X$ along $X-Z_{m}$ has an equivariant splitting sequence of length less than $m$ and so by induction is an equivariant Nisnevich cover. Similarly the pullback of $f$ along $\widetilde{Y}\to X$ equivariantly splits and is thus also an equivariant Nisnevich cover. It follows that $f$ itself is an equivariant Nisnevich covering.

\item[(2)$\Leftrightarrow$(3)] Suppose that $f$ has an equivariant splitting sequence. Then $x\in U_{k}= Z_{k}-Z_{k+1}$ for some $k$. Let $s$ be a section of $f$ over $U_{k}$ and let $y = s(x)$. Then one immediately verifies that $f$ induces an isomorphism $k(x) \iso k(y)$ and $G_{y}\iso G_{x}$.

For the other direction, by Noetherian induction it suffices to show that if for each generic point $\eta\in X$ there is $\eta'\in Y$ so that $f$ induces $k(\eta)\iso k(\eta')$ and $G_{\eta}\iso G_{\eta'}$ then there is an equivariant dense open $U\subset X$ such that $Y\times_{X}U\to U$ has an equivariant splitting. To show this it suffices to assume that $X$ is equivariantly irreducible. Let $\eta\in X$ be a generic point. Then there is an $\eta'\in Y$ such that 
$f:\eta'\iso \eta$ and $G_{\eta'}\iso G_{\eta}$. This implies that 
$G\smash\cdot \eta'\to G\smash\cdot \eta$ is an equivariant isomorphism. 
We have that $G\smash\cdot \eta'=\cap W'$  (resp.~$G\smash\cdot \eta_{i}$) is the intersection over all invariant opens $W'$ in $Y$ containing $\eta'$ (resp.~ all invariant opens in $X$) and so there is some invariant open $W'\subseteq Y$ such that $W'\to f(W')$ is an equivariant isomorphism. Setting $U=f(W')$ we obtain our equivariant splitting. 
\end{enumerate}
\end{proof}

Changing the condition above on stabilizers leads to the variant of the equivariant Nisnevich topology defined in \cite{Herrmann:EMHT}.
\begin{definition}\label{def:fpnis}
An equivariant \'etale map $f:Y\to X$ is a \textit{fixed point Nisnevich cover} if for each point point $x\in X$, there is a point $y\in Y$ such that $f$ induces  an isomorphism $k(x){\iso} k(y)$ of residue fields and an isomorphism $I_{y}{\iso} I_{x}$ of inertia groups. 
\end{definition}

By \cite[Lemma 2.12]{Herrmann:EMHT}, a map $Y\to X$ is a fixed point Nisnevich cover if and only if for every subgroup $H\subseteq G$, the map on fixed points $f^{H}:Y^{H}\to X^{H}$ is a Nisnevich cover. The following simple example illustrates an important difference between these two topologies. 
\begin{example}
 Let $X$ be a smooth $G$-scheme with free action. Consider the action map $G\times X^{triv}\to X$, where $X^{triv}$ is the $G$-scheme $X$ considered with trivial action. 
This is a fixed point Nisnevich cover. However, it is not equivariantly locally split and so is not an equivariant Nisnevich cover.  
\end{example}

We refer to \cite{VV:cd} for the definition of a complete, regular, and bounded cd-structure. 

\begin{theorem}
 The equivariant Nisnevich cd-structure on $G\Sm/k$ is complete, regular and bounded.
\end{theorem}
\begin{proof}
The argument is similar to that of \cite[Theorem 2.2]{Voev:unstable} for the usual Nisnevich topology. We provide a brief sketch of the details. First, since the equivariant distinguished squares are closed under pull back, it follows from \cite[Lemma 2.4]{VV:cd} that the equivariant Nisnevich cd-structure is complete. For regularity, one needs to see that for an equivariant distinguished square the square
$$ 
\xymatrix{
B\ar[r]\ar[d] & Y \ar[d] \\
B\times_{A}B \ar[r] & Y\times_{X} Y
}
$$
is also distinguished, where the horizontal arrows are the diagonal. Because an equivariant distinguished square is a square whose maps are equivariant and which is nonequivariantly a distinguished square, this follows immediately from the nonequivariant case which is verified in \cite[Lemma 2.14]{Voev:unstable}.

It remains to see that the cd-structure is bounded. For this we use the equivariant analogue of the standard density structure. That is for a smooth $G$-scheme $X$, let $D_{q}(X)$ be the set of open, invariant embeddings $U\to X$ whose complement has codimension at least $q$. The arguments of \cite[Proposition 2.10]{Voev:unstable} carry over to the equivariant case to show that equivariant cd-structure is bounded by this density structure.
\end{proof}

\begin{corollary}\label{cor:cdbd}
 Let $\mcal{F}$ be a sheaf of abelian groups on $G\Sm/k$ in the equivariant Nisnevich topology and let $X$ be a smooth $G$-scheme over $k$. Then 
$$
H^{i}_{GNis}(X,\mcal{F}) = 0
$$
for $i>\dim(X)$.
\end{corollary}

\subsection{Points}\label{sub:points}

Let $A$ be a commutative ring and $I\subseteq A$ an ideal contained in the Jacobson radical of $A$. Recall  that $(A,I)$ is said to be a \textit{Henselian pair} if for every \'etale ring map $f:A\to B$ and any $p:B\to A/I$ such that the composition $pf:A\to A/I$ equals the quotient map,  there is a lifting of 
$p$ to an $A$-homomorphism $\overline{p}:B\to A$.
We say that the pair $(A,I)$ has a $G$-action, if $A$ has a $G$-action and the ideal $I$ is invariant.  
There is a functorial Henselization of the pair $(A,I)$, consisting of a ring map $A\to A^{h}$ such that  $(A^{h}, IA^{h})$ is a Hensel pair and  
$A/I\iso A^{h}/IA^{h}$, see e.g., \cite{Raynaud:hensel}. If $(A,I)$ has $G$-action then $G$ acts on $(A^{h}, IA^{h})$ as well because $(-)^{h}$ is functorial. 

\begin{definition}\label{defn:hensG}
 Say that $S$ is a \textit{semilocal Henselian affine $G$-scheme over $k$} if $S= \spec(A^h)$, where $A^{h}$ is the Henselization of a pair $(A,I)$ where $A$ is a semilocal ring with $G$-action which is essentially of finite type over $k$, and $I$ is the Jacobson radical. We say that $S$ is \textit{smooth over $k$} if $A$ is essentially smooth over $k$. 
\end{definition}

\begin{remark}\label{rem:nishens}
 Let $S$ be a semilocal Henselian affine $G$-scheme over $k$.  Let $Z\subseteq S$ be the set of closed points and suppose that $f:Y\to S$ is an equivariant \'etale map which admits an equivariant splitting over $Z$. Then $f$ admits an equivariant splitting. Indeed, since $S$ is Henselian there exists a splitting $s:S\to Y$ extending the one over $Z$. Then $s$ is both an open and closed immersion. Thus $s$ is an isomorphism of $S$ onto its image and determines a decomposition $Y= s(S)\coprod Y'$. Note that $s(S)$ is invariant, otherwise $gz\in Y'$ for some $g\in G$ and $z\in Z$ but $s$ is equivariant on $Z$. It follows that $s(S)\subseteq Y$ is invariant and so $s$ is equivariant, being the inverse of the equivariant isomorphism 
$f|_{s(S)}$.

In particular if $S$ is a semilocal Henselian affine $G$-scheme with a single closed orbit and $Y\to S$ is an equivariant Nisnevich cover, 
then it can be refined by the trivial covering. 
\end{remark}

Suppose that $X$ is a $G$-scheme over $k$ and $x\in X$ has an invariant open affine neighborhood. Then $\mcal{O}_{X,G x}$ is a semilocal ring with $G$-action and 
$\spec(\mcal{O}^{h}_{X,Gx})$ is a semilocal Henselian affine $G$-scheme over $k$ with a single closed orbit. Any semilocal Henselian affine $G$-scheme over $k$ with a single orbit is equivariantly isomorphic to  $\spec(\mcal{O}^{h}_{A,Gx})$, for some affine $G$-scheme $A$ and $x\in A$.

In general, a point $x\in X$ might  not be contained in any $G$-invariant affine neighborhood. We can however still consider  $G\times^{G_{x}}\spec(\mcal{O}_{X,x}^{h})$. Additionally, it is always the case that
$Gx = G\times^{G_{x}}\{x\}\subseteq G\times^{G_{x}}X$ has a $G$-invariant affine neighborhood. The canonical map $\pi:G\times^{G_{x}}\mcal{O}_{X,x}^{h}\to \mcal{O}_{X,x}^{h}$ is \'etale and $G\times^{G_{x}}\{x\}\to G\smash\cdot x$ is an isomorphism, so $\pi$ is equivariantly split over $G\smash\cdot x$.

For $x\in X$ write $N(Gx)$ for the filtering category whose objects are pairs $(p:U\to X, s)$ where $p$ is \'etale and $s:Gx\to U$ is a section of $p$ over $Gx$, and $U$ is the union of its connected components which contain an element of $s(Gx)$. A morphism $(U\to X,s)$ to $(V\to X, s')$ is a map $f:U\to V$ making the evident triangles commute. 
Write  $N_{G}(Gx)$ for the filtering category whose objects are pairs $(p:U\to X, s)$ where $U$ is 
an equivariantly irreducible $G$-scheme, $p$ is an equivariant \'etale map, and $s:Gx\to U$ is an  equivariant section of $p$ over $Gx$. A morphism $(U\to X,s)$ to $(V\to X, s')$ in $N_{G}(Gx)$ is 
a map $f:U\to V$ making the evident triangles commute. We sometimes write $N(X;Gx)$ and $N_{G}(X;Gx)$ for these indexing categories if we need to be explicit about the ambient $G$-scheme containing $Gx$.

\begin{remark}\label{rem:locaff}
 Let $U$ be a $G_{x}$-invariant affine neighborhood of $x\in X$.  
 The transition maps of $N_{G}(G\times^{G_{x}}U,Gx)$ are affine, so $\lim_{V\in N_{G}(G\times^{G_{x}}U,Gx)} V$ exists in the category of $k$-schemes.
 The map $G\times^{G_{x}}U\to X$ is an equivariant \'etale neighborhood of $Gx$. In particular the map  $N_{G}(G\times^{G_{x}}U,Gx)\to N_{G}(X;Gx)$ is initial and so 
 $\lim_{N_{G}(X;Gx)}V$ exists as well (and equals 
 $\lim_{N_{G}(G\times^{G_{x}}U,Gx)} V$).
\end{remark}

\begin{proposition}\label{prop:Oh}
The forgetful functor $\phi:N_{G}(Gx)\to N(Gx)$ is initial. In particular  
$$
 \lim_{U\in N_{G}(Gx)}U \iso \spec(\mcal{O}_{G\times^{G_{x}}X,x}^{h})
 \iso G\times^{G_{x}}\spec(\mcal{O}_{X,x}^{h}).
$$
If $x\in X$ has a $G$-invariant affine neighborhood then additionally we have a canonical $G$-isomorphism $G\times^{G_{x}}\spec(\mcal{O}^{h}_{X,x}) \xrightarrow{\iso} \spec(\mcal{O}^{h}_{X,Gx})$.
\end{proposition}
\begin{proof}
We need to show that the comma category $(\phi/(p,s))$ is nonempty and connected for any $(p:U\to X, s)\in N(Gx)$. It suffices to check that it is nonempy since if $\phi(q_1,s_1)\to (p,s)$ and $\phi(q_2,s_2)\to (p,s)$ are two arrows in $(\phi/(p,s))$, there is $(q_3:V_3\to X,s_3)$ in $N_{G}(Gx)$ which maps to $(q_1,s_1)$ and $(q_2,s_2)$. The two maps $(q_3,s_3)\to (p,s)$ obtained from composition  agree on each point of $s(Gx)$ and induce the same map on the residue fields of these points. Each connected component of $V_3$ contains a point of $s(Gx)$ and so both maps $(q_3,s_3)\to (p,s)$ are equal.

Let $(p:U\to X, s)\in N(Gx)$. For $g\in G$, define $p^{g}:g_{*}U\to X$ to be the \'etale $X$-scheme given by $g_{*}U= U$ and $p^{g} := gp$. The identity $U=U$ can be viewed as a map $g_{*}U \to (hg)_{*}U$ over $h:X\to X$.

Label the elements of $G$ by $e = g_{0}, g_{1},\ldots g_{n}$. Define $W$ to be the $(n+1)$-fold fiber product 
$$
W= U\times_{X}(g_{1})_{*}U\times_{X}\cdots \times_{X}(g_{n})_{*}U.
$$
 Write $\pi_{g_{i}}:W\to (g_{i})_{*}U$ for the projection and consider $W$ as an $X$-scheme via the composite $p\pi_{e}:W\to U \to X$. Note that $hp\pi_{h} = p\pi_{e}$.

Now $W$ has a $G$-action given by permuting the factors. In other words we define $h:W\to W$ to be the map determined by the formula $\pi_{g_{i}}h = \pi_{h^{-1}g_{i}}$. This determines a map since $p^{g_{i}}\pi_{h^{-1}g_{i}}= p^{g_{j}}\pi_{h^{-1}g_{j}}$. Moreover $p^{g_{i}}\pi_{h^{-1}g_{i}} = hp\pi_{e}$ and thus $p\pi_{e}:W\to X$ is an equivariant \'etale map.

Define now $s':Gx\to W$ to be the map determined in the $g_{i}$th coordinate by $sg_{i}^{-1}:Gx \to (g_{i})_{*}U$. This defines an equivariant section of $p\pi_{e}:W\to X$ over $Gx\subseteq X$. 
and so $(W\to X, s')\in N_{G}(Gx)$. Now $\pi_{e}$ determines a map $(W\to X, s')\to (U\to X,s)$ in $N(Gx)$ and so $(\phi/(p,s))$ is nonempty, which completes the proof.
\end{proof}

For a $G$-scheme $X$ and $x\in X$  write $p^*_xF=F(\mcal{O}^{h}_{G\times^{G_{x}}X,Gx}) = \colim_{U\in N_{G}(Gx)}F(U)$. 
This defines a fiber functor from the category of sheaves to sets, i.e., it commutes with colimits and finite products and so determines a point of the equivariant Nisnevich topos. Every affine semilocal $G$-scheme with a single closed orbit determines such a point. By the previous paragraphs, any such $S$ is of the form $G\times^{G_{x}}\spec(\mcal{O}^{h}_{X,x})= \spec(\mcal{O}^{h}_{G\times^{G_{x}}X,Gx})$ for an appropriate $G$-scheme $X$. 
We now verify that these points form a conservative set of points.

\begin{theorem}
 A map $\phi:F_1\to F_{2}$ of sheaves of sets on $G\Sch/k$ (resp.~ $G\Sm/k$)  is an isomorphism if and only if 
$F_{1}(S) \to F_{2}(S)$ 
is an isomorphism for all (resp.~ all smooth) semilocal affine $G$-schemes $S$ over $k$ with a single closed orbit. 
\end{theorem}
\begin{proof}
%
%
%
%

Let $S$ be a semilocal affine $G$-scheme over $k$ with a single closed orbit. 
If $\phi$ is an isomorphism of sheaves then  it induces an isomorphism 
$F_1(S)\iso F_{2}(S)$ as $S_{GNis}$ is trivial.

Suppose that $\phi$ induces isomorphisms 
$F_{1}(S) \iso F_{2}(S)$ for all $S$.
We first show that $\phi$ is a monomorphism. Suppose that $\phi(\alpha) = \phi(\beta)$ for some $\alpha,\beta \in F_{1}(X)$. Then $[\alpha]=[\beta]\in F_1(\mcal{O}^{h}_{G\times^{G_{x}}X,Gx})$ for all $x\in X$. Thus for each $x\in X$, there is some equivariant \'etale map $V_{x}\to X$ which admits an equivariant section over $Gx \to V$ and $\alpha|_{V} = \beta|_{V}$. Let $\eta_{i}\in X$ be generic points. Then $V_{\eta_{i}}\to X$ has an equivariant section over an invariant open $U_{1}\subseteq X$. Consider $V_{\alpha_{i}}$ corresponding to generic points of $Z_{1}=X-U_{1}$ and let $Z_{2}\subseteq Z_{1}$ the complement of the set where $\coprod V_{\alpha_{i}}$ has a section. Proceeding in this way we obtain a finite number of equivariant \'etale maps $V_{x_{i}}\to X$ such that $V= \coprod V_{x_{i}} \to X$ is a Nisnevich cover.
Moreover $V$ has the property that $\alpha|_{V} = \beta|_{V}$ and because $F_1$ is a sheaf this means that $\alpha= \beta$ in $F_1(X)$.

Now we show that $\phi$ is a surjection. Let $\alpha\in F_{2}(X)$. For any $x\in X$ there is $[\beta]\in F_1(\mcal{O}^{h}_{G\times^{G_{x}}X,Gx})$ such that $\phi([\beta]) = [\alpha]\in F_2(\mcal{O}^{h}_{G\times^{G_{x}}X,Gx})$. Thus for each $x\in X$ there is an equivariant \'etale map $f_{x}:V_{x}\to X$, which admits an equivariant section $Gx\to V_{x}$, and $\beta_{x}\in F_1(V_{x})$ such that $\phi(\beta_{x}) = \alpha|_{V_{x}}$.  As in the previous paragraph we can find a finite number of points $x_{1},\ldots, x_{n}$ such that $V=\coprod V_{x_{i}}\to X$ is an equivariant Nisnevich cover. The elements $\beta_{x_{i}}$ determine the element $\beta\in F_1(V)$ with the property that $\phi(\beta) = \alpha|_{V}$ and thus $\phi$ is surjective as well.
\end{proof}

\begin{remark}\label{rem:fppoints}
 By \cite{Herrmann:EMHT}, the points of the fixed point Nisnevich topology are the semilocal affine $G$-schemes of the form $G/H\times \spec(\mcal{O}^{h}_{X,x})$ where $H\subseteq G$ is a subgroup, and $x$ is a point of a smooth scheme $X$ equipped with trivial action.
\end{remark}

Since $(G\Sm/k)_{GNis}$ has enough points we can form the Godement resolution of a presheaf.

\begin{definition}\label{def:god}
Let $F(-)$ be a presheaf of chain complexes of abelian groups on $G\Sm/k$.  Let
\begin{equation*}
G^{0}F(U) = \prod_{u\in U}F(\mcal{O}^{h}_{G\times^{G_{u}}U,u}).
\end{equation*}
Define $G^{n}F =G\circ\cdots \circ GF$ to be the $(n+1)$-fold composition of $G$. The inclusions and projections of the various factors give $n\mapsto G^{n}F(U)$ a cosimplicial structure.  
The \textit{Godement resolution} $F(-) \to \mcal{G}F(-)$ is defined by 
$$
\mcal{G}F(U) := \tot G^{\bullet}F(U).
$$ 
\end{definition}

Then $F\to \mcal{G}F$ is a flasque resolution of $F$ on $(G\Sm/k)_{GNis}$. Consequently $\mcal{G}F(U)$ computes the hypercohomology with coefficients in $F$:
\begin{equation*}
H^{k}\mcal{G}F(U) = H^{k}_{GNis}(U,F).
\end{equation*}

\subsection{Change of groups}\label{sub:grpchg}

Let $H\subseteq G$ be a subgroup. We have an adjoint pair of functors $\epsilon: H\Sch/k \rightleftarrows G\Sch/k: \rho$ where $\epsilon(X) = G\times^{H} X$ and $\rho (W) = W$. These restrict to an adjoint pair
\begin{equation}\label{eqn:grpchg}
\epsilon: H\Sm/k \rightleftarrows G\Sm/k: \rho.
\end{equation}
 Both functors commute with fiber products and send covering families to covering families for the equivariant Nisnevich topologies. 
We thus have adjoint functors
$$
\epsilon^{*}: \Shv_{HNis}(H\Sm/k) \rightleftarrows \Shv_{GNis}(G\Sm/k): \epsilon_{*}
$$
and 
$$
 \rho^{*}: \Shv_{GNis}(G\Sm/k) \rightleftarrows \Shv_{HNis}(H\Sm/k): \rho_{*},
$$
where $\epsilon_{*}F(X) = F(G\times^{H} X)$ and $\rho_{*}K(W) = K(W)$, and similarly for the categories of sheaves on $G\Sch/k$ and $H\Sch/k$. Additionally we have that $\rho^{*} = \epsilon_{*}$. It follows that $\epsilon_{*}$ is exact and so we have the following.
\begin{lemma}
 Let $X$ be an $H$-scheme. Then 
$$
H^{i}_{GNis}(G\times^{H} X, F) = H^{i}_{HNis}(X, \epsilon_*F).
$$
\end{lemma}

If we restrict our attention to the category $GQP/k$ of quasiprojective schemes of finite type over $k$ equipped with left $G$-action we have an adjoint pair
\begin{equation*}
\lambda: GQP/k \rightleftarrows QP/k: \eta,
\end{equation*}
where $\lambda(X) = X/G$ is the quotient by the $G$-action and $\eta(W) = W^{triv}$, where $W^{triv}$ is the scheme $W$ equipped with the trivial action. The functor $\eta$ commutes with fiber products and sends covering families to covering families. By \cite[Proposition 43]{Deligne:V} the functor $\lambda$ induces a continuous map of sites $(QP/k)_{Nis}\to (GQP/k)_{GNis}$, i.e., the presheaf $X\mapsto F(X/G)$ on $(GQP)_{GNis}$ is a sheaf  whenever $F$ is a sheaf on $(QP/k)_{Nis}$. 
We thus have adjoint functors
$$
\eta^{*}: \Shv_{Nis}(QP/k) \rightleftarrows \Shv_{GNis}(GQP/k): \eta_{*}
$$
and 
$$
 \lambda^{*}: \Shv_{GNis}(GQP/k) \rightleftarrows \Shv_{Nis}(QP/k): \lambda_{*},
$$
where $\eta_{*}F(W) = F(W^{triv})$ and $\lambda_{*}K(X) = K(X/G)$. Additionally we have that $\eta^{*} = \lambda_{*}$. It follows that $\lambda_{*}$ is exact and so we have the following.
\begin{lemma}\label{lem:qpquot}
 Let $X$ be a quasiprojective $G$-scheme. Then 
$$
H^{i}_{GNis}(X, \lambda_*F) = H^{i}_{Nis}(X/G, F).
$$
\end{lemma}


%

\section{Presheaves with equivariant transfers}\label{sec:pretran}
Let $G\Cor_{k}$ denote the category whose objects are smooth $G$-varieties and morphisms are equivariant finite correspondences, that is
$$
G\Cor_{k}(X,Y) := \Cor_{k}(X,Y)^{G}.
$$ 

An \textit{elementary equivariant correspondence} from $X$ to $Y$ is a correspondence of the form $\mcal{Z} = Z + g_{1}Z + \cdots + g_{k}Z$, where $Z\subseteq X\times Y$ is an elementary correspondence and $g_{i}$ range over a set of coset representatives for $\stab(Z) = \{g\in G | g(Z) = Z\}$.
The group $G\Cor_{k}(X,Y) $ is the free abelian group generated by the elementary equivariant correspondences.

There is an embedding of categories $G\Sm/k\subseteq G\Cor_{k}$ which sends an equivariant map $f:X\to Y$ to its graph $\Gamma_{f}\subseteq X\times Y$.

\begin{definition}
A \textit{presheaf with equivariant transfers} on $G\Sm/k$ is an additive presheaf $F:G\Cor_{k}^{op}\to \Ab$. 
\end{definition}

\begin{definition}
 \begin{enumerate}
  \item An \textit{elementary equivariant $\A^{1}$-homotopy} between two maps in $G\Sm/k$ (resp.~in $G\Cor_{k}$) $f_0,f_1:X\to Y$ is a map $H:X\times \A^{1}\to Y$ in $G\Sm/k$ (resp.~in $G\Cor_{k}$) such that $H|_{X\times\{i\}} = f_{i}$. 
  
  \item A map $f:X\to Y$ is said to be an \textit{elementary equivariant 
$\A^1$-homotopy equivalence} if there is a map $g:Y\to X$ such that both $fg$ and $gf$ are elementary equivariant $\A^{1}$-homotopic to the identity.

\item If $F$ is a presheaf on $G\Sm/k$ or on $G\Cor_k$, we say that $F$ is \textit{homotopy invariant} if the projection $p:X\times\A^{1}\to X$ induces an isomorphism $p^{*}:F(X)\xrightarrow{\iso} F(X\times \A^{1})$. 
 \end{enumerate}
\end{definition}

A simple but useful consequence of homotopy invariance is that all representations are contractible.
\begin{proposition}
 Let $F$ be homotopy invariant presheaf of abelian groups and $V$  a finite dimensional representation. Then $p^{*}:F(X\times \A(V)) \to F(X)$ is an isomorphism. 
\end{proposition}
\begin{proof}
 The map $\A(V)\times \A^{1} \to \A(V)$, $(v,t) \mapsto tv$ is an equivariant homotopy between the identity on $\A(V)$ and $\A(V)\to \{0\} \subseteq \A(V)$. 
\end{proof}

Every smooth $G$-scheme $Y$ represents a presheaf with equivariant transfers which we write as 
$$
\Z_{tr,G}(Y)(-) = G\Cor_{k}(-,Y).
$$
Note that this is in fact a sheaf in the equivariant Nisnevich topology.
If $Y$ is quasi-projective and $X$ is normal, then the map
$$
\Cor_{k}(X,Y) \xrightarrow{}  \Hom_{\Sch/k}(X, \Sym (Y))^{+}
$$  
becomes an isomorphism after inverting the exponential characteristic,
see e.g., \cite[Theorem 6.8]{SV:hom} or \cite[Proposition 2.1.3]{BV:V}.

\begin{example}\label{ex:norm}
 The sheaf $(\mcal{O}^{*})^{G}$ of invariant invertible functions is a presheaf with equivariant transfers which can be seen using Lemma \ref{lem:eqnorm}. Alternatively one may describe the transfer structure as follows. 
The sheaf $\mcal{O}^{*}$ is represented by the group scheme $\G_{m}$. We have an induced monoid morphism $\rho:\Sym(\G_{m}) \to \G_{m}$. Let $\mcal{W}:X\to \Sym (Y)$ be an effective finite correspondence. Define  $\mcal{W}^{*}:\G_{m}(Y) \to \G_{m}(X)$ by $\mcal{W}^{*}(\phi) = \rho\Sym(\phi) \mcal{W}$. It is immediate from this definition that $\mcal{W}^{*}(\phi)$ is equivariant  whenever $\mcal{W}$ and $\phi$ are  equivariant. 

\end{example}

\begin{lemma}\label{lem:corchange}
 Let $H\subseteq G$ be a subgroup. 
There is an adjunction
\begin{equation}
\epsilon: H\Cor_k \rightleftarrows G\Cor_k: \rho
\end{equation}
where $\epsilon(X) = G\times^{H} X$ and $\rho (W) = W$.

\end{lemma}
\begin{proof}

We need to show that if $X$ is a smooth $H$-scheme and $W$ is a smooth $G$-scheme then  $H\Cor_{k}(X,W)\iso G\Cor_{k}(G\times^{H}X,W)$. 
We have an $H$-equivariant map $i:X\to G\times ^{H}X$, induced by $x\mapsto (e,x)$. This gives $i^{*}:G\Cor_{k}(G\times^{H}X,W) \to H\Cor_{k}(X,W)$. It is straightforward to check that this is an isomorphism. 
\end{proof}

A presheaf $F$ with equivariant transfers is an \textit{equivariant Nisnevich sheaf with transfers} provided that the restriction of $F$ to $G\Sm/{k}$ is a sheaf in the equivariant Nisnevich topology. We finish this section with a discussion of the relation between transfers and sheafification.

\begin{theorem}\label{thm:gnisres}
 Let $X$ be a smooth $G$-scheme and $p:Y\to X$ an equivariant Nisnevich cover. Then
$$
\cdots \xrightarrow{p_{0}-p_{1}+p_{2}}\Z_{tr,G}(Y\times_{X}Y)\xrightarrow{p_{0}-p_{1}}\Z_{tr,G}(Y) \xrightarrow{p}\Z_{tr,G}(X)\to 0 
$$
is exact as a complex of equivariant Nisnevich sheaves.
\end{theorem}
\begin{proof}
The argument is similar to the nonequivariant argument \cite[Proposition 6.12]{MVW}. 
It suffices to check that the complex 
\begin{equation}\label{eqn:exact1}
\cdots \to \Z_{tr,G}(Y\times_{X} Y)(S)\to \Z_{tr,G}(Y)(S) \to \Z_{tr,G}(X)(S) \to 0
\end{equation}
is exact for  every semilocal Henselian affine $G$-scheme $S$ with a single closed orbit.
Let $Z\subseteq X\times S$ be an invariant closed subscheme which is quasi-finite over $S$.
 Define $L(Z/S)$ to be the free abelian group generated by the irreducible components of $Z$ which are finite and surjective over $S$. The assignment $Z\to L(Z/S)^{G}$ is covariantly functorial for equivariant maps of quasi-finite $G$-schemes over $S$. The sequence (\ref{eqn:exact1}) is a filtered colimit of sequences of the form
\begin{equation}\label{eqn:exact2}
\cdots \to L(Z_{Y}\times_{Z}Z_{Y}/S)^{G} \to L(Z_{Y}/S)^{G} \to L(Z/S)^{G} \to 0
\end{equation}
where the colimit is taken over all invariant closed subschemes $Z\subseteq X\times S$ which are finite and surjective over $S$. It therefore suffices to show that (\ref{eqn:exact2}) is exact. Since $S$ is a semilocal Henselian affine $G$-scheme over $k$ with a single closed orbit and $Z$ is finite over $S$ it is also Hensel semilocal. The equivariant Nisnevich covering $Z_{Y}\to Z$ therefore splits equivariantly (see Remark \ref{rem:nishens}). Let $s_{1}:Z\to Z_{Y}$ be a splitting. Set $(Z_{Y})^{k}_{Z} = Z_{Y}\times_{Z}\cdots \times_{Z}Z_{Y}$. We obtain contracting homotopies $s_{k}:L((Z_{Y})^{k}_{Z}/S)^{G} \to L((Z_{Y})^{k+1}_{Z})^{G}$ by letting $s_{k}$ be the map induced by $s_{1}\times_{Z}id _{(Z_{Y})^{k}_{Z}}$, which completes the proof.
\end{proof}

The previous statement fails when we replace the equivariant Nisnevich topology with the fixed point Nisnevich topology. 
(The following is also an example of a fixed point Nisnevich covering for which equivariant $K$-theory does not satisfy descent).
\begin{example}\label{ex:fpfail}
Consider the $\Z/2$-schemes $X=\spec(\C)$ over $\spec(\R)$ with conjugation action and $X^{triv}$ the scheme $X$ with trivial action. 
Let $Y=\Z/2\times X^{triv}$. 
The action map $Y\to X$ is then a fixed point Nisnevich cover.
In \cite{Herrmann:EMHT} it is shown that the points of the fixed point Nisnevich topology are of the form $G/H\times \mcal{O}^{h}_{W^{H},w}$. 
In particular, 
if 
\begin{equation*}\label{eqn:fp}
\cdots \to\Z_{tr,G}(Y\times_{X}Y)\to\Z_{tr,G}(Y) \to\Z_{tr,G}(X)\to 0
\end{equation*}
were to be exact in the fixed point Nisnevich topology, 
then its restriction  to $\Sm/\R$ would be exact in the Nisnevich topology.
But its restriction to $\Sm/\R$ is 
$$
\cdots \to\Z_{tr}(\C\times_{\R}\C)\to\Z_{tr}(\C) \to\Z_{tr}(\R)\to 0,
$$
which is not exact in the Nisnevich topology. 
Indeed if it were exact then applying $\Z/2(n)$, 
the complex computing weight $n$ motivic cohomology with mod-$2$ coefficients, 
would imply a quasi-isomorphism $\Z/2(n)(\R) = \Z/2(n)(\C)^{h\Z/2}$ and then we would have
$$
H^{i}_{\mcal{M}}(\R,\Z/2(n)) = H^{i}\Z/2(n)(\R) = H^{i}\Z/2(n)(\C)^{h\Z/2} = H^{i}_{et}(\R,\Z/2)
$$
for all $i\geq 0$, 
which is not true.
\end{example}

\begin{lemma}\label{lem:niscor}
 Let $p:U\to Y$ be an equivariant Nisnevich cover and $f:X\to Y$ an equivariant finite correspondence. Then there is an equivariant Nisnevich covering $p':V\to X$ and an equivariant finite correspondence $f':V\to U$ which fit into the following commutative square in $G\Cor_{k}$,
$$
\xymatrix{
V\ar[d]_{p'}\ar[r]^{f'} & U \ar[d]^{p} \\
X\ar[r]^{f} & Y .
}
$$
\end{lemma}
\begin{proof}
 We may assume that $f$ is an equivariant elementary correspondence. Write $Z= \supp(f)$ and consider the pullback $Z'= Z\times_{Y} U\subseteq X\times U$. Then $Z'\to Z$ is an equivariant Nisnevich cover and $\pi:Z\to X$ is finite. We can find an equivariant Nisnevich cover $V\to X$ such that $V\times_{X}Z'\to V\times_{X} Z$ has a section. 
 
Let $s$ be such a section. Then $s(V\times_{X}Z)\subseteq V\times U$ is finite and equivariant over $V$ and its associated equivariant correspondence defines the required $f':V\to U$.
\end{proof}

\begin{theorem}\label{thm:gnistrans}
 Let $F$ be a presheaf with equivariant transfers on $G\Sm/k$. Then $F_{GNis}$ has a unique structure of a presheaf with equivariant transfers such that $F\to F_{GNis}$ is a morphism of presheaves with equivariant transfers.
\end{theorem}
\begin{proof}
 The proof is parallel to the proof of \cite[Theorem 6.17]{MVW}. We begin with uniqueness. Let $F_{1}$ and $F_{2}$ be two presheaves with transfers with a map of presheaves with equivariant transfers 
 $F\to F_{i}$ whose underlying map of presheaves is the canonical map 
 $F\to F_{GNis}$. 
 Let $f:X\to Y$ be a map in $G\Cor_{k}$ and $y\in F_{1}(Y) = F_{2}(Y) = F_{GNis}(Y)$. Choose an equivariant Nisnevich covering $U\to Y$ such that $y|_{U}$ is in the image of $u\in F(U)$. Applying Lemma \ref{lem:niscor} we have a commutative square in $G\Cor_{k}$
$$
\xymatrix{
V\ar[d]_{p'}\ar[r]^{f'} & U \ar[d]^{p} \\
X\ar[r]^{f} & Y 
}
$$
where $p':V\to X$ is an equivariant Nisnevich cover. It is straightforward to verify, using this square, that $F_{1}(f)(y) = F_{2}(f)(y)$ and so $F_{1}=F_{2}$ as presheaves with equivariant transfers.

For existence, we first define a map
$F_{GNis}(Y) \to \Hom_{Sh}(\Z_{tr,G}(Y), F_{GNis})$, natural for maps in $G\Cor_{k}$ and such that the following square commutes
$$
\xymatrix{
F(Y) \ar[r]\ar[d] & \Hom_{\Pre(G\Cor_{k})}(\Z_{tr,G}(Y), F) \ar[d] \\
F_{GNis}(Y) \ar[r] & \Hom_{\Shv(G\Cor_{k})}(\Z_{tr,G}(Y), F_{GNis}).
}
$$

Given $y\in F_{GNis}(Y)$ there is an equivariant Nisnevich cover $U\to Y$ such that $y|_{U}$ is the image of $u\in F(U)$. The element $u$ determines a morphism $\Z_{tr,G}(U)\to F$ of presheaves with equivariant transfers. By shrinking $U$ we may assume that $u$ restricts to the zero map  $\Z_{tr,G}(U\times_{Y}U) \to F$ under the difference map $(p_1)_{*}-(p_2)_{*}:\Z_{tr,G}(U\times_{Y}U) \to \Z_{tr,G}(U)$. This in turn implies that the induced morphism of sheaves 
$\Z_{tr,G}(U) \to F_{GNis}$ also restricts to zero under the difference map.

Now Theorem \ref{thm:gnisres} implies that that $u$ determines a map $[y]:\Z_{tr,G}(Y) \to F_{GNis}$ and it is straightforward to verify that this is independent of the choice of $U$ and $u$. We now define $G\Cor_{k}(X,Y)\otimes F_{GNis}(Y) \to F_{GNis}(X)$ as follows. Let $f:X\to Y$ be a finite equivariant correspondence and $y\in F_{GNis}(Y)$. Consider the composition $[y]f:\Z_{tr,G}(X)\to \Z_{tr,G}(Y) \to F_{GNis}$ and define the pairing by sending $(f,y)$ to the image of the identity in $\Z_{tr,G}(X)(X)$ in $F_{GNis}(X)$ under $[y]f$. 
\end{proof}

A presheaf $F$ with equivariant transfers is said to be an \textit{equivariant Nisnevich sheaf} with equivariant transfers if the restriction of $F$ to $G\Sm/k$ is a sheaf in the equivariant Nisnevich topology. 
We write $\Shv(G\Cor_{k})$ for the category of sheaves with equivariant transfers in the equivariant Nisnevich topology.

\begin{corollary}
 The category $\Shv(G\Cor_{k})$ is an abelian category with enough injectives and the inclusion $i:\Shv(G\Cor_{k}) \to \Pre(G\Cor_{k})$ has a left adjoint $a_{GNis}$ which is exact and commutes with the forgetful functor to (pre)sheaves on $G\Sm/k$. 
\end{corollary}

\begin{theorem}\label{thm:gniscoh}
 Let $F$ be an equivariant Nisnevich sheaf with equivariant transfers. Then
\begin{enumerate}
 \item the cohomology presheaves $H^{n}_{GNis}(-,F)$ are presheaves with equivariant transfers,
\item for any smooth $X$, $F(X) = \Hom_{\Shv(G\Cor_{k})}(\Z_{tr}(X), F)$, and 
\item for any smooth $X$,
$$
\mathrm{Ext}^{n}_{\Shv(G\Cor_{k})}(\Z_{tr,G}(X), F) = H^{n}_{GNis}(X, F).
$$
\end{enumerate}
\end{theorem}
\begin{proof}
Suppose that $F$ is a sheaf with equivariant transfers. Then the  Godement resolution $F\to \mcal{G}F$ in Definition \ref{def:god} is a resolution of sheaves with equivariant transfers by the same reasoning as in \cite[Example 6.20]{MVW}. This implies the first statement. The second statement follows from the previous corollary together with the Yoneda lemma. For the third item, note that if $F$ is an injective sheaf with equivariant transfers then $F\to G^{0}F$ is split and so $H^{n}_{GNis}(X, F)$ is a summand of $H^{n}_{GNis}(X,G^{0}F)=0$. It follows that $H^{n}_{GNis}(X,F) = 0$ whenever $F$ is an injective and the result follows. 
\end{proof}

\section{Bredon motivic cohomology}\label{sec:bred}
In this section we introduce our Bredon motivic cohomology, explain how Mackey functors naturally appear in this setting, and give some examples of our theory. We will often need to assume the following condition.

\begin{condition}\label{conmod}
All irreducible $k[G]$-modules are one dimensional.
\end{condition}

If $G$ satisfies this condition it is necessarily abelian. Note that if $G$ satisfies Condition \ref{conmod} over $k$ then  it also does so over every field extension of $k$ as  do all of its subgroups.
\begin{lemma}
 Let $G$ be a finite abelian group and suppose that $k$ contains a primitive $n$th-root of unity where $n$ is the exponent of $G$ (i.e. the least common multiple of the orders its elements). Then Condition \ref{conmod} is satisfied.
\end{lemma}
\begin{proof}
For an abelian group $G$ Condition \ref{conmod} is equivalent to the condition that $k$ is a splitting field for $G$ (i.e. if $W$ is a simple $k[G]$-module then $W_L$ is a simple $L[G]$-module for any field extension $L/k$).  
 The lemma is thus a special case of a theorem of Brauer
 \cite[Theorem 41.1,Corollary 70.24]{CR}.
\end{proof}
%

\subsection{Definition and first properties}
If $F$ is a presheaf of abelian groups on $G\Sm/k$, write $C_{n}F(X) = F(X\times\Delta^{n}_{k})$, where $\Delta^{n}_{k}$ is  the standard algebraic simplex.
This gives a presheaf of simplicial abelian groups $n\mapsto C_{n}F(X)$ and thus yields an associated chain complex $C_{*}F$. We write $C^{*}F$ for the associated 
cochain complex,  $C^{-k}F = C_{k}F$. If $A$ is a cochain complex then the shifted complex $A[n]$ is the complex $A[n]^{i} = A^{i+n}$.
\begin{definition}
\begin{enumerate}
\item Let $V$ be a finite dimensional representation. Define $\Z_{G}(V)$ to be the complex of presheaves with equivariant transfers given by
$$
\Z_{G}(V) := C^{*}\left(\Z_{tr,G}(\P(V\oplus 1))/\Z_{tr,G}(\P(V))\right)[-2|V|].
$$
\item When $V=k[G]^{n}$ we adopt the notation
$$
\Z_{G}(n) = \Z_{G}(k[G]^{\oplus n})
$$
\end{enumerate}
\end{definition}

By virtue of their definition, the complexes $\Z_{G}(V)$ are acyclic in degrees larger than $2|V|$. In particular $\Z_{G}(n)$ is acyclic in degrees larger than $2n|G|$.

\begin{definition}
Let $X$ be a smooth $G$-variety. Define the \textit{Bredon motivic cohomology} of $X$ to be
$$
H^{n}_{G}(X,\Z(m)) = H^{n}_{GNis}(X,\Z_{G}(m)).
$$
More generally we write
$$
H^{n}_{G}(X,\Z(V)) = H^{n}_{GNis}(X,\Z_{G}(V)).
$$
\end{definition}

\begin{remark}
 By Corollary \ref{cor:cdbd} all objects of $G\Sm/k$ have finite equivariant Nisnevich cohomological dimension. This implies that the displayed hypercohomology groups (whose coefficients are unbounded complexes) are well defined, see \cite[Corollary 10.5.11]{Weibel}. 
\end{remark}

\begin{lemma}\label{lem:basicC*}
 \begin{enumerate}
  \item If $F$ is a presheaf and $f_{0},f_{1}:X\to Y$ are  elementary equivariant $\A^{1}$-homotopic then the maps $f_0^{*}, f^{*}_{1}:C^{*}F(Y) \to C^{*}F(X)$ are chain homotopic. 
\item The cohomology presheaves $X\mapsto H^{i}C^{*}F(X)$ are homotopy invariant.
\item If $f:X\to Y$ is an elementary $\A^{1}$-homotopy equivalence then the induced map of complexes $f_{*}:C_{*}\Z_{tr,G}(X)\to C_{*}\Z_{tr,G}(Y)$ is a chain homotopy equivalence.
 \end{enumerate}
\end{lemma}
\begin{proof}
 The proofs of all these statements are exactly as in the nonequivariant setting. See e.g., \cite[Lecture 2]{MVW}.
\end{proof}

\begin{proposition}
 \begin{enumerate}
  \item Let $U\coprod Y \to X$ be an equivariant distinguished cover. There is a Mayer-Vietoris long exact sequence 
\begin{align*}
\cdots \to H^{n}_{G}(X,\Z(m)) \to  & H^{n}_{G}(U,\Z(m))\oplus H^{n}_{G}(Y,\Z(m)) \to  H^{n}_{G}(U\times_{X}Y,\Z(m)) \\
\to & H^{n+1}_{G}(X,\Z(m)) \to \cdots
\end{align*}
(and similarly for coefficients in $\Z_{G}(V)$).
\item If $G$ satisfies Condition \ref{conmod}  then 
$$
H^{n}_{G}(X\times \A^{1}, \Z(m)) \iso H^{n}_{G}(X,\Z(m))
$$
 (and similarly for coefficients in $\Z_{G}(V)$). 
 \end{enumerate}
\end{proposition}
\begin{proof}
 The first statement follows immediately from the fact that the Bredon motivic cohomology is defined as hypercohomology in the equivariant Nisnevich topology. The cohomology presheaves of $\Z(m)$ are homotopy invariant presheaves with transfers and so the second item follows from Theorem \ref{thm:hi} together with the standard hypercohomology spectral sequence.
\end{proof}

\begin{theorem}\label{thm:C*acy}
 Suppose that  $G$ satisfies Condition \ref{conmod} and $F$ is a homotopy invariant presheaf with equivariant transfers.  
 If $F_{GNis} = 0$ then $(C_{*}F)_{GNis} \wkeq 0$.
\end{theorem}
\begin{proof}
Using the equivariant homotopy invariance result Theorem \ref{thm:hi}, the argument is the same as in \cite[Theorem 13.12]{MVW}.
\end{proof}

\begin{proposition}\label{prop:PV}
 Suppose that $G$ satisfies Condition \ref{conmod}. Let $V$ be a finite dimensional representation. There is a quasi-isomorphism 
$$
C_{*}\big(\Z_{tr,G}(\A(V))/\Z_{tr,G}(\A(V)-0)\big)\xrightarrow{\wkeq}C_{*}\big(\Z_{tr,G}(\P(V\oplus 1)/\Z_{tr,G}(\P(V))\big)
$$ 
of complexes of equivariant Nisnevich sheaves.
\end{proposition}
\begin{proof}

It follows from Example \ref{ex:P1} that the map
$$
\Z_{tr,G}(\A(V))/\Z_{tr,G}(\A(V)-0) \to \Z_{tr,G}(\P(V\oplus 1)/\Z_{tr,G}(\P(V\oplus 1)- \P(1))
$$
becomes an isomorphism after equivariant Nisnevich sheafification. 
The inclusion $\P(V\oplus 1)-\P(1)\subseteq \P(V\oplus 1)$ is equivariantly $\A^{1}$-homotopic to the inclusion $\P(V)\subseteq \P(V\oplus 1)$, the requisite deformation being given by  $([x_{0}:\cdots:x_{n+1}], t)\mapsto [x_{0}:\cdots:x_{n}:tx_{n+1}]$. 

The result now follows from an application of Theorem \ref{thm:C*acy} and Lemma \ref{lem:basicC*}.
\end{proof}
\begin{corollary}\label{cor:V0}
 Under the assumptions above, there is a quasi-isomorphism 
$$
\cone\big(C_{*}\Z_{tr,G}(\A(V)-0)\to \Z \big)
\wkeq C_{*}\big(\Z_{tr,G}(\P(V\oplus 1)/\Z_{tr,G}(\P(V))\big).
$$
\end{corollary}
\begin{proof}
 The map $\A(V)\to \spec(k)$ is an equivariant $\A^{1}$-homotopy equivalence. The result thus follows from the previous proposition together with Theorem \ref{thm:C*acy} and Lemma \ref{lem:basicC*}.
\end{proof}

\subsection{Coefficient systems}\label{sub:coefficients}
Let $\mcal{O}_{G}$ denote the category of finite left $G$-sets with equivariant maps. 

 A \textit{Bredon coefficient system} is an additive functor $M:\mcal{O}_{G}^{op}\to \Ab$. 
Let $\mcal{B}_G$ denote the Burnside category of $G$, its objects are the same as $G\Gamma$ and $\Hom_{\mcal{B}_{G}}(A,B)$ consists of isomorphism classes of diagrams of equivariant maps of finite $G$-sets of the form $A \leftarrow X \rightarrow B$. The composition of $A \leftarrow X \rightarrow B$ and $B \leftarrow X' \rightarrow C$ is given by $A \leftarrow X\times_{B}X' \rightarrow C$. A \textit{Mackey functor} is an additive functor $M:\mcal{B}_G^{op}\to \Ab$. 

The Hecke category $\mcal{H}_{G}$ has the same objects as $\mcal{O}_{G}$ and morphisms are given by $\Hom_{\mcal{H}_{G}}(S,T) = \Hom_{\Z[G]}(\Z[S], \Z[T])$. A \textit{cohomological Mackey functor} is an additive functor $M:\mcal{H}_{G}^{op}\to \Ab$. There is a Hurewicz functor $H:\mcal{B}_{G}\to \mcal{H}_{G}$ given by sending the map $S\xleftarrow{f} X\xrightarrow{g} T$ to the map $\Z[S]\to \Z[T]$ given by $s\mapsto \sum_{x\in f^{-1}(s)}g(x)$. 
Maps between orbits in $\mcal{B}_{G}$ may be interpreted as the maps between orbits in the stable equivariant homotopy category. From this point of view the functor $H$ above is the result of applying the usual Hurewicz functor in stable equivariant homotopy theory.

We have an embedding $\mcal{O}_{G}\subseteq G\Sm/k$ given by $S\mapsto \coprod_{S} \spec(k)$. The composition $\mcal{O}_{G}\subseteq G\Sm/k \subseteq G\Cor_{k}$ factors through a faithful embedding $\mcal{H}_{G}\subseteq G\Cor_{k}$. 
We thus have an embedding of the category of cohomological Mackey functors into the category of presheaves with equivariant transfers. The category of (pre)sheaves with transfers has a tensor structure and by tensoring the complexes $\Z_{G}(n)$ with a cohomological Mackey functor $M$, we obtain  Bredon  motivic cohomology theory with coefficients in $M$.

\begin{example}\label{ex:repsph}
We have an embedding of topological representation spheres into our setting. 
Let $X$ be a based $G$-$CW$ complex write $\underline{C}_{n}(X)$ for the Mackey functor given by $G/H \mapsto C_{n}(X^{H})$, where $C_{*}(W)$ denotes the
reduced cellular chain complex of a $CW$-complex. Let $V$ be a real orthogonal representation and write $S^{V}$ for the one-point compactification of $V$.  Define
$\Z_{top}(V) =  \underline{C}_{*}(S^{V})$.
Note that this definition depends on the choice of equivariant cellular decomposition of $S^{V}$, however two different choices give quasi-isomorphic complexes.
\end{example}

An important case of the topological spheres in the previous example occurs when $G=\Z/2$. Write $S^{\sigma}$ for the topological representation sphere associated to the sign representation. Note that in this case, $\Z_{top}(S^{\sigma}) = \cone(\Z_{tr,G}(\Z/2) \to \Z)$.
As in \cite[Section 8]{MVW}, there is a tensor product $\otimes_{tr}$ on $D^{-}(G\Cor_{k})$ which is induced by $\Z_{tr,G}(X)\otimes_{tr}\Z_{tr,G}(Y) = \Z_{tr}(X\times Y)$. The following will be useful in Section \ref{sec:cancellation}. 
\begin{lemma}\label{lem:invert}
 The complex $\Z_{top}(S^{\sigma})$ is invertible in $(D^{-}(G\Cor_{k}),\otimes_{tr})$.
\end{lemma}
\begin{proof}
 Write $q_{i}:\Z/2\times \Z/2\to \Z/2$ for the projection to the $i$th factor and write  $p:\Z/2\to *$ 
for the projection to a point. The complex $\Z_{top}(S^{\sigma})$ is the complex $0\to \Z_{tr,G}(\Z/2) \xrightarrow{p} \Z \to 0$ (with $\Z = \Z_{tr,G}(*)$ in degree zero). We claim that the inverse $\Z_{top}(S^{-\sigma}))$ is given by $0 \to \Z \xrightarrow{p^{t}} \Z_{tr,G}(\Z/2)\to 0$ (again with $\Z$ in degree zero, and $(-)^{t}$ denotes the transpose). The tensor product $\Z_{top}(S^{\sigma}) \otimes_{tr}^{\mathbb{L}} \Z_{top}(S^{-\sigma})$ is the complex
$$
0 \to \Z_{tr,G}(\Z/2) \xrightarrow{(p,-q_1^t)} \Z\oplus \Z_{tr,G}(\Z/2\times \Z/2) \xrightarrow{p^t\oplus q_2} \Z_{tr,G}(\Z/2) \to 0.
$$
Write $E_*$ for this complex. We have a chain homotopy $s:E_{*}\to E_{*+1}$ between the identity on $E_{*}$ and the composite $E_{*}\to \Z \to E_{*}$ (where $\Z$ is concentrated in degree zero) given by $s_{0} = -\Delta$, $s_{-1} = \Delta$ and $s_{i}= 0$ for $i\neq 0, -1$, where $\Delta:\Z/2\to \Z/2\times\Z/2$ is the diagonal.
\end{proof}

\subsection{Examples}\label{sub:examples}

We record a few simple examples. The first one is straightforward.
\begin{lemma}
Suppose that $G$ satisfies Condition \ref{conmod}.
 Then there is a quasi-isomorphism $\Z_{G}(0) \wkeq \Z$ of complexes of equivariant Nisnevich sheaves, where $\Z$ is the complex consisting of the constant sheaf $\Z$ in degree zero.
\end{lemma}

\begin{proposition}\label{prop:v11}
Let $V$ be a one dimensional representation. Then we have an isomorphism
$$
C_{*}\Z_{tr,G}(\A(V)-0) \wkeq (\mcal{O}^{*})^{G}\oplus \Z,
$$
in the derived category of equivariant Nisnevich sheaves on $G\Sm/k$, 
where $(\mcal{O}^*)^{G}$ is the sheaf of invariant invertible functions viewed as a complex concentrated in degree zero.
\end{proposition}
\begin{remark}
 This is not a decomposition of complexes of sheaves with equivariant transfers except in the case when $V$ is the trivial representation. 
\end{remark}
\begin{proof}

The homology of $C_{*}\Z_{tr,G}(\A(V)\setm 0)(X)$ is identified with equivariant Suslin homology defined in \S \ref{sec:relcart}, 
i.e.,  
$$
H_{i}(C_{*}\Z_{tr,G}(\A(V)\setm 0)(X)) = 
H_{i}^{Sus}(G;X\times (\A(V)\setm 0)/X). 
$$
By Theorem \ref{thm:shcurve} we thus have 
$$
H_{i}(C_{*}\Z_{tr,G}(\A(V)\setm 0 )(X)) = \begin{cases}
                                              \Div^{G}_{rat}(X\times\P(V\oplus 1), X\times \{0, \infty\}) & i = 0 \\
0 & i\neq 0.
                                             \end{cases}
$$
Write $K$ for the kernel of the action of $G$ on $X\times \P(V\oplus 1)$. Then $G/K$ acts faithfully on $X\times \P(V\oplus 1)$, 
$$ 
\Div^{G}_{rat}(X\times\P(V\oplus 1), X\times \{0, \infty\}) = \Div^{G/K}_{rat}(X\times\P(V\oplus 1), X\times \{0, \infty\}),
$$ 
and by Proposition \ref{prop:divinj}, 
$$
\Div^{G/K}_{rat}(X\times\P(V\oplus 1), X\times \{0, \infty\}) = 
\Pic^{G/K}(X\times\P(V\oplus 1), X\times \{0, \infty\}).
$$
Using the exact sequence (\ref{eqn:picseq}) for the relative equivariant Picard group and that for $X$ smooth $\Pic^{G}(X\times\P(V\oplus 1)) = \Pic^{G}(X)\times \Z$, we have
$$
\Pic^{G/K}(X\times\P(V\oplus 1), X\times \{0, \infty\}) = \mcal{O}^{*}(X)^{G/K}\oplus \Z = \mcal{O}^{*}(X)^{G} \oplus \Z,
$$
from which the result follows.
\end{proof}

\begin{proposition}\label{prop:qpG}
 Let $X$ be a smooth, quasi-projective $G$-scheme. Then 
$$
H^{i}_{GNis}(X, (\mcal{O}^*_X)^{G}) = \begin{cases}
                                       \Gamma(X,\mcal{O}^{*})^{G} & i=0 \\
\Pic(X/G) & i = 1 
                                      \end{cases}
$$
\end{proposition}
\begin{proof}
Since $(\mcal{O}^*_X)^{G})$ is the sheaf on $X_{GNis}$ given by $U\mapsto \mcal{O}_{U/G}^{*}$, the proposition follows from Lemma \ref{lem:qpquot}. 
\end{proof}

\begin{corollary}
 Suppose that $G$ satisfies Condition \ref{conmod} and $V$ is a one dimensional representation. Then $\Z_{G}(V)\wkeq (\mcal{O}^*)^G[-1]$. In particular if  $X$ is a smooth, quasi-projective $G$-scheme then 
$$
H^{i}_{G}(X, \Z(V)) = \begin{cases}
                       \Gamma(X, \mcal{O}^*)^{G} & i =1 \\
\Pic(X/G) & i = 2 
                      \end{cases}
$$
\end{corollary}
\begin{proof}
 By Proposition \ref{prop:v11} and Corollary \ref{cor:V0} we have that $\Z_{G}(V)\wkeq (\mcal{O}^*)^G[-1]$. The second statement follows immediately from the previous proposition.
\end{proof}

For  one dimensional representations $V$, $V'$, the chain complexes $\Z_{G}(V)$ and $\Z_{G}(V')$ are quasi-isomorphic as complexes of sheaves in the equivariant Nisnevich topology. The following example makes explicit that for higher dimensional representations, distinct representations give rise to distinct chain complexes. For a representation $V$, we write $\Z_{tr,G}(T^{V}) := \Z_{tr,G}(\P(V\oplus 1))/\Z_{tr,G}(\P(V))$.

\begin{example}
Let $p$ be a prime and $G=\Z/p$ and $k$ a field which admits resolution of singularities. For any representation $V$, we have that 
$$
H^{i}_{G}(k,\Z(V)) =  H^{i-2pn}(C_{*}\Z_{tr,G}(T^{V})(k)) = H^{i-2pn}(z_{equi}(\A(V),0)(\Delta^{\bullet}_{k})^{G}),
$$
 where $z_{equi}(X,0)$ is the presheaf of equidimensional cycles of relative dimension zero. 
We compare the complexes $\Z_{G}(V)[2|V|]= C_{*}\Z_{tr,G}(T^{V})$ for $V=nk[G]$ and $1^{np}$. We have that 
$H^{i}(C_{*}\Z_{tr,G}(1^{np})(k)) = H^{i+2np}_{\mcal{M}}(k,\Z(np))$. On the other hand,  $C_{*}z_{equi}(\A(nk[G]),0)^{G}=\oplus_{j=n}^{np-1}\Z/p(j)[2j] \oplus \Z(np)[2np]$ in $DM(k)$ by \cite[Theorem 5.4]{Nie}. 
Therefore we have that
$$
H^{i}C_{*}\Z_{tr,G}(1^{np})(k) = H^{i+2np}_{\mcal{M}}(k, \Z(np)).
$$
while
$$
H^{i}C_{*}\Z_{tr,G}(T^{nk[G]})(k) =  H^{i+2np}_{\mcal{M}}(k, \Z(np))\oplus (\oplus_{j=n}^{np-1}H^{i+2j}_{\mcal{M}}(k,\Z/p(j)) ).
$$
We see that $C_{*}\Z_{tr,G}(1^{np})$ and $C_{*}\Z_{tr,G}(T^{nk[G]})$ are not quasi-isomorphic in general because there are values of $i$ so that the group $\oplus_{j=n}^{np-1}H^{i+2j}_{\mcal{M}}(k,\Z/p(j))$ is nonzero (e.g. $H^{0}_{\mcal{M}}(k,\Z/p(n)) = \Z/p$ and so the above group is nonzero whenever $i+2j=0$).
\end{example}

%
%

We finish this section by relating our construction to Edidin-Graham's equivariant higher Chow groups \cite{EG:int}. Recall that these are constructed as follows. Consider a pair $(V,U)$ where $V$ is  a faithful representation  and  $U\subseteq \A(V)$ is an open subscheme on which $G$ acts freely.  Then $U/G$ exists as a scheme and it is an algebro-geometric  approximation to $BG$. Such pairs always exist, moreover one can arrange that $\dim V$ and  $\codim_{\A(V)}(\A(V)-U)$ are arbitrarily large. The equivariant higher Chow group of a quasi-projective $G$-scheme $X$ in bidegree $(n,m)$ is defined to be $CH_{G}^{n}(X,m) = CH^{n}(X\times ^{G} U, m)$ for a pair $(V,U)$ with $\A(V)\setm U$ of sufficiently large codimension. We refer to loc.~cit.~for full details. 

\begin{theorem}\label{thm:EGcomp}
 Let $X$ be a smooth quasi-projective $G$-scheme. There is a natural map
$$
H^{i}_{G}(X, \Z(1^{q})) \to CH_{G}^{q}(X, 2q-i).
$$
\end{theorem}
\begin{proof}
We have a natural isomorphism $H^{n}_{\mcal{M}}(X,\Z(q))\iso CH^{k}(X,2q-n)$. 
The complex $\Z(q)$ on $\Sm/k$ is $C_{*}(\Z_{tr}(\P^{q})/\Z_{tr}(\P^{q-1}))[-2n]$. 
If $Y$ has trivial action then $G\Cor_{k}(X,Y) = \Cor_{k}(X/G,Y)$ for a $G$-scheme $X$.
Therefore we have the natural identification  $\Z_{G}(1^{q})(X) = \Z(q)(X/G)$. Using this identification,  
 Lemma \ref{lem:qpquot}, and the projection $X\times U \to X$ we thus have the comparison map 
\begin{align*}
H^{i}_{G}(X, \Z(1^{q})) & \to H^{i}_{G}(X\times U,\Z(1^q)))
\\ & = H^i_{\mcal{M}}(X\times^G U,\Z(q)) = CH^{k}(X\times^{G} U, 2q-i).
\end{align*}
Taking $(V,U)$ such that $\A(V)\setm U$ has sufficiently large codimension yields the result.

\end{proof}

\begin{remark}
 The map of the previous theorem can be seen to be an isomorphism when $X$ has free action. It is not an isomorphism in general. For example, if $X$ has trivial action then $H^{i}_{G}(X,\Z(1^q))$ is isomorphic to ordinary motivic cohomology groups, which in turn is isomorphic to Bloch's higher Chow groups. Under these isomorphisms, the comparison map case just constructed is identified in this case with the comparison map $CH^{q}(X,2q-i)\to CH^{q}_{G}(X, 2q-i)$ between ordinary higher Chow groups and equivariant higher Chow groups, which is not an isomorphism.
\end{remark}

\section{Relative equivariant Cartier divisors}\label{sec:relcart}
In this section we introduce an equivariant version of Suslin homology and relate it to the group of relative equivariant Cartier divisors.

Let $f:X\to S$ be smooth. Recall that $C_{0}(X/S)$ denotes the group of cycles on $X$ which are finite and surjective over a component of $S$.
If $f:X\to S$ is equivariant then we have an equivariant inclusion $C_{0}(X/S)\subseteq \Cor_{k}(S,X)$, induced by  $\langle f,id_{X}\rangle : X\hookrightarrow S\times X$.  
Let $F:G\Cor_{k}^{op}\to \Ab$ be a presheaf with equivariant transfers. 
Define the pairing
\begin{equation}\label{eqn:tran}
Tr:C_{0}(X/S)^{G}\otimes F(X) \to F(S)
\end{equation}
to be the composite
$$
\xymatrix{
C_{0}(X/S)^{G}\otimes F(X) \ar[r]^-{Tr}\ar@{^{(}->}[d]& F(S) \\
G\Cor_{k}(S,X)\otimes F(X) \ar[ur]_-{\textrm{evaluate}} . & 
}
$$

Define  
$$
C_{n}(X/S) = C_{0}(X\times\Delta^{n}/S\times\Delta^{n}).
$$ 
The assignment $n\mapsto C_{n}(X/S)^{G}$ is a simplicial abelian group and hence gives rise to an associated chain complex.

\begin{definition}
The $n$th \textit{equivariant Suslin homology} of $X/S$ is defined to be
$$
H_{n}^{Sus}(G;X/S) = H_{n}C_{\bullet}(X/S)^{G}.
$$
\end{definition}

\begin{lemma}\label{lem:sustr}
Let $F$ be a homotopy invariant presheaf with equivariant transfers.
The map (\ref{eqn:tran}) factors through the zeroth Suslin homology group to yield the pairing
$$
Tr:H_{0}^{Sus}(G;X/S)\otimes F(X) \to F(S).
$$
\end{lemma}
\begin{proof}
Because $F(X\times \A^{1}) = F(X)$ we have the commutative diagram
$$
\xymatrix{
C_{0}(X\times\A^{1}/S\times\A^{1})^{G}\otimes F(X) \ar[r]^-{Tr}\ar[d]_{\partial_{0}-\partial_{1}} & F(S\times\A^{1}) \ar[d]^-{i_{0}-i_{1}} \\
C_{0}(X/S)^{G}\otimes F(X) \ar[r]^-{Tr} & F(S),
}
$$ 
which implies the lemma.
\end{proof}

Our next goal is to compute the equivariant Suslin homology of equivariant relative curves. 
Recall that an equivariant Cartier divisor on $X$ is an element of $\Gamma(X,\mcal{K}_{X}^{*}/\mcal{O}_{X}^{*})^{G}$, see Section \ref{sub:div} for a recollection.

\begin{definition}
\begin{enumerate}
\item Let $X$ be a $G$-scheme and $Y\subseteq X$ an invariant closed subscheme. 
A \textit{relative equivariant Cartier divisor} on $X$ is an equivariant Cartier divisor on $X$ (see Definition \ref{def:ecart}) such that $\supp(D)\cap Y = \emptyset$. Write $\Div^{G}(X,Y)$ for the group of relative equivariant Cartier divisors, where the group operation is induced by that on $\Div^{G}(X)$.
\item A \textit{principal relative equivariant Cartier divisor} is an invariant rational function $f\in \Gamma(X,\mcal{K}^{*})^{G}$ on $X$ such that $f$ is defined and equal to 1 on $Y$.
\item Write $\Div^{G}_{rat}(X,Y)$ for the group of relative equivariant Cartier divisors modulo the principal relative equivariant Cartier divisors.
\end{enumerate}
\end{definition}

Let $i:Y\hookrightarrow X$ be an equivariant closed embedding of $G$-schemes. Define an \'etale sheaf on $X$ by
$$
\G_{X,Y} = \ker(\mcal{O}^{*}_{X} \to i_{*}\mcal{O}^{*}_{Y}).
$$
Since $\mcal{O}^{*}_{X}$ and $\mcal{O}^{*}_{Y}$ are \'etale $G$-sheaves, so is $\G_{X,Y}$.
See Section \ref{sub:Gshv} for a recollection on $G$-sheaves and their cohomology.
\begin{definition}
Define the \textit{relative equivariant Picard group} by
$$
\Pic^{G}(X,Y) = H^{1}_{et}(G;X, \G_{X,Y}).
$$
\end{definition}

From the definition of $\G_{X,Y}$ we have a natural exact sequence
\begin{equation}\label{eqn:picseq}
\Gamma(X, \mcal{O}^{*}_{X})^{G} \to \Gamma(Y,\mcal{O}^{*}_{Y})^{G} \to 
\Pic^{G}(X,Y) \to \Pic^{G}(X) \to \Pic^{G}(Y).
\end{equation}
Theorem \ref{thm:H90} and the above exact sequence imply that $\Pic^{G}(X,Y) = H^{1}_{Zar}(G;X, \G_{X,Y})$. The following lemma is straightforward to verify.
\begin{lemma}
The group $\Pic^{G}(X, Y)$ is isomorphic to the group consisting of isomorphism classes of pairs $(\mcal{L},\phi)$ where $\mcal{L}$ is a $G$-line bundle on $X$ and $\phi$ is an equivariant isomorphism  
$\phi: \mcal{L}|_{Y} \iso  \mcal{O}_{Y} $ of $G$-line bundles on $Y$ and group operation induced by tensor product of $G$-line bundles.
\end{lemma}

\begin{proposition}\label{prop:divinj}
There is a natural injection
$\iota:\Div^{G}_{rat}(X,Y) \hookrightarrow \Pic^{G}(X,Y)$. If in addition $X$ has faithful action  and $Y$ has an invariant affine open neighborhood, then $\iota$ is an isomorphism.  
\end{proposition}
\begin{proof}
 Let $D\in \Div^{G}(X,Y)$. 
Since $Y\cap\supp(D) = \emptyset$, there is a canonical equivariant trivialization 
$s_{D}:\mcal{O}_X(D)|_{Y} \iso \mcal{O}_{Y}$. 
In particular we have a natural homomorphism $\Div^{G}(X,Y) \to \Pic^{G}(X,Y)$. 
If $(\mcal{O}_X(D),s_{D}) = (\mcal{O}_{X}, \id)$ then there is an equivariant isomorphism $\psi:\mcal{O}_{X}\iso \mcal{O}_X(D)$ such that $\psi|_{Y} = (s_{D})^{-1}$. We have an induced isomorphism $\psi:\Gamma(X,\mcal{O}_X)^{G} \iso \Gamma(X,\mcal{O}_X(D))^{G}$ and letting $f = \psi(1) \in \Gamma(X,\mcal{O}_X(D))^{G}\subseteq \Gamma(X,\mcal{K}_X)^{G}$ we have 
that $D= \div(f^{-1})$ and  $D|_Y = 1$ which implies that $\iota$ is injective.

The image of $\iota$ consists of pairs $(\mcal{L},\phi)$ such that $\mcal{L}$ admits an equivariant injection into $\mcal{K}_{X}$ and $\phi$ extends to an equivariant trivialization on an invariant open neighborhood of $Y$. When the action on $X$ is faithful, every $G$-line bundle on $X$ admits an equivariant injection into $\mcal{K}_{X}$ by Proposition \ref{prop:rat}. When $Y$ has an invariant open affine neighborhood, every equivariant trivialization $\phi$ extends to an invariant open neighborhood of $Y$. 
\end{proof}

\begin{lemma}\label{lem:divA1}
 If $X$ is normal and $Y$ is reduced then  
$$
\Pic^{G}(X\times\A^{1},Y\times \A^{1}) \iso \Pic^{G}(X,Y).
$$ 
If $Y$ has an invariant affine open neighborhood then 
$$
\Div^{G}_{rat}(X\times\A^{1},Y\times \A^{1}) \iso \Div_{rat}^{G}(X,Y).
$$
\end{lemma}
\begin{proof}
The first statement follows from the exact sequence (\ref{eqn:picseq}), Lemma \ref{lem:picGA1}, and the five lemma. For the second statement,  observe that if  $K = \ker(G\to \Aut(X))$, then $G/K$ acts faithfully on $X$ and $\Div_{rat}^{G}(X,Y) = \Div_{rat}^{G/K}(X,Y)$. The result  then follows from the first part together with Proposition \ref{prop:divinj}.
\end{proof}

\begin{definition}\label{def:good}
Let $X\to S$ be a smooth relative curve in $G\Sm/k$ (i.e. an equivariant smooth map of relative dimension one). An \textit{equivariant good compactification} of $X$ over $S$ is an equivariant open embedding $X\subseteq \overline{X}$ of $G$-schemes over $S$ where $\overline{X}\to S$ is a proper normal (not necessarily smooth) curve with $G$-action and $X_{\infty}=(\overline{X}-X)_{red}$ has an invariant open affine neighborhood in $\overline{X}$. 
\end{definition}

If $X\to S$ is an equivariant smooth relative curve with equivarian good compactification then  we have an isomorphism $\cyc:\Div(\overline{X}, X_{\infty})\iso C_{0}(X/S)$. Indeed, if $D\in \Div(\overline{X}, X_{\infty})$ 
then $\cyc(D)$ is supported on $X$ and the assumptions above guarantee that it is finite and surjective over a component of $S$. It is straightforward to check this is an equivariant isomorphism (see Lemma \ref{lem:divagr}) and so we immediately conclude the following.

\begin{lemma}\label{lem:C0}
Let $X\to S$ be an equivariant smooth curve  with good equivariant compactification $\overline{X}\to S$. Then $\cyc$ induces a natural isomorphism
$$
\cyc:\Div^{G}(\overline{X}, X_{\infty})\xrightarrow{\iso}C_{0}(X/S)^{G}.
$$
\end{lemma}

With these definitions, Suslin-Voevodsky's fundamental computation of the Suslin homology of relative curves extends to the equivariant setting.

\begin{theorem}\label{thm:shcurve}
 Let $p:X\to S$ be an equivariant smooth quasi-affine curve with equivariant good compactification $\overline{X}\to S$. Then
$$
H^{Sus}_{n}(G;X/S) \cong \begin{cases}
                  \Div^{G}_{rat}(\overline{X},X_{\infty})    & n=0 \\
		    0 & n>0.
                      \end{cases}
 $$
\end{theorem}
\begin{proof}
The argument is similar to \cite[Theorem 3.1]{SV:hom}. Define  
$$
\mcal{M}_{n}(X/S) = \{f\in \Gamma(\overline{X}\times\Delta^{n}, \mcal{K}^{*})\,|\, \textrm{$f$ is defined and equal to $1$ on $S\times\Delta^{n}$}\}.
$$
As shown in  \cite[proof of Theorem 3.1]{SV:hom}, the natural map $\mcal{M}_{n}(X/S)\to C_{n}(X/S)$ is an injection. 
We thus have an exact sequence of complexes
$$
0\to \mcal{M}_{\bullet}(X/S)^{G} \to C_{\bullet}(X/S)^{G} \to \Div^{G}_{rat}(\overline{X}\times\Delta^{\bullet}, Y\times\Delta^{\bullet})\to 0.
$$

By Lemma \ref{lem:divA1} the result follows once we show that $\mcal{M}_{\bullet}(X/S)^{G}$ is acyclic. We work with the normalized chain complex. Let $f\in \mcal{M}_{n}(X/S)^{G}$ and suppose that $\partial_{i}(f) = 1$ for $i=0,\ldots, n$. We need to show that there is $g\in \mcal{M}_{n+1}(X/S)^{G}$ such that $\partial_{i}(g) = 1$ for $i=0,\ldots, n$ and $\partial_{n+1}(g) = f$. 
Following \cite[Theorem 3.1]{SV:hom}, we consider
$$
g_{i} = (t_{i+1} + \cdots + t_{n+1}) + (t_{0} + \cdots + t_{i})s_{i}(f). 
$$
Since $f$ is equivariant it follows that $g_{i}$ is equivariant. The function
$$
g = g_{n}g_{n-1}^{-1}\cdots g_{0}^{(-1)^{n}}
$$
is then equivariant  and by loc.~cit.~it has the required properties.

\end{proof}

We finish with a discussion of the compatibility of the isomorphism in the previous theorem with respect to push-forwards along finite morphisms. 

\begin{lemma}\label{lem:eqnorm}
Let $W\to X$ be an equivariant finite surjection between normal $G$=schemes over $k$. Then the norm map $N:\mcal{K}^{*}(W) \to \mcal{K}^{*}(X)$ is equivariant.
\end{lemma}
\begin{proof}
If $g:Y'\to Y$ is an isomorphism then the norm  $N:\mcal{K}^{*}(Y') \to \mcal{K}^{*}(Y)$ is just the inverse of the induced isomorphism $\tilde{g}:\mcal{K}^{*}(Y) \to \mcal{K}^{*}(Y')$. Thus the lemma follows from functoriality of the norm map.
\end{proof}

Suppose that $f:\overline{X}\to \overline{Y}$ is a finite surjective equivariant map between normal $G$-schemes which restricts to a finite surjective equivariant map  $X_{\infty} \to Y_{\infty}$, where $X_{\infty}\subseteq \overline{X}$, $Y_{\infty}\subseteq Y$ are invariant closed, reduced subschemes.
The norm map induces a map $\Div^{G}(\overline{X}, X_{\infty}) \to \Div^{G}(\overline{Y}, Y_{\infty})$ which factors through rational equivalence to give
$$
f_{*}:\Div^{G}_{rat}(\overline{X}, X_{\infty}) \to \Div^{G}_{rat}(\overline{Y}, Y_{\infty}).
$$

If $X_{\infty}\subseteq \overline{X}$ has an invariant affine neighborhood, every invertible invariant regular function $\alpha$ on $X_{\infty}$ extends to an invariant rational function $\tilde{\alpha}$ on $\overline{X}$.  The difference of two different extensions is a principal relative equivariant divisor and so we have a well-defined homomorphism
$$
\mcal{O}^{*}(X_{\infty})^{G} \to \Div^{G}_{rat}(\overline{X},X_{\infty}).
$$
Additionally we can define $f_{*}:\mcal{O}^{*}(X_{\infty})^{G}\to \mcal{O}^{*}(Y_{\infty})^{G}$ by extending $\alpha$ to $\tilde{\alpha}$ as above and then define  $f_{*}(\alpha) = N(\tilde{\alpha})|_{Y_{\infty}}$. It is easily checked that $N(\tilde{\alpha})|_{Y_{\infty}}$ lies in $\mcal{O}^{*}(Y_{\infty})$ and that this value does not depend on the choice of extension.

\begin{lemma}\label{lem:normcomm}
 Let $(\overline{Y},Y_{\infty})$ and $(\overline{X},X_{\infty})$ be good equivariant compactifications of $Y$ and $X$. Let $f:\overline{Y}\to \overline{X}$ be a finite map which restricts to a map $f:Y\to X$. Then the following diagram commutes
$$
\xymatrix{
\mcal{O}^{*}(X_{\infty})^{G} \ar[r]\ar[d]^{f_{*}} & \Div^{G}_{rat}(\overline{X},X_{\infty})\ar[d]^{f_{*}}\ar[r]^{\cong} & H_{0}^{Sus}(G;X/S) \ar[d]^{f_{*}} \\
\mcal{O}^{*}(Y_{\infty})^{G} \ar[r]& \Div^{G}_{rat}(\overline{Y},Y_{\infty})\ar[r]^{\cong} & H_{0}^{Sus}(G;Y/S),
}
$$
where the left hand and middle vertical maps are induced by the norm map  and the right hand vertical map is  push forward of cycles.
\end{lemma}
\begin{proof}
The commutativity of the left hand square is immediate from the definitions. If $D$ is a Cartier divisor on $\overline{Y}$ then $f_{*}\cyc(D) = \cyc(f_{*}D)$ by \cite[Proposition 21.10.17]{EGAIVpt4}, this implies that the right hand square commutes.
\end{proof}

\section{Equivariant triples}\label{sec:etrip}
In this section we introduce and study an equivariant generalization of Voevodsky's standard triples and establish  equivariant analogues of their basic properties. From now on $k$ is assumed to be perfect. 
Additionally we will usually assume that $G$ satisfies Condition \ref{conmod}, i.e. all irreducible representations are assumed to be one dimensional.

\begin{definition}
An \textit{equivariant standard triple} $(\overline{X} \xrightarrow{\overline{p}} S, X_{\infty}, Z)$ consists of a proper equivariant map $\overline{p}$ of relative dimension one between $G$-schemes and closed invariant subschemes $Z$, $X_{\infty}$ of $\overline{X}$ such that
\begin{enumerate}
 \item $S$ is smooth and $\overline{X}$ is normal
\item $X = \overline{X}-X_{\infty}$ is quasi-affine and smooth over $S$
\item $Z\cap X_{\infty} = \emptyset$
\item $X_{\infty} \cup Z$ has an invariant affine neighborhood in $\overline{X}$.
\end{enumerate}
\end{definition}
Note that $\overline{X}$ is an equivariant  good compactification of both $X$ and $X-Z$.

\begin{remark}
By \cite[Remark 11.6]{MVW}, these conditions imply that $S$ is affine and that $Z$ and $X_{\infty}$ are finite over $S$.
\end{remark}
 
Nonequivariantly any smooth quasi-projective scheme fits into a triple, locally around any finite set of points. Equivariantly this is more delicate. If $f:X\to S$ is an equivariant curve which is smooth at a point $x\in X$ then the induced map 
$\Omega^1_{S/k,f(x)}\otimes_{\mcal{O}_{S,f(x)}}k(x)\to \Omega^1_{X/k, x}\otimes_{\mcal{O}_{X,x}}k(x)$ 
is an injection of $I_{x}$-representations over $k(x)$. However it could happen that 
$\Omega^1_{X/k, x}\otimes k(x)$ has no codimension 1 summand, in which case there can be no such equivariant curve 
$X\to S$ which is smooth at $x$. 
Under the assumption of Condition \ref{conmod} 
we can  construct enough equivariant triples around an orbit in order to establish Theorem \ref{thm:shfinj} below.

Write $T_{x}X:=\Hom_{k(x)}(\Omega^{1}_{X/k,x}\otimes k(x), k(x))$ for the tangent space at $x\in X$. 
\begin{proposition}\label{prop:square}
Let $V$ be a finite dimensional representation, 
$Y\subseteq X\subseteq \A(V)$ equivariant closed embeddings of smooth $G$-schemes, and $x\in Y$ a closed point. Suppose that there are $G$-representations  $W_{2}\subseteq W_{1}$ such that there is an $I_{x}$-equivariant isomorphism $f:T_{x}X\iso (W_{1})_{k(x)}$ which restricts to an $I_{x}$-equivariant isomorphism $T_{x}Y\iso (W_{2})_{k(x)}$. Fix an equivariant splitting $W_1\to W_2$ of the inclusion $W_2\subseteq W_1$. 
Then there is a $G$-equivariant linear projection 
$V\to W_{1}$ such that the composition $X\subseteq\A(V)\to \A(W_{1})$ is \'etale  
and the induced map $Y\to \A(W_{2})$ is also \'etale. 
\end{proposition}
\begin{proof}
Equivariant linear projections $V\to W_1$, which satisfy the above conditions, are parameterized by an open subset $U\subseteq \A(\Hom_{k[G]}(V,W_1))$ of the affine space associated to the $k$-vector space of $G$-equivariant linear maps. More precisely, a point $p\in U$ corresponds to an equivariant $k(p)$-linear map $\A(V)_{k(p)}\to \A(W_{1})_{k(p)}$ such that the induced maps $X_{k(p)}\to \A(W_1)_{k(p)}$ and $Y\to \A(W_2)_{k(p)}$ are \'etale at any point $x'\in Y_{k(p)}$ which lies over $x\in Y$. We need to check that $U$ is nonempty, which implies the result as any nonempty open subset of an affine space has a rational point.

We first treat the case when  $x$ is a rational point.
Consider the diagram
$$
\xymatrix@R=1.5em{
& W_2 & W_1 \ar[l] & \\
& T_{x}Y \ar[u]^{f'}\ar@{^{(}->}[r]\ar[dl] 
& T_{x}X \ar[u]^{f}\ar@{^{(}->}[r]\ar[dl]_{i}
& V , \ar@{-->}@/^1pc/[dll]^-{\rho} \\
\Ind(T_xY)\ar[r] \ar@/^/[uur] & \Ind(T_{x}X) 
\ar@/^/[uur] \ar[urr] & & 
} 
$$
where the inclusion $T_{x}X\subseteq V$ (resp.  $T_{x}Y\subseteq V$) is the one induced by $X\subseteq \A(V)$ (resp. $Y\subseteq \A(Y)$) and $\Ind(M) := k[G]\otimes_{k[I_x]} M$ is the $G$-representation which is induced by the $I_x$-representation $M$.
The $G$-equivariant maps $\Ind(T_{x}X)\to W_1$ and $\Ind(T_{x}Y)\to W_{2}$ are induced respectively by the 
$I_{x}$-equivariant maps $f:T_{x}X\to W_1$ and $f':T_{x}Y\to W_{2}$. Similarly the $I_{x}$-equivariant inclusion
$T_{x}X\subseteq V$ induced by $X\subseteq \A(V)$ induces the $G$-equivariant linear map $\Ind(T_{x}X)\to V$. Choose a $G$-equivariant linear map $\rho:V\to \Ind(T_{x}X)$ so that the composition $T_{x}X \to V \xrightarrow{\rho} \Ind(T_{x}X)$ agrees with the canonical $I_{x}$-equivariant linear map $i$. Then the composition $V\to \Ind(T_{x}X) \to W_1$ has the required properties and thus $U$ is nonempty in this case.

Now suppose that $x\in Y$ is a nonrational closed point. 
Consider the $G$-equivariant embeddings $Y_{k(x)}\subseteq X_{k(x)}\subseteq\A(V)_{k(x)}$ and the points $y_{i}\in Y_{k(x)}$ which lie over $x\in Y$. Consider the open subset $U'\subseteq \A(\Hom_{k(x)[G]}(V_{k(x)},(W_1)_{k(x)}))$ consisting of $p'$ such that the corresponding equivariant linear projections 
$V_{k(p')}\to (W_{1})_{k(p')}$ induces maps $X_{k(p')}\to \A(W_{1})_{k(p')}$  and $Y_{k(p')}\to \A(W_2)_{k(p')}$ which are \'etale at any point $y'\in Y_{k(p')}$ lying over a $y_i\in Y_{k(x)}$. Note that $I_{y_i}=I_{x}$, $T_{y_i}(X_{k(x)}) = T_{x}X$, and 
$T_{y_i}(Y_{k(x)}) = T_{x}Y$ and so the hypothesis of the proposition apply to $y_i\in Y_{k(x)}$. Since these are rational points, the previous paragraph shows that $U'$ is nonempty. Consider the image 
$p\in \A(\Hom_{k[G]}(V,W_1))$ of a $p'\in U'$. The squares
$$
\xymatrix{
X_{k(p')}\ar[r]^-{p'}\ar[d] &\A(W_1)_{k(p')} \ar[d] \\
X_{k(p)}\ar[r]^-{p} & \A(W_1)_{k(p)}
}
\;\;\;\text{and}\;\;\;
\xymatrix{
Y_{k(p')}\ar[r]^-{p'}\ar[d] &\A(W_2)_{k(p')} \ar[d] \\
Y_{k(p)}\ar[r]^-{p} & \A(W_2)_{k(p)}
}
$$
are cartesian and so by faithfully flat descent we conclude that the lower horizontal arrows are \'etale at any point $x'\in Y_{k(p)}$ lying over $x\in Y$. In other words $U'$ maps to $U$ under the projection $\A(\Hom_{k[G]}(V,W_1))_{k(x)}\to \A(\Hom_{k[G]}(V,W_1))$ and so $U$ is nonempty as well.
\end{proof}

\begin{theorem}\label{thm:loctrip}
 Assume that $G$ satisfies Condition \ref{conmod}. Let $X$ be a smooth quasi-projective $G$-scheme over $k$ of pure dimension $d$. Let $Y\subseteq X$ be a smooth invariant closed subscheme containing no component of $X$ and $x\in Y$ a closed point. 
Then there exists an invariant open affine neighborhood  $U$ in $Y$ of $y$ and an equivariant standard triple $(\overline{U}\to S, U_{\infty}, Z)$ such that $(U,U\cap Y) \iso (\overline{U}-U_{\infty}, Z)$.
\end{theorem}
\begin{proof}
First we claim that there are $G$-representations $W_2\subseteq W_1$ defined over $k$, such that there is an isomorphism 
$(W_1)_{k(x)}\iso T_{x}X$ of $I_{x}$-representations over $k(x)$, which restricts to an $I_{x}$-equivariant isomorphism $T_{x}Y\iso (W_{2})_{k(x)}$. Indeed, if $G$ satisfies Condition \ref{conmod} then so does $I_{x}$. From the fact that $k(x)[I_{x}]=k[I_{x}]\otimes_{k} k(x)$ is the sum of the irreducible representations (over $k(x)$), which are one-dimensional, we see that for any representation $M'$ of $I_{x}$ defined over $k(x)$ there is a representation $M$ defined over $k$ such that $M'=M_{k(x)}$. Similarly, Condition \ref{conmod} implies that for every $I_{x}$-representation $N$ there is a $G$-representation  $N'$ such that $N'=N$ as $I_{x}$-representations. 
These observations easily imply the claim. 

Since $X$ is quasi-projective, there is an open invariant affine neighborhood of $x$ and so we may shrink $X$ equivariantly around $x$ and assume that it is affine.
Embed $X$ in some representation $\A(V)$. 
Fixing a choice of equivariant projection $W_1\to W_2$ and applying Proposition \ref{prop:square} we obtain an equivariant linear projection $V\to W_1$ inducing maps $X\to \A(W_1)$ and $Y\to \A(W_2)$ which are \'etale at $x$ (and hence at all points of $G\smash\cdot x$). Let $W\subseteq W_1$ be a codimension-one $G$-representation containing $W_2$ and choose an equivariant projections $W_1\to W$ and $W\to W_{2}$ factoring $W_1\to W_2$. 
The induced equivariant map $p:X\to \A(W)$ is smooth at every point of $G\smash \cdot x$. 
Shrinking $X$ equivariantly around $G\smash\cdot x$, we may assume that 
$p:X\to \A(W)$ is smooth and $Y\to \A(W_2)$ is \'etale. The map $Y\to \A(W)$ is then quasi-finite.

Let $V'\subseteq V$ be a complementary representation to $W$ so that $V = V'\oplus W$. 
Let $\overline{X}\subseteq \P(V'\oplus 1)\times \A(W)$ be the closure of $X$ and write $p:\overline{X}\to \A(W)$ for the induced equivariant map. The fiber of $X\to \A(W)$ over any point of its image is one-dimensional. It follows that $\overline{X}-X$ is finite over $\A(W)$.
Write $\Sigma\subseteq \overline{X}$ for the set of singular points of 
$p:\overline{X}\to \A(W)$. Then $\Sigma\subseteq \overline{X}$ is closed, invariant and is finite over each point of $p(y)$, for any $y\in X$. 
Therefore there is an invariant affine open neighborhood $S$ of $p(G\smash\cdot x)$ in $\A(W)$ over which $\Sigma$ is finite and over which $Y$ has finite fibers. Note that $\Sigma$ and $Y$ are disjoint. 
Define $U$ to be $p^{-1}(S)\cap (X-\Sigma)$. By construction $p:U\to S$ is smooth and equivariant. Define $\overline{U}\subseteq \overline{X}$ to be the preimage of $S$ and set $U_{\infty} = \overline{U}-U$.

It remains to see that we may arrange that $U_\infty\coprod (U\cap Y)$ has an invariant affine neighborhood in $\overline{U}$.
Since $\overline{U}$ is projective over $S$
there is a global section of some very ample line bundle $\mcal{L}$ whose divisor $D$ misses the finite set of points of $U_{\infty}$ and $U\cap Y$ over $G\smash\cdot y$. As $S$ is affine and $\mcal{L}$ is very ample, $\overline{U}-D$ is affine. Intersecting all of the translates of this affine neighborhood, we obtain an invariant open affine neighborhood of all of the points of  $U_{\infty}$ and $U\cap Y$ over $G\smash\cdot x$. Replacing $S$ by a smaller invariant open affine neighborhood  of $p(G\smash\cdot x)$ we may assume that $D$ misses all of $U_{\infty}$ and $U\cap Y$. We thus obtain an invariant affine neighborhood of $U_{\infty}$ and $U\cap Y$.

\end{proof}

Let $(\overline{X}\to S, X_{\infty}, Z)$ be an equivariant triple. Write $\Delta X$ for the equivariant Cartier divisor associated to the diagonal $X\subseteq X\times_{S} X$. 

\begin{definition}
An equivariant standard triple is \textit{equivariantly split} over an invariant open $U\subseteq X$ if 
$\Delta X|_{U\times_{S} Z}$ is an equivariant principal divisor.
\end{definition}

The proof of the following is straightforward.

\begin{lemma}\label{lem:pulllem}
 Let $f:S'\to S$ be an equivariant map between smooth affine $G$-schemes over $k$ and $T=(\overline{X}\to S, X_{\infty}, Z)$ an equivariant triple over $S$. Then $f^{*}T = (\overline{X}\times_{S} S'\to S', X_{\infty}\times_{S} S', Z\times_{S} S')$ is an equivariant triple over $S'$. If $T$ is equivariantly split over $U$ then $f^{*}T$ is equivariantly split over $U\times_{S}S'$.
\end{lemma}

In the equivariant case, the question of a triple being locally split is more delicate than its nonequivariant analog. Nonequivariantly, all divisors on $X\times X$ are locally principal when $X$ is smooth. The nonequivariant argument requires more work as an equivariant Weil divisor (equivalently by Lemma \ref{lem:divagr}, an equivariant Cartier divisor) on a smooth $G$-scheme might not be locally equivariantly principal. This can be seen for example from Proposition \ref{prop:rat} together with the fact that $\Pic^{G}(S)$ can be nonzero for local rings $S$.

If $\pi:A\to A/G$ is a quotient and $B\subseteq A$ is an invariant closed subscheme then since $|G|$ is coprime to $\ch(k)$ the canonical map $B/G\to \pi(B)$ is an isomorphism. In particular we have a Cartesian square
$$
\xymatrix{
X \ar@{^{(}->}[r]^-{\Delta}\ar[d] & X\times_{S}X \ar[d]^{\pi} \\ 
X/G \ar@{^{(}->}[r] & (X \times_{S} X)/G,
}
$$
whenever the right hand vertical quotient exists.

\begin{proposition}
 Let $J$ be a smooth $G$-scheme which is finite over $k$ and let $C\to J$ be a smooth equivariant 
curve. Then the Weyl divisor $(\Delta C)/G \hookrightarrow (C\times_{J} C)/G$ is locally principal.
\end{proposition}
\begin{proof}
Consider the coherent sheaf $\mcal{O}(\Delta C/G)$ on  
$(C\times_{J} C)/G$ associated to the Weyl divisor $(\Delta C)/G$.
The condition that the divisor $(\Delta C)/G$ is locally principal is equivalent to the condition that the coherent sheaf  $\mcal{O}(\Delta C/G)$ is locally free. 
Let $\overline{k}$ be an algebraic closure of $k$. The sheaf $\mcal{O}(\Delta C/G)$ is locally free if it is so after base change to $\overline{k}$. The base change of $\mcal{O}(\Delta C/G)$ to $\overline{k}$ is the sheaf associated to $(\Delta C/G)_{\overline{k}}$. Since $(\Delta C/G)_{\overline{k}} = (\Delta C_{\overline{k}})/G$ it is enough to consider the case when $k$ is algebraically closed.
We may also assume that $J$ is equivariantly irreducible. Since $k$ is algebraically closed, $J= G/H$ for some subgroup $H$ and so 
$C = G\times^{H}C'$ for some smooth  $H$-curve $C'\to \spec(k)$. Since $(C\times_{J} C)/G = (C'\times C')/H$ we may replace $C$ by $C'$ and $G$ by $H$. In other words, we may assume that $J=\spec(k)$ and it suffices to show that for a smooth $G$-curve over the algebraically closed field $k$, the coherent sheaf associated to $\Delta C/G\hookrightarrow (C\times C)/G$ is locally free.

Let $c\in C$ be a closed point. Write $\mcal{O}(\Delta C/G)_{[c]}$ for the restriction of $\mcal{O}(\Delta C/G)$ to $(C\times Gc)/G$. The sheaf $\mcal{O}(\Delta C/G)_{[c]}$ is the coherent sheaf associated to the divisor 
$(\Delta (G\smash\cdot c))/G\hookrightarrow (C\times G  c)/G$. 
We have that $G\smash\cdot c \iso G/I_{c}$ and so the divisor $(\Delta G c)/G\hookrightarrow (C\times G  c)/G$ is identified with $[c] \hookrightarrow  C/I_{c}$ under the equivariant isomorphisms
$(C\times G\sdot c)/ G\iso (C\times G/I_{c})/G\iso (C\times^{I_{c}}G)/G \iso C/I_{c}$, the second isomorphism arising from $(c,[g])\mapsto (g^{-1}c, g)$. Normality is preserved under taking quotients and so $C/I_{c}$ is a normal curve and therefore it is also  smooth and so $\mcal{O}(\Delta C/G)_{[c]}$ is locally free of rank one. Every closed point of $(C\times C)/G$ is in some $(C\times G \smash\cdot c)/G$ and so $\rank_{x}[\mcal{O}(\Delta C/G)]=1$ (where $\rank_{x}\mcal{F} = \dim_{k(x)}\mcal{F}_{x}\otimes k(x)$)
for every closed point $x\in (C\times C)/G$. But the collection of points where the rank of a coherent sheaf takes on a fixed value is constructible and so $\rank_{x}[\mcal{O}(\Delta C/G)]=1$ for every $x\in (C\times C)/G$. A coherent sheaf $\mcal{F}$ on a reduced scheme $X$ is locally free exactly when the function $x\mapsto \rank_{x}\mcal{F}$ on $X$ is  locally constant. We conclude that $\mcal{O}(\Delta C/G)$ is locally free of rank one.
\end{proof}

\begin{corollary}\label{cor:curve}
 Let $X\to S$ be a smooth equivariant curve, with $X$ and $S$ quasi-projective $G$-schemes. Then the equivariant Cartier divisor $\Delta X \hookrightarrow X\times_{S}X$ is equivariantly locally principal.
\end{corollary}
\begin{proof}
 We also write $\Delta X$ for the associated equivariant Weyl divisor on $X\times_{S}X$. We have the Cartesian square of normal schemes
$$
\xymatrix{
\Delta X \ar[r]\ar[d] & X\times_{S}X \ar[d]^{\pi} \\ 
\Delta X/G \ar[r] & (X \times_{S} X)/G.
}
$$
The equivariant Weyl divisor $\Delta X \hookrightarrow X\times_{S} X$ is equivariantly locally principal if the Weyl divisor $(\Delta X)/G \hookrightarrow (X\times_{S} X)/G$ is locally principal. This Weyl divisor is locally principal exactly when the associated coherent sheaf $\mcal{O}((\Delta X)/G)$ is locally free. Consider the map $p:(X\times_{S} X)/G \to S/G$.  The fibers over a point $[s]\in S/G$ are $p^{-1}([s])= (X_{Gs}\times_{Gs} X_{Gs})/G$. The restriction
 $\mcal{O}((\Delta X)/G)_{[s]}$ of $\mcal{O}((\Delta X)/G)$ to the fiber $p^{-1}([s])$ is the coherent sheaf associated to $(\Delta X_{Gs})/G\hookrightarrow (X_{Gs}\times_{Gs} X_{Gs})/G$. By the previous proposition $\mcal{O}((\Delta X )/G)_{[s]}$ is locally free of rank one for all closed points $[s]\in S/G$. Every closed point of $\supp(\Delta X/G)$ is in some $p^{-1}([s])$. 
The sheaf $\mcal{O}((\Delta X)/G)$ is isomorphic to the trivial line bundle at all points not in $\supp(\Delta X/G)$. We conclude that the coherent sheaf  $\mcal{O}((\Delta X)/G)$ has rank one at all closed points and hence has rank one at all points. It follows that it is locally free of rank one.
\end{proof}

The following is an important example.

\begin{lemma}\label{lem:splittrip4}
Let $J$ be an equivariantly irreducible, smooth zero-dimensional $G$-scheme over $k$ and $W$ a $G$-representation and $L:=J\times \A(W)$.  Let $X_{\infty}$ and $Z$ be disjoint, invariant nonempty finite subsets of $\P(L\oplus 1)$. Then 
$$
T:=(\P(L\oplus 1)\to J, X_{\infty} , Z)
$$ 
is an equivariant standard triple which is equivariantly split over any invariant open $U\subseteq L$.
\end{lemma}
\begin{proof}
That $T$ is an equivariant standard triple is clear.  

Let $x\in J$ be a point. Then $J=G\times^{G_{x}}\{x\}$,
$\P(L\oplus 1) \iso G\times ^{G_{x}}(\{x\}\times \P(W\oplus 1))$, 
$U= G\times^{G_{x}} (\{x\}\times U')$ for a $G_{x}$-invariant open $U'\subseteq \A(W)$, and $X_{\infty} = G\times^{G_{x}}(\{x\}\times
X_{\infty}')$ and $Z = G\times^{G_{x}}(\{x\}\times Z')$ for $G_x$-invariant disjoint subsets $X_{\infty}'$ and $Z$ of $\P(W\oplus 1)$. It is thus enough to show that the $G_{x}$-equivariant triple $(\P(W\oplus 1), X_{\infty}', Z')$ is split over any $G_{x}$-invariant open $U'\subseteq  \A(W)$.

We show that $\Delta\A(W)$ is equivariantly principal on $\A(W)\times\A(W)$. To show this it suffices to show that 
$\mcal{O}(\Delta\A(W))$ is the trivial $G_{x}$-line bundle. 
By Corollary \ref{cor:curve} the equivariant Cartier divisor $\Delta \A(W) \subseteq \A(W)\times \A(W)$
 is equivariantly locally principal which implies in turn that  $\Delta \A(W)/G_{x} \subseteq (\A(W)\times \A(W))$ is locally prinicipal. Therefore
$\mcal{O}(\Delta \A(W)/G_x)$, the coherent sheaf on $(\A(W)\times \A(W))/G_x$ associated to the Weil divisor $(\Delta \A(W))/G_x$, is a line bundle. 
By \cite[Theorem 2.4]{Kang:pic}, $\Pic((\A(W)\times \A(W))/G_x)=0$ and so the $G_x$-line bundle $\mcal{O}(\Delta \A(W))= \pi^*\mcal{O}(\Delta \A(W)/G_x)$ is trivial as needed.
\end{proof}

\begin{theorem}\label{thm:tripsplit}
 Let $(\overline{X}\to S, X_{\infty}, Z)$ be an equivariant standard triple. Then any finite set of points in $X$ has an invariant open neighborhood $U$ over which this triple splits.  
\end{theorem}
\begin{proof}
Let $\mcal{P}\subseteq X$ be a finite set of points in $X$. Replacing $\mcal{P}$ by $G\sdot\mcal{P}$ we may assume that $\mcal{P}$ is invariant. The equivariant map $\pi:X\times_{S}Z\to X$ is finite and so $\pi^{-1}\mcal{P}\subseteq X\times_{S}Z$ is also an invariant finite set of points. It follows from Corollary \ref{cor:curve} that
$D/G$ is locally principal on $(X\times_{S} Z)/G$. Let $W\subseteq (X\times_{S} Z/G)$ be a neighborhood of $\mcal{P}/G$ on which $D/G$ is principal and let $V\subseteq X\times_{S} Z$ be its preimage. There is some equivariant neighborhood $U$ of $\mcal{P}$ such that $U\times_{S} Z\subseteq V$. The equivariant triple is split over this $U$.
\end{proof}

Two equivariant finite correspondences $\lambda_{0},\lambda_{1}\in G\Cor_{k}(X,Y)$ are said to be \textit{equivariantly $\A^{1}$-homotopic} provided there is an 
$H\in G\Cor_{k}(X\times\A^{1}, Y)$ such that $H|_{X\times\{i\}} = \lambda_{i}$, $i=0,1$.

\begin{proposition}\label{prop:lambda}
 Let $(\overline{X} \xrightarrow{\overline{p}} S, X_{\infty}, Z)$ be an equivariant standard triple which is split over an open affine $U\subseteq X$. Then there is an equivariant finite correspondence 
$$
\lambda:U\to X\setm Z
$$ 
such that $\lambda$ composed with $j:X\setm Z \subseteq X$ is equivariantly $\A^{1}$-homotopic to the inclusion $i:U\subseteq X$. In particular for any homotopy invariant presheaf with equivariant transfers $F$, we have the commutative diagram
$$
\xymatrix{
F(X) \ar[r]^{j^{*}}\ar[d]_{i^{*}} & F(X\setm Z) \ar[dl]^{\lambda^{*}} \\
F(U) . & 
}
$$
\end{proposition}
\begin{proof}
 We write $\overline{X}_{U} = U\times_{S}\overline{X}$. Pulling back to $U$ gives the equivariant triple $(\overline{p}':\overline{X}_{U}\to U, 
(X_{\infty})_{U}, Z_{U})$. The diagonal $\Delta:U\to X_{U}$ is an equivariant 
section of $\overline{p}'$, so is an element of $C_{0}(X_{U}/U)^{G}$. By Theorem \ref{thm:shcurve} it thus determines the class
 $\Delta U\in \Div^{G}_{rat}(\overline{X}_{U}, (X_{\infty})_{U})$. 
By assumption, $\Delta U$ restricted to $Z_{U}$ is equivariantly principal, say $\Delta U|_{Z_{U}} = \div(r_{U})$, where $r_{U}$ is an invariant regular function. Since $Z_{U}\coprod (X_{\infty})_{U}$ has an invariant affine neighborhood in $\overline{X}_{U}$, we can use the Chinese remainder theorem to find an invariant rational function $\phi$ on $\overline{X}_{U}$ which is defined in an invariant neighborhood of $Z_{U}\coprod (X_{\infty})_{U}$ and is equal to $1$ on $(X_{\infty})_{U}$ and equal to $r_{U}$ on $Z_{U}$. Note that $\div(\phi)$ is zero in $\Div^{G}_{rat}(\overline{X}_{U}, (X_{\infty})_{U})$.
We lift the class $\Delta U$ to a class $[\lambda'] \in \Div^{G}_{rat}(\overline{X}_{U}, (X_{\infty})_{U}\coprod Z_{U})$ by setting $[\lambda'] = \Delta U - \div(\phi)$.

 Let $F$ be any homotopy invariant presheaf with transfers. The diagram
$$
\xymatrix{
F(X_{U}) \ar[r] \ar[d]_-{Tr([\Delta])} & F((X\setm Z)_{U})\ar[dl]^-{Tr([\lambda])} \\
F(U) & 
}
$$
is commutative, where the vertical and diagonal maps are those obtained from Lemma \ref{lem:sustr}.
Let $\lambda' \in C_{0}((X\setm Z)_{U}/U)^{G}\subseteq G\Cor_{k}(U,(X\setm Z)_{U})$ be any representative of $[\lambda']$ and $\lambda:U \to X\setm Z$ be the composition of $\lambda'$ together with the projection to $X\setm Z$. It is easily verified that $j\lambda$ and $i$ are equivariantly $\A^{1}$-homotopic.
\end{proof}

\begin{corollary}\label{cor:pttriag}
 Assume that $G$ satisfies Condition \ref{conmod}. 
 Let $F$ be a homotopy invariant presheaf with equivariant transfers, $Z\subseteq X$  a closed embedding of smooth quasi-projective $G$-schemes over $k$, and $x\in X$ a  closed point. 
 Then there exists an open invariant neighborhood $U$ of $x$, and  a map $\phi:F(X-Z)\to F(U)$ such that the following triangle commutes,
$$
\xymatrix{
F(X) \ar[d]\ar[dr] & \\
F(X-Z) \ar[r]^{\phi} & F(U).
}
$$
\end{corollary}
\begin{proof}
 If $x\notin Z$ there is nothing to prove. If $x\in Z$ then by Theorem \ref{thm:loctrip} there is an invariant open neighborhood $X'$ of $x$ and an equivariant triple $(\overline{X'}, X'_{\infty}, Z')$ such that $(X', X'-Z) = (\overline{X'}- X'_{\infty}, Z')$. By Theorem \ref{thm:tripsplit} there is an invariant open neighborhood $U$ of $x$ such that this triple splits over $U$. Applying Proposition \ref{prop:lambda} to $U$ yields the corollary. 
\end{proof}

\begin{theorem}\label{thm:shfinj}
 Assume that $G$ satisfies Condition \ref{conmod}. Let $F$ be a homotopy invariant presheaf with equivariant 
 transfers, $S$ a smooth semilocal affine $G$-scheme over $k$ with a single closed orbit and $S_0\subseteq S$ a dense invariant open subscheme. Then the restriction map 
 $F(S)\to F(S_0)$ is injective.
\end{theorem}
\begin{proof}
Write the $G$-scheme $S$ as the intersection $\cap X_{i}$ and $S_0= \cap V_i$ where $X_i$ are invariant open neighborhoods of a point $x$ of a smooth affine $G$-scheme $X$, $V\subseteq X$ an invariant open, and $V_i = V\cap X_i$.

Write $Z= (X-V)_{red}$. First observe that we may assume that $Z$ is smooth. Indeed, since $k$ is perfect, there is a filtration $\emptyset=Z(n+1)\subseteq Z(n)\subseteq\cdots\subseteq Z(1)\subseteq Z(0)=(X-V)_{red}$ by closed invariant subschemes such that $Z(r)-Z(r-1)$ is smooth (take $Z(r)\subseteq Z(r+1)$  to be the set of singular points). Write $Z(r)_{i}= X_i\cap Z(r)$. Each $X_i-Z(r-1)_i\subseteq X_i-Z(r)_i$ is the complement of an invariant smooth closed  subscheme.  If the morphism 
$F(\cap (X_i-Z(r)_i))\to F(\cap (X_i-Z(r-1)_i))$ is injective for all $r$, then  $F(\cap X_i)\to F(\cap V_i)$ is injective. Thus we may assume that $Z$ is smooth. Consequently $Z_{i}:= X_{i}-V_{i}$ is also smooth.

Now the $U_{i}$ given by Corollary \ref{cor:pttriag} is contained in some $X_{j}$ and so the kernel of $F(X_i)\to F(V_i)$ vanishes in $F(X_{j})$. Thus 
the map $F(S)\to F(S_{0})$ is injective.
\end{proof}

Recall that a $G$-scheme $W$ is called \textit{equivariantly irreducible} if there is an irreducible component $W_{0}$ of $W$ such that $G\smash\cdot W_{0} = W$. The underlying scheme of an essentially smooth, zero dimensional $G$-scheme $J$ over $k$ is a disjoint union of the Zariski spectra of finitely generated field extenstions of $k$.

\begin{corollary}\label{cor:vanfield}
  Assume that $G$ satisfies Condition \ref{conmod}. Suppose that $F$ is a homotopy invariant presheaf with equivariant transfers and  that $F(J) = 0$ for any essentially smooth, zero dimensional $G$-scheme $J$ over $k$. Then $F_{GNis} = 0$.
\end{corollary}

\begin{definition}
An \textit{equivariant covering morphism}  $f:T_{Y}\to T_{X}$, of two  equivariant standard triples
 $T_{Y} = (\overline{Y}\to S, Y_{\infty}, Z_{Y})$ and $T_{X} = (\overline{X}\to S, X_{\infty}, Z_{X})$, is an equivariant finite map $f:\overline{Y}\to \overline{X}$ such that 
\begin{enumerate}
 \item $f(Y)\subseteq X$,
\item $f|_{Y}:Y\to X$ is \'etale,
\item $f$ induces and isomorphism $Z_{Y}\xrightarrow{\cong} Z_{X}$, and $Z_{Y} = f^{-1}Z_{X}\cap Y$.
\end{enumerate}
\end{definition}

Write $Q(X,Y,A)$ for the equivariant distinguished square
$$
\xymatrix{
 B \ar[r]\ar[d] & Y \ar[d]^{f} \\
 A \ar[r]^{i} & X,
 }
$$
where $i$ is an equivariant open embedding and $f:Y\to X$ is an equivariant \'etale morphism. 

\begin{definition}
Let $f:T_{Y}\to T_{X}$ be an equivariant covering morphism of equivariant standard triples as above. The \textit{associated equivariant distinguished square} to this morphism is $Q= Q(X,Y, X-Z_{X})$ and we say that the square $Q$ \textit{comes from} this covering morphism.
\end{definition}

The following is an important class of examples.
\begin{example}\label{ex:cov}
 Suppose that $X$ is affine, has an equivariant good compactification $\overline{X}$ over some smooth $S$ (see Definition \ref{def:good}), and $X=U\cup V$ is an open cover by invariant open subschemes such that $X-(U\cap V)$ has an invariant open affine neighborhood. Then 
$$
\xymatrix{
U\cap V \ar[r]\ar[d] & U \ar[d] \\
V\ar[r] & X 
}
$$
comes from the morphism of triples $(\overline{X}, \overline{X} - U, X-V) \to (\overline{X}, X_{\infty}, X- V)$ defined by the identity on $\overline{X}$.
\end{example}

The proof of the following theorem (and the lemmas below on which it depends) are similar to the arguments in the nonequivariant case. We include complete details for the reader's convenience.
 
\begin{theorem}\label{thm:complex}
Let $X$ be a smooth equivariantly irreducible $G$-scheme over $k$. Let $Q'=Q(X', Y', A')$ and $Q= Q(X,Y, A)$  be  equivariant distinguished squares such that $Q'$ is the restriction of $Q$ along an invariant open subscheme $X'\subseteq X$. Write $j:Q'\hookrightarrow Q$ for the inclusion.
Assume that $X'$ and $Y'$ are affine and  that 
$Q$ comes from an equivariant covering map 
$$
T_{Y} =(\overline{Y}, Y_{\infty}, Z_{Y}) \to T_{X}=(\overline{X}, X_{\infty}, Z_{X})
$$ 
of equivariant standard triples and that $T_{X}$ splits over $X'$.

 Let $F:G\Cor_{k}^{op}\to \Ab$ be a homotopy invariant presheaf with equivariant transfers. Then the map of complexes
$$
\xymatrix{
0 \ar[r] & F(X) \ar[d]^-{j_{X}}\ar[r]^-{(i,f)} & F(A)\oplus F(Y) \ar[r]^-{(-f,i)}\ar[d]^-{\Big(\stackrel{j_{A}}{j_{Y}}\Big)}
& F(B) \ar[r]\ar[d]^-{j_{B}} & 0 \\
0 \ar[r] & F(X') \ar[r]^-{(i',f')} & F(A')\oplus F(Y') \ar[r]^-{(-f',i')} & F(B') \ar[r] & 0
}
$$
is chain homotopic to zero.

In particular if $Q'=Q$ then the Mayer-Vietoris sequence
$$
0 \to F(X)\to F(A)\oplus F(Y) \to F(B) \to 0
$$
is split-exact.
\end{theorem}
\begin{proof}
By Lemmas \ref{lem:lambda2} and \ref{lem:lem3}  we have maps $s_{1}= (\lambda_{A},0):F(A)\oplus F(Y) \to F(X')$ and $s_{2}=(\psi,\lambda_{B}):F(B)\to F(A')\oplus F(Y')$. For these maps to form a chain homotopy from $j$ to zero we need that $sd+ds=j$. This boils down to six equations. Three come from the commutativity of the trapezoid in Lemma \ref{lem:lambda2}. The remaining three which involve $\psi$ are $\psi i \wkeq 0$, $j_{A}\wkeq i'\lambda_{A} - \psi f$ and $j_{B} \wkeq i'\lambda_{B} - f'\psi$. These follow from Lemma \ref{lem:lem3}. 
\end{proof}

\begin{lemma}\label{lem:covsplit}
 Let $f:T_{Y}\to T_{X}$ be an equivariant covering morphism of equivariant standard triples. If $T_{X}$ is equivariantly split over $V$ then $T_{Y}$ is equivariantly split over $f^{-1}(V)\cap Y$.
\end{lemma}
\begin{proof}
 By assumption the equivariant Cartier divisor $\Delta X |_{V\times_{S}Z_{X}}$ is an equivariant principal divisor, say  $\Delta X|_{V\times_{S} Z_{X}}= \div(\phi)$. Then $(f\times f)^{*}(\Delta X) = \Delta Y + Q$, where the support of $Q$ is disjoint from that of $\Delta Y$. Since $Z_{Y}\iso Z_{X}$, $\supp(Q)$ is also disjoint from $Y\times_{S} Z_{Y}$ and therefore $Q|_{Y\times_{S} Z_{Y}} = 0$. Since $(\Delta Y + Q)|_{(f^{-1}V\cap Y)\times_{S} Z_{Y}}= \div(\phi \circ (f\times f))$, it follows that 
 $\Delta Y |_{(f^{-1}V\cap Y)\times_{S} Z_{Y}} = \div(\phi \circ (f\times f))$ as well.
\end{proof}

\begin{lemma}\label{lem:lambda2}
 Let $j:Q'\hookrightarrow Q$ be as above. Then there are finite equivariant correspondences  $\lambda_{A}:X'\to A$ and $\lambda_{B}:Y'\to B$ such that the following diagram in $G\Cor_{k}$ is commutative up to equivariant $\A^{1}$-homotopy,
$$
\xymatrix@+1pc{
 & Y' \ar@{_{(}->}[dl]_{j_{Y}}\ar[d]^{\lambda_{B}}\ar[r]^{f'} & X' \ar[d]_{\lambda_{A}}\ar@{^{(}->}[dr]^{j_{X}} &  \\
Y & B \ar@{_{(}->}[l]^{i}\ar[r]^{f} & A\ar@{^{(}->}[r]^{i} & X.
}
$$

\end{lemma}
\begin{proof}
 The equivariant triple $T_{X}$ is split over $X'$. 
By Lemma \ref{lem:covsplit}, $T_{Y}$ splits over $Y'$. Proposition \ref{prop:lambda} gives the existence of $\lambda_{A}$ and $\lambda_{B}$ making the triangles commute up to $\A^{1}$-homotopy. The square is easily seen to commute up to $\A^{1}$-homotopy by the construction used in the proof of Proposition \ref{prop:lambda}.
\end{proof}

\begin{lemma}\label{lem:lem3}
 Let $j:Q'\hookrightarrow Q$ be as above. 
There is an equivariant correspondence $\psi\in G\Cor_{k}(A',B)$ such that the square 
$$
\xymatrix@+1.5pc{
B'\ar[r]^-{\lambda_{B}\circ i' - j_{B}} \ar[d]_{f'} & B \ar[d]^-{f} \\
A' \ar@{-->}[ur]_-{\psi}\ar[r]_-{\lambda_{A}\circ i' -j_{A}} & A
}
$$
is homotopy commutative in $G\Cor_{k}$, where $\lambda_{A}\in G\Cor_{k}(X',A)$, $\lambda_{B}\in G\Cor_{k}(Y',B)$ are the equivariant correspondences from Lemma \ref{lem:lambda2}. Moreover the composite $i\psi:A'\to Y$ is equivariantly $\A^{1}$-homotopic to zero.
\end{lemma}
\begin{proof}
First we define the equivariant correspondence $\psi\in G\Cor_{k}(A',B)$. Write $\Delta X'$ for the equivariant Cartier divisor on $X'\times_{S} \overline{X}$ corresponding to the graph of $X'\hookrightarrow \overline {X}$. Similarly, write $\Delta Y'$ for the equivariant Cartier divisor corresponding to the graph of $Y'\hookrightarrow \overline{Y}$. Write $\mcal{M}$ for the pullback of $\Delta X'$ to $X'\times_{S} \overline{Y}$.

The support of $\Delta X'$ is disjoint from $A'\times_{S} Z_{X}$ so $\Delta X'|_{A'\times_{S} Z_{X}} = 0$. Similarly 
$\mcal{M}|_{A'\times_{S} Z_{X}} = 0$ and 
$\Delta Y'|_{B'\times_{S} Z_{Y}} = 0$.

By assumption the equivariant Cartier divisor 
$\Delta X' |_{X'\times_{S} Z_{X}}$ is equivariantly principal. By Lemma \ref{lem:covsplit}, $\Delta Y' |_{Y'\times_{S} Z_{Y}}$ is equivariantly principal as well. 
Write $\Delta X' |_{X'\times_{S} Z_{X}} = \div(r_{X})$, $\mcal{M}|_{X'\times_{S} Z_{Y}} = \div(r_{M})$, and $\Delta Y'|_{Y'\times_{S} Z_{Y}} = \div(r_{Y})$, where $r_{X}$, $r_{M}$, and $r_{Y}$ are invariant regular functions. Furthermore we have that $r_{X}$ is invertible on $A'\times_{S} Z_{X}$. Similarly $r_{M}$ is invertible on $A'\times_{S}Z_{Y}$ and $r_{Y}$ is invertible on $Y'\times_{S} Z_{Y}$. Under the isomorphism $Z_{Y}\iso Z_{X}$, $r_{X}$ becomes identified with $r_{M}$.

Let $U$ be an invariant affine neighborhood of $Y_{\infty}\coprod Z_{Y}$ in $\overline{Y}$. Then $X'\times_{S} U$ is an invariant affine neighborhood of $X'\times_{S}(Y_{\infty}\coprod Z_{Y})$. Since points of $(X'\times_{S} U)/G$ are orbits, $X'\times_{S} Y_{\infty}$ and $X'\times_{S} Z_{Y}$ remain disjoint in $(X'\times_{S} U)/G$. 
Since $\ch(k)$ doesn't divide $|G|$, we have $(X'\times_{S}Y_{\infty})/G = \pi(X'\times_{S} Y_{\infty})$ and similarly for $X'\times_{S}Z_{Y}$ where $\pi:X'\times_{S} U \to (X'\times_{S} U)/G$ is the quotient map. The invariant regular functions $1$ and $r_{M}$ on $X'\times_{S}Z_{Y}$ and $X'\times_{S} Y_{\infty}$ define invariant regular functions on their quotients. Since $(X'\times_{S}U)/G$ is affine, we may apply the Chinese remainder theorem, see e.g. \cite[Proposition B.1]{AG1}) to obtain a regular function $\overline{h}$ on $(X'\times_{S}U)/G$ that equals $\overline{r_{M}}$ on $(X'\times_{S} Z_{Y})/G$ and $1$ on $(X'\times_{S} Y_{\infty})/G$. We thus have an invariant regular function $h$ on $X'\times_{S} U$ which equals $r_{M}$ on $X'\times_{S} Z_{Y}$ and $1$ on $X'\times_{S} Y_{\infty}$.

View $h$ as an invariant rational function on $A'\times_{S} \overline{Y}$. The support of its associated divisor $\div(h)$ is disjoint from $A'\times_{S}(Z_{Y} \coprod Y_{\infty})$ and so is an element of $\Div^{G}(A'\times_{S} \overline{Y}, A'\times_{S}(Z_{X}\coprod Y_{\infty})) = C_{0}(A'\times_{S} B/A')^{G}$. Since $C_{0}(A'\times_{S}B/A')^{G} \subseteq G\Cor_{k}(A',B)$, the divisor $-\div(h)$ determines the equivariant correspondence $\psi:A'\to B$. It remains to verify its properties.

First $i\psi \in G\Cor_{k}(A',Y)$ corresponds to $-\div(h)$ in $\Div^{G}(A'\times_{S}\overline{Y} A'\times_{S} Y_{\infty})$. But since $h|_{A'\times_{S} Y_{\infty}} = 1$, $-\div(h)$ is a principal relative equivariant Cartier divisor and so represents $0$ in $H_{0}^{Sus}(G;A'\times_{S} Y/A')$. Thus $i\psi$ is equivariantly $\A^{1}$-homotopic to zero.

It remains to see that the diagram of the lemma is homotopy commutative. By the construction of $\lambda_{A}$ and $\lambda_{B}$ the composition $\lambda_{A}\circ i'\in G\Cor_{k}(A',\,A)$ and $\lambda_{B}\circ i' \in G\Cor_{k}(B',\,B)$ 
are represented by the classes $\Delta A' - \div(\phi_{X})$ and $\Delta B'- \div(\phi_{Y})$ in  
 $\Div^{G}_{rat}(A'\times_{S} \overline{X}, A'\times_{S} (X_{\infty}\coprod Z_{X}))$ and $ \Div^{G}_{rat}(B'\times_{S}\overline{Y},B'\times_{S}(Y_{\infty}\coprod Z_{Y}))$, where $\phi_{X}$ is an invariant rational function which is $1$ on $A'\times_{S} X_{\infty}$.

On the other hand the inclusions 
$j_{A}$ and $j_{B}$ are represented by the classes $\Delta A'$ and $\Delta B'$. 
 It follows that the differences 
$\lambda_{A}\circ i'-j_{A}\in G\Cor_{k}(A', A)$ and $\lambda_{B}\circ i' - j_{B}\in G\Cor_{k}(B',B)$ are represented by  the classes $\div(\phi_{X})$ and $\div(\phi_{Y})$ respectively.

The composition $\psi f' \in G\Cor_{k}(B',B)$ is represented by the divisor of the rational function $hf'$ which is $1$ on $B'\times_{S} Y_{\infty}$ and $r_{M}f' = r_{Y}$ on $B'\times_{S} Z_{Y}$.  We thus have $\psi f' = \lambda_{B}\circ i' - j_{B}$ in 
$\Div^{G}_{rat}(B'\times_{S}\overline{Y},\,B'\times (Y_{\infty}\coprod Z_{Y}))$.

Now the composition $f\psi \in G\Cor_{k}(A',A)$ represents the push forward of $\psi$ along $H_{0}^{Sus}(G;A'\times_{S}B/A') \to H_{0}^{Sus}(G;A'\times_{S} A/A'$. By Lemma \ref{lem:normcomm} this is represented by the norm $N(h^{-1})$. Since $h^{-1}$ is $1$ on $f^{-1}(X_{\infty})\subseteq Y_{\infty}$, $N(h)=1$ on $A'\times_{S} X$. By the following lemma we have that $N(h) = r_{X}$ on $A'\times_{S} Z_{X}$ which yields the desired equality $f\psi = \lambda_{A}\circ i' - j_{A}\in G\Cor_{k}(A',\, A)$.
\end{proof}

\begin{lemma}
 Let $f:U\to V$ be a finite equivariant map with $U$ and $V$ normal. Suppose that $Z\subseteq V$ and $Z'\subseteq U$ are reduced closed subschemes such that the induced map $Z'\to Z$ is an isomorphism and $U\to V$ is \'etale in a neighborhood of $Z'$. If $h\in \mcal{O}^{*}(U)^{G}$ is $1$ on $f^{-1}(Z)-Z'$ then $N(h)|_{Z}$ and $h|_{Z'}$ are identified by $Z'\iso Z$.
\end{lemma}
\begin{proof}
 This follows immediately from the nonequivariant statement [MVW, Lemma 21.10]. (i.e., forget the $G$-action then [MVW, Lemma 21.10] tells us that $N(h)|_{Z}$ and $h|_{Z'}$ are identified by $Z'\iso Z$.
\end{proof}

We finish this section with the following useful application of Theorem \ref{thm:complex}.

\begin{theorem}\label{thm:hizero2}
 Let $F$ be a homotopy invariant presheaf, $J$ a smooth equivariantly irreducible zero-dimensional $G$-scheme and $W$ a $G$-representation. Then for any open invariant 
 $U\subseteq L:=J\times W$ we have
$$
H^{i}_{GNis}(U, F_{GNis}) = \begin{cases}
                                     F(U) & i = 0 \\
				      0  & i>0 .
                                    \end{cases}
$$
\end{theorem}
\begin{proof}
Corollary \ref{cor:cdbd} implies that $H^{i}_{GNis}(U,F) = 0$ for $i>1$. 
Consider an equivariant distinguished square $Q= Q(U,V,A)$,
$$
\xymatrix{
B\ar[d]\ar[r] & V \ar[d] \\
A \ar[r] & U .
}
$$

There is an equivariant embedding of $V$ into a smooth projective curve $\overline{V}$ with $G$-action which is finite over $\P(L\oplus 1)$. Indeed, ignoring the group action on $V$ there is an embedding into a smooth projective curve $\overline{V}$. Rational maps between smooth projective curves extend uniquely to morphisms which implies that $\overline{V}$ inherits a $G$-action from $V$ and maps equivariantly and finitely to $\P(L\oplus 1)$.

The square $Q$ comes from the equivariant covering morphism of equivariant standard triples, $(\overline{V}, V_{\infty}, Z) \to (\P(L\oplus 1), U_{\infty}, Z)$ where $V_{\infty} = \overline{V} \setm V$, $U_{\infty} = \P(L\oplus 1)\setm U$, and $Z =- U\setm A$. The triple  $(\P(L\oplus 1), U_{\infty}, Z)$ is split over $U$ by Lemma \ref{lem:splittrip4}. Applying Theorem \ref{thm:complex} with $Q=Q'$ we see that the Mayer-Vietoris sequence 
$$ 
0\to F(U) \to F(A) \oplus F(V) \to F(B) \to 0
$$
is split exact. This implies that $F$ is a sheaf in the equivariant Nisnevich topology on $U$ and that $\check{H}^{1}(\mcal{U}/U, F) = 0$ for any cover $\mcal{U}$ coming from a distinguished Nisnevich square.

We claim that any equivariant Nisnevich cover of $U$ can be refined by one coming from an equivariant distinguished square, and consequently $\check{H}^{1}(\mcal{U}, F) = 0$. This will finish the proof since $H^{1}_{GNis}(U,F) = \check{H}^{1}(U,F)$.
First, since $F$ takes disjoint unions to sums we can replace a cover $\{V_{i}\to U\}$ by a single cover $f:V'\to U$. Indeed there is a dense invariant open $A\subseteq U$ over which $f$ has a splitting. The complement $Z = U\setm A$ is a finite set of closed points and we choose a splitting $Z\subseteq V_{Z}$. Let $V= V'\setm(V_{Z}\setm Z)$ then $Q(U,V,A)$ is a distinguished square and the associated cover refines $f:V'\to U$.  
\end{proof}

\section{Homotopy invariance of cohomology}\label{sec:hi}
In this section we show that 
under the assumption of Condition \ref{conmod}
equivariant Nisnevich cohomology with coefficients in a homotopy invariant presheaf with transfers is again a homotopy invariant  presheaf with equivariant transfers. This result is the equivariant analogue of the fundamental technical result in Voevodsky's machinery of presheaves with transfers.   

Unless specified otherwise,  $G$ is assumed to satisfy Condition \ref{conmod} throughout this section.

\begin{proposition}\label{prop:gnishi}
 Let $F$ be a homotopy invariant presheaf with equivariant transfers. Then $F_{GNis}$ is also a homotopy invariant presheaf with equivariant transfers.
\end{proposition}
\begin{proof}
By Theorem \ref{thm:gnistrans}, $F_{GNis}$ is a presheaf with equivariant transfers. To show homotopy invariance it suffices to show that $i^{*}:F_{GNis}(X\times\A^{1})\to F_{GNis}(X)$ is injective for any equivariantly irreducible $X$, where $i:X\to X\times \A^{1}$ is the inclusion at $0\in \A^{1}$. It suffices to do this locally in the equivariant Nisnevich topology, so we may assume that $X$ is affine semilocal with a single orbit. 
 Let $\mcal{Z}\subseteq X$ be the set of generic points with induced $G$-action.  We have a commutative square
$$
\xymatrix{
F_{GNis}(X\times\A^{1})\ar[r]\ar[d] & F_{GNis}(X) \ar[d] \\
F_{GNis}(\mcal{Z}\times\A^{1}) \ar[r]^-{\iso} & F_{GNis}(\mcal{Z}).
}
$$
We may view $F$ as a homotopy invariant presheaf with equivariant transfers on $G\Sm/K$, where $K=k(X)^{G}$. 
Theorem \ref{thm:hizero2} implies that $F$ is an equivariant Nisnevich sheaf on $\mcal{Z}\times \A^{1}$ and therefore the bottom horizontal arrow is an isomorphism.
The vertical arrows are injective by Theorem \ref{thm:shfinj} and so $i^{*}$ is injective and thus $F_{GNis}$ is homotopy invariant.
\end{proof}

\subsection{Equivariant contractions}
If $F$ is a presheaf with transfers on $\Sm/k$ the contraction 
$F_{(-1)}(X) := F(X\times\A^{1}\setm\{0\})/F(X\times\A^{1})$ plays an important role in the study of presheaves with transfers. We introduce an equivariant analogue and establish a few basic results concerning equivariant contractions.

\begin{definition}
 Let $F$ be a presheaf on $G\Sm/k$ and $W$ a representation of $G$. Define the presheaf $F_{(-W)}$ by
$$
F_{(-W)}(X) = \coker(F(X\times \A(W)) \to F(X\times \A(W)\setm\{0\})).
$$
\end{definition}

When $F$ is a presheaf with equivariant transfers then so is $F_{(-W)}$ since it is the quotient of such presheaves. Similarly if $F$ is homotopy invariant then $F_{(-W)}$ is as well.

Nonequivariantly the  projection $X\times \A^1 \to X$ is split by including at $1\in \A^{1}$, inducing 
a decomposition $F(X\times \A^{1}\setm 0) = F(X) \oplus F_{(-1)}(X)$ whenever $F$ is homotopy invariant. When $W$ is a representation with $W^{G} = 0$ then there is no such equivariant splitting. Nonetheless when $F$ is a presheaf with equivariant transfers we still obtain this decomposition, at least for affine $X$.

\begin{proposition}\label{prop:decomp}
 Let $F$ be a homotopy invariant presheaf with equivariant transfers on $G\Sm/k$. Let $S$ be a smooth affine $G$-scheme over $k$ and $W$ be a one-dimensional representation. Then there is an equivariant finite correspondence $\lambda:S\times \A(W) \to S\times \A(W)\setm 0$ inducing a decomposition
$$
F(S\times \A(W)\setm\{ 0\}) = F(S) \oplus F_{(-W)}(S).
$$ 
Moreover, this decomposition is natural for equivariant maps $S'\to S$, where $S'$ is affine.
\end{proposition}
\begin{proof}
We have an equivariant standard triple $(S\times \P(W\oplus 1) \to S, S\times \infty, S\times 0)$. By Lemmas \ref{lem:splittrip4} and \ref{lem:pulllem} this equivariant triple is equivariantly split over $X = S\times \A(W)$. Applying Proposition \ref{prop:lambda} yields the correspondence $\lambda$ which induces the splitting $F(S\times \A(W)\setm 0) \to F(S\times \A(W))$. 
\end{proof}

\begin{proposition}\label{prop:shfcon}
 Let $F$ be a homotopy invariant presheaf with equivariant transfers and $W$ a one-dimensional representation. Then 
$$
(F_{GNis})_{(-W)}(S) = (F_{(-W)})_{GNis}(S)
$$ 
 for any smooth affine Henselian $G$-scheme over $k$ with a single closed orbit. 
\end{proposition}
\begin{proof}
By Proposition \ref{prop:gnishi}, $F_{GNis}$ is a homotopy invariant sheaf with equivariant transfers. 
Observe that $(F_{(-W)})_{GNis}\to (F_{GNis})_{(-W)}$ is a morphism of presheaves with equivariant transfers. Applying Corollary \ref{cor:vanfield}  to the kernel and cokernel of this map, it suffices to show that 
$(F_{(-W)})_{GNis}(J)= (F_{GNis})_{(-W)}(J)$ for any essentially smooth zero dimensional $G$-scheme over $k$. The left-hand side is by definition $F(J\times \A(W)\setm\{0\})/F(J\times \A(W))$. The right-hand side is $F_{GNis}(J\times \A(W)\setm\{0\})/F_{GNis}(J\times \A(W))$. Applying Proposition \ref{prop:gnishi} and Theorem \ref{thm:hizero2} shows the two sides are equal.
\end{proof}

\begin{definition}
 Let $i:Z\hookrightarrow Y$ be an invariant closed embedding with open complement $j:V\subseteq Y$ and $F$ a presheaf. Define the equivariant Nisnevich sheaf $F_{(Y,Z)}$ on $Z$ as in the nonequivariant case. That is let $K_{(Y,Z)} = K$ be the cokernel of $F \to j_{*}j^{*}F$ and define $F_{(Y,Z)} = (i^{*}K)_{GNis}$. 
\end{definition}

Since sheafification is exact we have an exact sequence
\begin{equation}\label{eqn:excon}
F_{GNis} \to (j_{*}j^{*}F)_{GNis} \to i_{*}F_{(Y,Z)} \to 0.
\end{equation}

\begin{lemma}\label{lem:coh}
 For $n\geq 0$ we have $H^{n}_{GNis}(-,F)_{(Y,Z)} = i^{*}R^{n}j_{*}F$.
\end{lemma}
\begin{proof}
 The same argument as in \cite[Example 2.3.8]{MVW} works here. Namely, we have $H^{n}_{GNis}(-,F)_{GNis} = 0$ and therefore $i_{*}H^{n}(F)_{(Y,Z)} \iso (j_{*}j^{*}H^{n}(F))_{GNis} = R^{n}j_{*}F$. Since $i^{*}i_{*} = \id$ the result follows.
\end{proof}

Let $i:S\hookrightarrow S\times \A(W)$ be the invariant closed embedding determined by $0\in \A(W)$. Then $F_{(-W)}(U) = K(U\times \A(W))$. We obtain by adjunction the map $K(U\times \A(W)) \to i_{*}i^{*}K(U\times \A(W)) = i^{*}K(U)$. Therefore we have the map of sheaves on $S$
$$
(F_{(-W)})_{GNis} \to F_{(S\times \A(W), S\times 0)}.
$$

\begin{proposition}\label{prop:zerocon}
 Let $F$ be a homotopy invariant presheaf with equivariant transfers, $W$ a one-dimensional $G$-representation, and $S$ a smooth $G$-scheme. Then we have an isomorphism
$$
(F_{(-W)})_{GNis}|_{S} \xrightarrow{\iso} F_{(S\times \A(W), S\times 0)}.
$$
\end{proposition}
\begin{proof}
We use the argument of \cite[Proposition 23.10]{MVW}. We need to compare $F_{(-W)}$ and $j_{*}j^{*}F/F$ in an invariant neighborhood of an orbit $Gs$ of a point $s$ in a smooth affine $G$-scheme $S$. The equivariant standard triple 
$T = (\P(W\oplus 1)_{S}, S\times \infty , S\times 0)$ is split over $S\times \A(W)$ by Lemma \ref{lem:splittrip4}. Let $U$ be an affine invariant neighborhood of $S\times 0$ in $S\times \A(W)$ and let $T_{U} = (\P(W\oplus 1)_{S}, \P(W\oplus 1)_{S} - U, S\times 0)$. We need to show that by shrinking $S$ there is an invariant open affine neighborhood of $(\P(W\oplus 1)_{S} - U) \cup S\times 0$. It follows that $T_{U}$ is an equivariant standard triple.

There is an invariant open $V\subseteq \P(W\oplus 1)$ so that $Gs\times V$ contains $Gs\times 0 $ and the finite invariant set $\P(W\oplus 1)_{Gs} - U_{Gs}$. The complements of $U$ and $S\times V$ intersect in a closed subset, disjoint from the fiber $\P(W\oplus 1)_{Gs}$. Since $\P(W\oplus 1)_{S}$ is proper over $S$ we may shrink $S$ around $Gs$ to assume that the complements are disjoint. Then $S\times V$ contains both $\P(W\oplus 1)_{S} - U$ and $S\times 0$ as needed.

The identity on $\P(W\oplus 1)$ is an equivariant covering morphism of triples $T_{U}\to T$. Let $U_{0} = U-S\times 0$. Consider the distinguished square $Q$
$$
\xymatrix{
U_{0} \ar[r]\ar[d] & U \ar[d] \\
S\times \A(W)\setm 0 \ar[r]^-{j} & S\times W.
}
$$
Write $Q'$ for the same square. The identity $Q'=Q$ comes from the map of triples $T_{U}\to T$, see Example \ref{ex:cov}. Applying Theorem \ref{thm:complex}  we have a split exact Mayer-Vietoris sequence
$$
0\to F(S\times \A(W)) \to F(S\times W\setm 0)\oplus F(U) \to F(U_{0})\to 0.
$$ 
This together with the homotopy invariance of $F$ implies we have a pushout square
$$
\xymatrix{
F(S\times \A(W)) \ar[r]\ar[d] & F(U) \ar[d]\\
F(S\times \A(W)\setm 0) \ar[r] & F(U_{0}).
}
$$
In particular $F(U) \to F(U_{0})$ is injective and $F_{(-W)}(S) = F(U_{0})/F(U)$.

Note that $j:S\times \A(W)\setm 0 \to S\times W$ has $j_{*}j^{*}F(U) = F(U_{0})$ and thus 
$j_{*}j^{*}F/F(U) = F(U_{0})/F(U)$ and the result follows by passing to the limit over $U$ and $S$. 
\end{proof}

\begin{lemma}\label{lem:nisbr}
 Let $f:Y\to X$ be an equivariant \'etale morphism and $Z\subseteq X$ an invariant closed subscheme such that $f^{-1}(Z) \to Z$ is an isomorphism. Then for any presheaf $F$ we have 
$$
F_{(X,Z)} \xrightarrow{\iso} F_{(Y,f^{-1}(Z))}.
$$
\end{lemma}
\begin{proof}
It is enough to check the isomorphism on stalks. We may thus assume that $Y$, $X$ are semilocal Henselian $G$-schemes, $X$ has a single closed orbit $Gx$, and $Z$ is nonempty. Since $f^{-1}(Z) \iso Z$ and $Gx\subseteq Z$, it follows that $Y$ also has a single closed orbit and that $Y\iso X$. 
\end{proof}

Recall if $G$ acts on the ring $R$, we write 
$R^{\#}[G]$ for the twisted group ring (see Remark \ref{rem:skew}). If $H$ acts on the field $L$ and $W$ is a $k[H]$-module then $W_{L}$ is a $L^{\#}[H]$-module via $r[g](w\otimes x) = gw\otimes r(gx)$.

\begin{lemma}
 Let  $Z\subseteq X$ be an equivariant closed embedding of smooth affine $G$-schemes over $k$ and $x\in Z$ a closed point.  Suppose that there are $G$-representations $W_{2}\subseteq W_{1}$  and isomorphism $f:T_{x}X\iso (W_{1})_{k(x)}$ of $k(x)^{\#}[G_{x}]$-modules which restricts to a $k(x)^{\#}[G_{x}]$-module isomorphism $T_{x}Y\iso (W_{2})_{k(x)}$. Then there is an invariant open neighborhood $U$ of $x$ and an equivariant Cartesian diagram
 $$
 \xymatrix{
 U\cap Z \ar[r]\ar[d] & U \ar[d] \\
 \A(W_{2}) \ar[r] & \A(W_{1})
 }
 $$
 with \'etale vertical arrows.
\end{lemma}
\begin{proof}
Let $L/k(x)$ be a finite extension such that the composite $L/k$ is Galois with Galois group $\Gamma$. The schemes $X_{L}$ and $Z_{L}$ will be considered as $G\times\Gamma$-schemes over $k$ via the diagonal action. We will construct a $G\times\Gamma$-equivariant maps $\phi$, $\phi'$ which fit into a commutative square
 $$
 \xymatrix{
 Z_L \ar@{^{(}->}[r]\ar[d]_{\phi'} & X_L \ar[d]^{\phi}\\
 \A(W_{2})_L\ar@{^{(}->}[r] & \A(W_1)_L
 }
 $$
and are equivariant at each point of the $G\times\Gamma$-orbit of $y$, where $y\in X_L$ lies over $x$. The set of points at which $\phi$ is \'etale is an open (and invariant) subset of $X_L$ which contains the orbit of $y$. Further shrinking this set equivariantly if necessary, we find an invariant open subset $\tilde{U}\subseteq X_L$ such $\phi$ is \'etale on $\tilde{U}$ and $\phi'$ is \'etale on $\tilde{U}\cap Z_{L}$. Now $\tilde{U}\times_{\A(W_1)_L}\A(W_2)_L$ is a disjoint union $\tilde{U}\cap Z_{L} \coprod C$. Replacing $\tilde{U}$ by $\tilde{U}-C$ we may assume that $\tilde{U}$ satisfies $\tilde{U}\cap Z\subseteq \tilde{U}\times_{\A(W_1)_L}\A(W_2)_L$. Galois descent then yields the desired $G$-equivariant square of the lemma.

It remains to construct the desired $G\times\Gamma$-equivariant square above. 
Let $y_1\in Z_{L}$ be a point lying over $x$ and let $\{y_1,y_2,\ldots, y_n\}$ be the $G\times\Gamma$-orbit of $y_1$. 
Let $R$ be the coordinate ring of $X$ and $I\subseteq R$ the defining ideal of $Z$.
Write $m_i$ and $\overline{m}_i$ respectively for the maximal ideals $y_i$ in $R_L=R\otimes_{k} L$ and $(R/I)_L=R/I\otimes_k L$. 
The ideals $\cap m_{i}\subseteq R_L$ and $\cap \overline{m}_i\subseteq (R/I)_L$ are $G\times\Gamma$-invariant. Consider the morphism $\cap m_i\to m_1/m_1^2\times \cdots m_n/m_n^2$, induced by the quotients $m_i\to m_i/m_i^2$. Using the Chinese Remainder Theorem, we see that it is surjective.

Let $S\subseteq G\times \Gamma$ be the set-theoretic stabilizer of $y_1$. It is the subgroup $S\subseteq G_{x}\times \Gamma$ consisting of pairs $(g,\gamma)$ such that the two maps 
$k(x)\to L$ given by $\iota g$ and $\gamma \iota$ are equal (where $\iota:k(x)\subseteq L$ is the embedding chosen at the beginning of the proof).
Let  $\alpha_1,\ldots, \alpha_n$ be left coset representatives for $(G\times\Gamma)/S$. For an element $\beta=(g,\gamma)$ of  $G\times \Gamma$ write $\beta \alpha_i = \alpha_{j(i)} s_{i}$, for appropriate indices $j(i)$ and $s_i\in S$. 
If $M$ is an $L^{\#}[S]$-module we obtain an induced $L^{\#}[G\times\Gamma]$-module (here the action of $G\times \Gamma$ on $L$ is via the projection to $\Gamma$). As in the case of an ordinary group ring, we may describe the induced module 
$\Ind(M):=L^{\#}[G\times\Gamma]\otimes_{L^{\#}[S]}M$ as the direct sum $\oplus ([\alpha_{i}]\otimes M)$ of copies of $M$ with basis $\{[\alpha_i]\}$.  The $L^{\#}[G\times \Gamma]$ module structure on $\Ind(M)$ is determined by the equations
$[\beta]([\alpha_i]\otimes m_i) = ([\alpha_{j(i)}]\otimes s_{i}m_{i})$, for $\beta\in G\times\Gamma$ and $r([\alpha_i]\otimes m_i) = ([\alpha_{i}]\otimes (\alpha_{i}^{-1}r)m_{i})$ for $r\in L$.

We have an isomorphism 
$m_1/m_1^2\times\cdots\times m_n/m_n^2 \iso 
\Ind(m_1/m_1^2)$ given by sending $\overline{r}_i\in m_i/m_i^{2}$ to $[\alpha_i]\otimes \overline{\alpha_{i}^{-1}r_i} $. We thus obtain a surjection $\cap m_i\to \Ind(m_1/m_1^2)$ which is a surjection of $L^{\#}[G\times\Gamma]$-modules. In a similar fashion we obtain a surjection $\cap \overline{m}_i\to \Ind(\overline{m}_1/\overline{m}_1^2)$ of $L^{\#}[G\times\Gamma]$-modules. 

Now the isomorphism $T_xX\iso (W_1)_{k(x)}$ of $k(x)^{\#}[G_x]$-modules yields the isomorphism $T_{y_1}(X_{L})=T_{x}(X)\otimes_{k(x)}L \iso (W_1)_L$ of $L^{\#}[S]$-modules, which restricts to an isomorphism $T_{y_1}Z_{L}\iso (W_2)_{L}$. Since $m_1/m_1^2\iso (T_{y_1}X_L)^{\vee}$ and $\overline{m}_1/\overline{m}_1^2 \iso (T_{y_1}Z_L)^{\vee}$ we obtain a commutative diagram of $L^{\#}[G\times \Gamma]$-modules
$$
\xymatrix{
R_L \ar@{->>}[d]
& \cap m_i \ar[d]\ar@{->>}[r]\ar@{_{(}->}[l]
& \Ind(m_1/m_1^2) \ar@{->>}[d]\ar[r]^{\iso}
& \Ind((W_1)_L^{\vee}) \ar@{->>}[d]
\\
(R/I)_L 
& \cap \overline{m}_{i}  \ar@{->>}[r]\ar@{_{(}->}[l] 
& \Ind(\overline{m}_1/\overline{m}_1^2) \ar[r]^{\iso}
& \Ind((W_2)_L^{\vee}).
}
$$

The kernel of the action of $G\times \Gamma$ on $L$ is equal to  $G$. Since $|G|$ is invertible in $L$, the ring 
$L^{\#}[G\times \Gamma]$ is semi-simple, see e.g. \cite[Lemma 1.3]{Kunzer}. We may therefore choose compatible splittings $\Ind((W_{1})_{L}^{\vee})\to \cap m_{i}$ and $\Ind((W_{2})_{L}^{\vee})\to \cap\overline{m}_{i}$ to the horizontal arrows. We have as well a $L^{\#}[G\times\Gamma]$-module map $(W_1)_L^{\vee}\to \Ind((W_1)_L^{\vee})$ given by $\omega \mapsto ([\alpha_i]\otimes \alpha_i^{-1}\omega)_i$ and similarly for $(W_2)_L^{\vee}\to \Ind((W_2)_L^{\vee})$.  
We thus obtain a commutative square 
$$
\xymatrix{
(R)_L \ar[d] & \Sym_{L}((W_1)_L^{\vee}) \ar[l]\ar[d] \\
(R/I)_L  & \Sym_{L}((W_2)_L^{\vee}) \ar[l]
}
$$
of $k$-algebras with $G\times\Gamma$-action.
Tracing through the construction of these maps, we see that the compositions $(W_1)_L^{\vee} \to \cap m_i\to m_i/m_i^2$ and  $(W_2)_L^{\vee}\to \overline{m}_i/\overline{m}_{i}^{2}$ are $L^{\#}[G\times\Gamma]$-module isomorphisms. This implies that the horizontal arrows are \'etale at the $m_i$ and therefore applying $\spec(-)$ to this square we obtain the desired square of $G\times\Gamma$-schemes over $k$.

\end{proof}

\begin{theorem}\label{thm:loc2}
Let $X$ be a smooth affine $G$-scheme and $Z\subseteq X$ a closed invariant smooth $G$-scheme of codimension one. Let $x\in Z$ be a closed point. Let $W$ be a $G$-representation defined over $k$ such that there is an $I_{x}$-equivariant isomorphism $W_{k(x)} \iso T_{x}X/T_{x}Z$.  Then there is an invariant open neighborhood $U\subseteq X$  of $x$ such that for any smooth $G$-scheme $T$ we have isomorphisms of sheaves on $(U \cap Z) \times T$
$$
F_{(U\times T, (U\cap Z)\times T)} \iso (F_{(-W)})_{GNis}.
$$
\end{theorem}
\begin{proof}
We need to see that Condition \ref{conmod} implies that the hypothesis of the previous lemma are satisfied. It suffices to see that every irreducible $k(x)^{\#}[G_{x}]$-module is isomorphic to $M\otimes_{k}{k(x)}$ for some $k[G]$-module $M$. 
Condition \ref{conmod} implies that there are irreducible $k[G]$-modules $M_{1},\dots, M_{d}$ which form a complete set of irreducible $k[I_{x}]$-representations. Each $M_{i}$ is one-dimensional and $d=|I_{x}|$.  Each $M_{i}':=M_{i}\otimes_{k}k(x)$ is an irreducible $k(x)^{\#}[G_{x}]$-module. Any $k(x)^{\#}[G_{x}]$-module isomorphism $M_{i}'\iso M_{j}'$ is also a 
$k(x)[I_{x}]$-module isomorphism. Since the $\{M_{i}'\}$ form 
$d$-distinct irreducible $k(x)[I_{x}]$-modules, they are also $d$-distinct irreducible $k(x)^{\#}[G_{x}]$-modules. We claim that this is a complete list of irreducible $k(x)^{\#}[G_{x}]$-modules. First note that we have $\mathrm{End}_{k(x)^{\#}[G_{x}]}(M_{i}') = F$, where 
$F= k(x)^{G_{x}}$ is the fixed field. Write $n = [L:F] = |G_{x}/I_{x}|$. Since $k(x)^{\#}[G_{x}]$ is semi-simple, we may write it as a direct sum involving all of the irreducible modules. Since  $\mathrm{End}_{F}(M_{i}') = (M_i')^{n}$, each $M_{i}'$ appears with multiplicity $n$ in this decomposition. Comparing dimensions (as $k(x)$-vectorspaces) we see that the $M_{i}'$ form a complete list of irreducible $k(x)^{\#}[G_{x}]$.

Let $W_1$ and $W_2$ be $G$-representations satisfying the hypothesis of the previous lemma. Set $W=W_1/W_2$. 
By the previous lemma, after shrinking $X$ around $x$, there is  an equivariant Cartesian square, 
$$
\xymatrix{
Z\ar[r]\ar[d] & X \ar[d] \\
 \A(W_2)\ar[r] & \A(W_1)
}
$$ 
where the vertical maps are \'etale. Proceeding as in 
\cite[Theorem 23.12]{MVW} yields the result.
\end{proof}

\subsection{Proof of homotopy invariance}\label{sub:hi}

We remind the reader that $G$ is assumed to satisfy Condition \ref{conmod}.

\begin{theorem}\label{thm:hi}
 Let $F$ be a homotopy invariant presheaf with equivariant transfers on $G\Sm/k$. Then 
$H^{n}_{GNis}(- ,F_{GNis})$ is also a homotopy invariant presheaf with equivariant transfers.
\end{theorem}
\begin{proof}
 By Theorem \ref{thm:gniscoh}, $H^{n}_{GNis}(- ,F_{GNis})$ is a presheaf with equivariant transfers and it remains to verify that it is homotopy invariant. 
The case $n=0$ is Proposition \ref{prop:gnishi}.
We may thus assume that $F = F_{GNis}$ and we proceed by induction on $n$. 
Let $X$ be a smooth $G$-scheme and consider the map $\pi:X\times\A^{1} \to X$ and  the Leray spectral sequence
$$
H^{p}_{GNis}(X, R^{q}\pi_{*}F) \Rightarrow H^{p+q}_{GNis}(X\times\A^{1}, F).
$$
We have that $\pi_{*}F = F$ since $\pi_{*}F(U) = F(U\times\A^{1}) \iso F(U)$. By induction we have that $R^{q}\pi_{*}F = 0$ for $0<q<n$. The spectral sequence collapses by Theorem \ref{thm:Rnzero}, yielding the desired isomorphism $H^{n}_{GNis}(X ;F_{GNis}) = H^{n}_{GNis}(X\times\A^{1} ;F_{GNis})$.
\end{proof}

\begin{lemma}\label{lem:filt}
 Let $X$ be a smooth $G$-scheme over $k$, $Z\subseteq X$ a closed invariant subset such that $\codim(Z)\geq 1$, and $x$ a point of $X$. Then there is an open invariant neighborhood $U\subseteq X$ of $x$ and a sequence of invariant reduced closed subschemes $\emptyset = Y_{-1} \subseteq Y_{0} \subseteq Y_{1} \subseteq \cdots \subseteq Y_{k}$ in $U$ satisfying the following two properties.
\begin{enumerate}
 \item The $G$-schemes $Y_{i}-Y_{i-1}$ are smooth invariant divisors on $U-Y_{i-1}$.
\item $U\cap Z \subseteq Y_{k}$.
\end{enumerate}
\end{lemma}
\begin{proof}
 The argument is similar to \cite[Lemma 4.31]{V:pth}. The key point is that 
under our assumptions, there is a smooth equivariant curve $p:U\to V$ and so the induction argument of loc.~cit.~applies here.

\end{proof}

Now we are ready to prove the vanishing of $R^{n}\pi_{*}F$.

\begin{theorem}\label{thm:Rnzero}
Let $X$ be a smooth $G$-scheme over $k$ and $F$  a homotopy invariant equivariant Nisnevich sheaf with equivariant transfers. Assume that 
 that $R^{q}\pi_{*}F = 0$ for $0<q<n$, and that $H^{p}_{GNis}(-, F)$ is  homotopy invariant 
for $p<n$.
Then $R^{n}\pi_{*}F = 0$ as well.
\end{theorem}
\begin{proof}
 We may assume that $X$ is equivariantly irreducible. We need to show that given an $\alpha\in H^{n}_{GNis}(X\times \A^{1}, F)$ it becomes zero on an equivariant Nisnevich cover of $X$.  Let $J$ denote the set of generic points of $X$. By Theorem \ref{thm:hizero2}, $H^{n}_{GNis}(J\times \A^{1}, F)=0$. This implies that 
 there is an open dense $V\subseteq X$ such that $\alpha|_{V}$ vanishes. Let $Z= X-V$ with its reduced structure. It now suffices to show that 
$$
H^{n}_{GNis}(X\times\A^{1}, F) \to H^{n}_{GNis}((X-Z)\times\A^{1}, F)
$$
is injective locally in the equivariant Nisnevich topology on $X$. 
Using Lemma \ref{lem:filt} we may assume that $Z$ is a smooth invariant divisor. We are thus reduced to showing that  
\begin{equation}\label{eqn:inj}
H^{n}_{GNis}(X'\times\A^{1}, F) \to 
H^{n}_{GNis}((X'-Z')\times\A^{1}, F)
\end{equation}
is injective where $X'$ is a smooth affine Heneselian semilocal $G$-scheme over $k$ with a single closed orbit and $Z'\subseteq X'$ is a smooth invariant divisor.

Write  $i:Z'\to X'$ and
$j:U'=X'-Z'\to X'$. The map (\ref{eqn:inj}) factors as
$$
H^{n}_{GNis}(X'\times\A^{1}, F)  \xrightarrow{\tau} 
H^{n}_{GNis}(X'\times\A^{1}, j_{*}j^{*}F) 
\xrightarrow{\eta} H^{n}_{GNis}(X'-Z'\times\A^{1}, j^{*}F)
$$
(where we view $F$ as a sheaf on $X'\times\A^{1}$). We show that each of these maps is injective.

First we show that $\eta$ is injective. We begin by showing that $R^{q}j_{*}F = 0$ for $0<q<n$. By the inductive hypothesis we have that $H^{q}_{GNis}(-,F)$ is a homotopy invariant presheaf with equivariant transfers. Since $q>0$ we have $H^{q}_{GNis}(-,F)_{GNis} = 0$.  
 By Theorem \ref{thm:loc2} there is a $G$-representation $W$ such that $(H^{q}(F)_{(-W)})_{GNis} \iso H^{q}(F)_{(X'\times\A^{1}, Z'\times\A^{1})}$. By Proposition \ref{prop:shfcon} we have that 
$(H^{q}(F)_{(-W)})_{GNis} = ((H^{q}(F)_{GNis})_{(-W)}))_{GNis} =0$. Finally, by Lemma \ref{lem:coh}, we have
$$
R^{q}j_{*}F \iso i_{*}H^{q}(F)_{(X'\times\A^{1}, Z'\times\A^{1})} \iso i_{*}(H^{q}(F)_{GNis})_{GNis} = 0
$$
and so $R^{q}j_{*}F = 0$ as claimed. Now consider the Leray spectral sequence
$$
H^{p}_{GNis}(X'\times \A^{1}, R^{q}j_{*}(j^{*}F))\Rightarrow H^{p+q}_{GNis}(X'-Z', j^{*}F).
$$
Since $R^{q}j_{*}F = 0$ for $0<q<n$ we obtain an exact sequence
$$
0 \to H^{n}_{GNis}(X'\times\A^{1}, j_{*}j^{*}F) 
\xrightarrow{\eta} H^{n}_{GNis}(U'\times\A^{1}, j^{*}F) \to H^{0}(X'\times\A^{1}, R^{n}j_{*}j^{*}F).
$$
In particular $\eta$ is injective as required.

It remains to show that $\tau$ is injective as well. By Theorem \ref{thm:shfinj} we have an injection $F\to j_{*}j^{*}F$. Combining this with (\ref{eqn:excon}) we have an exact sequence
$$
0\to F \to j_{*}j^{*}F \to i_{*}F_{(X'\times\A^{1},Z'\times\A^{1})} \to 0.
$$
As noted in the previous paragraph we have $((F)_{(-W)})_{GNis} \iso (F)_{(X'\times\A^{1}, Z'\times\A^{1})}$ as sheaves on $Z'\times\A^{1}$. Consider the long exact sequence associated to the above short exact sequence,
\begin{align*}
 H^{n-1}(X'\times\A^{1},j_{*}j^{*}F)  &\to H^{n-1}(Z'\times\A^{1}, F_{(-W)}) 
\to H^{n}(X'\times\A^{1},F) \\
& \to H^{n}(X'\times\A^{1},j_{*}j^{*}F) \to H^{n}(Z'\times\A^{1}, F_{(-W)}).
\end{align*}
It suffices to show that $H^{n-1}(X'\times\A^{1},j_{*}j^{*}F)  \to H^{n-1}(Z'\times\A^{1}, F_{(-W)})$ is onto. For $n>1$ we have that $H^{n-1}(Z'\times\A^{1}, F_{(-W)}) = H^{n-1}(Z',F_{(-W)}) = 0$, by homotopy invariance of $F_{(-W)}$, the induction hypothesis and that $Z'$ is a semilocal Hensel $G$-scheme with one orbit. It remains to consider $n=1$ and show that 
$$
F(U'\times \A^{1})= H^{0}(X'\times\A^{1},j_{*}j^{*}F)  \to H^{0}(Z'\times\A^{1}, F_{(-W)})
$$
is surjective. By homotopy invariance of $F$ and $F_{(-W)}$ this map is identified with the map 
$F(U') \to F_{(-W)}(Z')$. 
By Theorem \ref{thm:loc2} and (\ref{eqn:excon}) we have a surjection 
$$
j_*j^{*}F \to i_{*}(F_{(-W)})_{GNis} \to 0
$$
which shows that $F(U') \to F_{(-W)}(Z')$ is surjective because $X'$ is a Henselian semilocal $G$-scheme with a single orbit. We conclude that $\tau$ is injective and the proof of the theorem is complete.
\end{proof}

\subsection{Applications of homotopy invariance}
As before, we assume that $G$ satisfies Condition \ref{conmod}.

\begin{theorem}\label{thm:vancurve}
 Let $F$ be a homotopy invariant presheaf with equivariant transfers on $G\Sm/k$. Let $S$ be a smooth affine semilocal $G$-scheme over $k$ with a single closed orbit and $W$ a one-dimensional representation. Then for any invariant open $U\subseteq \A(W)\subseteq \P(W\oplus 1)$ and any $n>0$,
 $$
 H^{n}_{GNis}(S\times U,F) = 0.
 $$ 
 In particular, $H^{n}_{GNis}(S,F) = 0$ for $n>0$.
\end{theorem}
\begin{proof}
 By Theorem \ref{thm:hi} $H^{n}_{GNis}(-,F)$ is a homotopy invariant presheaf with equivariant transfers, in particular the second statement follows from the first. By Corollary \ref{cor:vanfield} it suffices to show that $H^{n}_{GNis}(J\times U,F) = 0$ for any equivariantly irreducible zero dimensional $G$-scheme $J$ over $k$. This follows from Theorem \ref{thm:hizero2}.

\end{proof}

\begin{theorem}\label{thm:cont}
Let $F$ be a homotopy invariant sheaf with equivariant transfers on $G\Sm/k$. Let $W$ be a one-dimensional representation and $X$ a smooth $G$-scheme over $k$. Then
$$
H^{n}_{GNis}(X\times (\A(W)\setm 0),F) \iso H^{n}_{GNis}(X,F) \oplus 
H^{n}_{GNis}(X,F_{(- W)}).
$$
\end{theorem}
\begin{proof}
  Write $\pi:X\times (\A(W)\setm 0) \to X$ for the projection.
Consider the Leray spectral sequence  
$H^{p}_{GNis}(X, R^{q}\pi_{*}F) \Rightarrow
  H^{p+q}_{GNis}(X\times(\A(W)\setm 0), F)$.
 We have by Theorem \ref{thm:vancurve} that $H^{q}_{GNis}(S\times (\A(W)\setm 0),F) = 0$ for any smooth affine semilocal $G$-scheme $S$ with a single closed orbit and any $q>0$ and therefore this spectral sequence collapses.
    We therefore have that
$H^{n}_{GNis}(X\times (\A(W)\setm 0),F) \iso H^{n}_{GNis}(X,\pi_{*}F)$. Since $F$ is a sheaf, Proposition \ref{prop:decomp} is seen to imply that there is a decomposition $\pi_{*}F = F\oplus F_{(-W)}$ of sheaves on $X$ and  we are done.
\end{proof}

Our final application in this section is to show that the equivariant Nisnevich and equivariant Zariski cohomology with coefficients in a homotopy invariant presheaf with transfers agree.

\begin{theorem}\label{thm:cohagr}
Let $F$ be a homotopy invariant, equivariant Nisnevich sheaf with equivariant transfers on $G\Sm/k$. Then for any smooth quasi-projective $G$-scheme $X$ we have an isomorphism
$$
H^{n}_{GZar}(X,F) \iso H^{n}_{GNis}(X,F).
$$
\end{theorem}
\begin{proof}
Consider the Leray spectral sequence 
$H^{p}_{GZar}(X, \mcal{H}^{q})\Rightarrow H^{p+q}_{GNis}(X,F)$ where $\mcal{H}^{q}$ is the equivariant Zariski sheafification of the presheaf $U\mapsto H^{q}_{GNis}(U,F)$
on $X_{GZar}$. The result will follow if we see that $\mcal{H}^{q} = 0$ for $q>0$. 
The points of $X_{GZar}$ are the semilocal rings $\mcal{O}_{X,Gx}$ of an orbit  $Gx\subseteq X$. 
By Theorem \ref{thm:hi}, the presheaf $ H^{q}_{GNis}(-,F)$
is a homotopy invariant presheaf with equivariant transfers on $G\Sm/k$ and so the vanishing of $\mcal{H}^{q}$ follows from Theorem \ref{thm:vancurve}.
\end{proof}

\subsection{The group {$G=\Z/2$}}
Let $G=\Z/2$ and write $\G_{m}^{\sigma}$ for the $\Z/2$-scheme consisting of $\G_{m}$ with $\Z/2$-action given by $x\mapsto x^{-1}$.  Note that $\G_{m}^{\sigma}$ is not an invariant open in any representation and so the considerations above do not immediately apply to $\G_{m}^{\sigma}$. This scheme will be important for the cancellation theorem in the next section. We first record a few useful analogues of the previous results for this $\Z/2$-scheme.

Define $F_{(-\sigma)}$ to be the presheaf given by
$$
F_{(-\sigma)}(X) := \coker(F(X)\xrightarrow{\pi^{*}} F(X\times\G_{m}^{\sigma}))
$$
where $\pi:\G_{m}^{\sigma}\to \spec(k)$ is the structure map. Note that inclusion at $\{1\}$ yields an equivariant section $i_{1}:\spec(k)\to \G_{m}^{\sigma}$. We thus have 
$$
F(X\times \G_{m}^{\sigma}) = F(X) \oplus F_{(-\sigma)}(X).
$$
In particular if $F$ is a presheaf with equivariant transfers, homotopy invariant, or a sheaf then so if $F$.

Write $\sigma$ for the sign representation. Then $\P(\sigma\oplus 1)$ is $\P^{1}$ equipped with the action $[a:b]\mapsto [-a:b]$ and $\G_{m}^{\sigma}$ embeds into $\P(\sigma \oplus 1)$ as an open invariant subscheme (but is not contained in $\A(\sigma)$).

\begin{lemma}\label{lem:Z/2splittrip}
Let $X_{\infty}$ be the complement of $\G_{m}^{\sigma}\subseteq \P(\sigma\oplus 1)$ and $Z$ a finite, invariant set of closed points, disjoint from $X_{\infty}$. The triple $(\P(\sigma\oplus 1), X_{\infty}, Z)$ is split over any  invariant open subscheme $U\subseteq \G_{m}^{\sigma}$.  
\end{lemma}
\begin{proof}
 The argument is a simpler version of the argument given in Lemma \ref{lem:splittrip4}. The key point is that $\Pic((\G_{m}^{\sigma}\times \G_{m}^{\sigma})/\Z/2)=0$ by \cite[Corollary 12]{Magid}.
\end{proof}

\begin{theorem}\label{thm:Z2hizero}
 Let $F$ be a homotopy invariant presheaf, $S$ a smooth semilocal $\Z/2$-scheme over $k$ with a single closed orbit, and $U\subseteq \G^{\sigma}$ an invariant open subscheme. Then
$$
H^{i}_{GNis}(S\times U, F_{GNis}) = \begin{cases}
                                     F(U) & i = 0 \\
				      0  & i>0 .
                                    \end{cases}
$$
\end{theorem}
\begin{proof}
 For each $i$, $H^{i}_{GNis}(-,F_{GNis})$ is a homotopy invariant presheaf with equivariant transfers. It thus suffices by 
 Corollary \ref{cor:vanfield} to treat the case when $S$ is a zero dimensional smooth $\Z/2$-scheme. This case follows exactly as in the argument for Theorem \ref{thm:hizero2}, replacing the use of Lemma \ref{lem:splittrip4} with Lemma \ref{lem:Z/2splittrip}.
\end{proof}

\begin{proposition}\label{prop:Z2sfhcon}
 Let $F$ be a homotopy invariant presheaf with equivariant transfers on $G\Sm/k$. Then 
$$
(F_{GNis})_{(-\sigma)} = (F_{(-\sigma)})_{GNis}
$$ 
for any smooth semilocal $G$-scheme $S$ with a single closed orbit. 
\end{proposition}
\begin{proof}
 The argument is the same as in Proposition \ref{prop:shfcon}.  
\end{proof}

\begin{theorem}\label{thm:Z2cont}
Let $F$ be a homotopy invariant sheaf with equivariant transfers on $G\Sm/k$.  Then
$$
H^{n}(X\times \G_{m}^{\sigma},F) \iso H^{n}(X,F) \oplus H^{n}(X,F_{(-\sigma)}).
$$
\end{theorem}
\begin{proof}
 Write $\pi:X\times\G_{m}^{\sigma} \to X$ for the projection. By Theorem \ref{thm:Z2hizero} we have that $H^{q}(S\times\G_{m}^{\sigma},F) = 0$ for any smooth semilocal $G$-scheme $S$ over $k$ and $q>0$. Therefore  $R^{q}\pi_{*}F = 0$  for $q>0$ and so the Leray spectral sequence degenerates yielding
$H^{n}(X\times\G_{m}^{\sigma},F) \iso H^{n}(X,\pi_{*}F)$. Since $\pi_{*}F = F\oplus F_{(-\sigma)}$ we are done.
\end{proof}

\section{Cancellation Theorem}\label{sec:cancellation}
We apply the machinery developed in the previous sections in order to establish a $\Z/2$-equivariant version of Voevodsky's Cancellation Theorem. 
The argument given here is an equivariant modification of Voevodsky's argument in \cite{Voev:can}.

Let $Z$ be a smooth $\Z/2$-scheme  and $z\in Z$ an invariant rational point. Write $e:Z \to Z$ for the equivariant  idempotent morphism defined by the composition
$Z \to z \to Z$.

Let $(Z,z)$ and $e$ be as above. For $\Z/2$-schemes $X$, $Y$ define  
$$
\Cor_{k}(X\wedge Z, Y\wedge Z) := \{\mcal{V}\in \Cor_{k}(X\times Z, Y\times Z) \,|\, \mcal{V}\circ (id_{X}\times e) = 0 = (id_{Y}\circ e)\circ \mcal{V} \}.
$$
Note that this group inherits a natural $\Z/2$-action from that on 
$\Cor_{k}(X\times Z, Y\times Z)$ and as usual we write 
$$
\Z/2\Cor_{k}(X\wedge Z, Y\wedge Z):=  \Cor_{k}(X\wedge Z, Y\wedge Z)^{\Z/2}.
$$

This construction applies in particular to the $\Z/2$-varieties $(\G_{m},1)$ and $(\G_{m}^{\sigma},1)$, where $\G_{m}$ is considered with trivial action and $\G_{m}^{\sigma}$ has action given by $x\mapsto 1/x$.

Write $f_{i}$, $i =1,2$ for the projection $f_{i}:X\times\G_{m} \times Y \times \G_{m}\to \G_{m}$ to the $i$th copy of $\G_{m}$. Similarly write $f^{\sigma}_{i}:X\times\G^{\sigma}_{m} \times Y \times \G^{\sigma}_{m}\to \G_{m}$ to the $i$th copy of $\G_{m}$ (considered with trivial action). 
Define the rational functions $g_{n}$ and $g_{n}^{\sigma}$ by 
$$
g_{n} = \frac{f^{n+1}_{1} - 1}{f_{1}^{n+1} - f_{2}} \;\;\; \textrm{and} \;\;\; g^{\sigma}_{n} = 
\frac{(f^{\sigma}_{1})^{n+1} - 1}{(f^{\sigma}_{1})^{n+1} - (f_{2}^{\sigma})} .
$$ 
We consider the associated divisors $D(g_{n})$ and $D(g^{\sigma}_{n})$. Of course these are exactly the same divisor nonequivariantly, only the $\Z/2$-actions differ.
Observe that both of these are invariant divisors.

\begin{notation}
We write $\G$ to refer to one of $\G_{m}$ or $\G_{m}^{\sigma}$ and $g_{n}$ will correspondingly refer to either $g_{n}$ and $g_{n}^{\sigma}$. Similarly $D_{n}$ refers to either $D(g_{n})$ or $D(g_{n}^{\sigma})$.
\end{notation}

\begin{lemma}[{\cite[Lemma 4.1]{Voev:can}}]\label{vwd}
For any $\mcal{Z}\in \Cor_{k}(X\times \G, Y\times \G)$ there exists an $N$ such that for all $n\geq N$ the divisor $D_{n}$ of $g_{n}$ intersects $\mcal{Z}$ properly over $X$. 
\end{lemma}

For $\mcal{Z}\in \Cor_{k}(X\times \G, Y\times \G)$ the intersection
$\mcal{Z}\smash\cdot D_{n}$ is  an equidimensional relative cycle once $n$ is
large enough. Define $\rho_{n}(\mcal{Z})\in \Cor_{k}(X,Y)$ to be the projection to $X\times Y$ of this intersection. Observe that $g\rho_{n}(\mcal{Z}) = \rho_{n}(g\mcal{Z})$ for $g\in \Z/2$. Therefore if 
$\mcal{Z}\in \Z/2\Cor_{k}(X\times \G, Y\times \G)$ then 
$\rho_{n}(\mcal{Z})\in \Z/2\Cor_{k}(X,Y)$.
If both $\rho_n(\mcal{Z})$ and $\rho_m(\mcal{Z})$ are defined, they differ only up to equivariant $\A^{1}$-homotopy, see \cite{Voev:can}.

\begin{lemma}[{\cite[Lemmas 4.3, 4.4, 4.5]{Voev:can}}]\label{vcp} 
\hskip 1cm 
\begin{enumerate}
\item For $\mcal{W}\in \Z/2\Cor_{k}(X, Y)$ and $n\geq 1$ we have $\rho_{n}(\mcal{W}\times \id_{\G})= \mcal{W}$.
\item Let $e$ denote the composition $\G\to \{1\} \to \G$. Then $\rho_{n}(\id_{X}\otimes e) = 0$ for all $n\geq 0$ and all $g\in G$.
\item Let $\mcal{Z}\in \Z/2\Cor_{k}(X\times\G, Y\times \G)$  such that $\rho_{n}\mcal{Z}$ is defined. Let $\mcal{W}\in c(X',X)$ be arbitrary. 
Then $\rho_{n}(\mcal{Z}\circ (\mcal{W}\otimes \id_{\G}))$ is defined and 
\begin{equation*}
\rho_{n}(\mcal{Z}\circ (\mcal{W}\otimes \id_{\G})) = \rho_{n}(\mcal{Z})\circ \mcal{W},
\end{equation*}
where $\circ$ denotes composition of correspondences.
\item Let $\mcal{Z}\in \Z/2\Cor_{k}(X\times\G, Y\times \G)$ be such that $\rho_{n}\mcal{Z}$ is defined. Let $f:X'\to Y'$ be a morphism of schemes. Then $\rho_{n}(\mcal{Z}\times f)$ is defined and 
\begin{equation*}
\rho_{n}(\mcal{Z}\otimes f) = \rho_{n}(\mcal{Z})\otimes f.
\end{equation*}

\end{enumerate}
\end{lemma}

Write $\mcal{I}\in \Z/2\Cor_{k}(\G, \G)$ for the finite correspondence given by $\mcal{I} = \id_{\G}-e$.

\begin{proposition}\label{prop:comm}
 There is an equivariant homotopy 
 $H\in \Z/2\Cor_{k}(\G\wedge \G\wedge \A^{1}, \G\wedge\G)$ 
 such that $H_{0}-H_{1} = \tau - \id_{\G\wedge \G}$, where $\tau$ is the endomorphism of $\G\wedge \G$ which switches the factors.
\end{proposition}
\begin{proof}
It suffices to treat the case $\G=\G_{m}^{\sigma}$ since the proof for $\G=\G_{m}$ is contained in \cite[Proposition 3.2]{BV:V}. 

There is a canonical map $p:\G_{m}^{\sigma}\times \G_{m}^{\sigma}\to \Sym^{2}\G_{m}^{\sigma}$ with transpose $p^{t}$. 
Then $p^{t}p\in \Cor_{k}(\G_{m}^{\sigma}\times\G_{m}^{\sigma},\G_{m}^{\sigma}\times\G_{m}^{\sigma})$ is equal to $\id + \tau$. 
Write $\alpha:\G_{m}^{\sigma}\times \G_{m}^{\sigma}\to \G_{m}^{\sigma}\times \G_{m}^{\sigma}$ for the map defined by $(x,y)\mapsto (xy,1)$.

Define the $\Z/2$-scheme $M$ to have underlying scheme $M= \G_{m}\times \A^{1}$ and the $\Z/2$-action is specified by letting the nontrivial element 
act by $(x,y) \mapsto (x^{-1}, x/y)$. The map $\G_{m}^{\sigma}\times \G_{m}^{\sigma} \to M$, $(x,y)\mapsto (xy, x+y)$ is an equivariant isomorphism.
Therefore we have an equivariant $\A^{1}$-homotopy between the correspondences $p$ and $p\alpha$. 
Explicitly, we have an equivariant homotopy $H:M\times\A^{1} \to M$, given by $(x,y, t)\mapsto (x, t(1+x) +(1-t)y)$ which induces the desired homotopy. 

We therefore have that $\id +\tau = p^{t}p\wkeq  p^{t}p\alpha = \alpha + \tau\alpha$. 
Now  $\alpha+\tau\alpha$ takes values in $1\times \G_{m}^{\sigma} \cup \G_{m}^{\sigma}\times 1$ and therefore $\id = \tau$ in 
$\Z/2\Cor_{k}(\G_{m}^{\sigma}\wedge \G_{m}^{\sigma}, \G_{m}^{\sigma}\wedge\G_{m}^{\sigma})/\sim_{\A^{1}}$
\end{proof}

For $\mcal{W}\in \Z/2\Cor_k(X\wedge \G, Y\wedge\G)$ define $\mcal{W}\otimes^{(\tau)}\mcal{I} \in c(X\wedge \G\wedge\G_, Y\wedge\G\wedge\G)$ by
$\mcal{W}\otimes^{(\tau)}\mcal{I} = (\id_{Y}\otimes \tau)\circ (\mcal{W}\otimes \mcal{I})\circ (\id_{X}\times\tau).$

\begin{lemma}\label{ehp}
Let $\mcal{W}\in \Z/2\Cor_{k}(X\wedge \G, Y\wedge\G)$. There is an equivariant homotopy 
$$
\phi=\phi_{\mcal{W}} \in 
\Z/2\Cor_{k}(X\times \A^{1}\wedge \G\wedge\G, Y\wedge\G\wedge\G)
$$ 
such that $\phi_{0}-\phi_{1} = \mcal{W}\otimes^{(\tau)}\mcal{I}-\mcal{W}\otimes \mcal{I}$.
\end{lemma}
\begin{proof}
Let $H\in 
\Z/2\Cor_{k}(\G\times \G\times \A^{1}, \G\times\G)$ be the homotopy as in the previous proposition. We proceed as in \cite[Lemma 4.70]{Suslin:GraySS}. Let $\phi= \phi_{\mcal{W}}$ be defined by
\begin{align*}
\phi =  (\id_{Y}\otimes H)\circ [(\pm(\mcal{W}\otimes\mcal{I})\otimes \id_{\A^{1}}] & + \\ +  (\id_{Y}\otimes\tau)\circ(\mcal{W}\otimes\mcal{I})\circ(\id_{X}\otimes H) & .
\end{align*}
If $\mcal{W}$ is invariant then $\phi$ is also invariant. 
\end{proof}

Recall that if $F$ is a presheaf we write $C_{n}F$ for the presheaf $X\mapsto F(X\times\Delta^{n}_{k})$.

\begin{theorem}
Let $X$, $Y$ be smooth $\Z/2$-schemes over $k$. The homomorphism of simplicial abelian groups
\begin{equation*}
\Z/2\Cor_{k}(X\times\Delta^{\bullet}_{k},Y) \to \Z/2\Cor_{k}(X\times\Delta^{\bullet}_{k}\wedge \G , Y\wedge \G )
\end{equation*}
given by $\mcal{Z}\mapsto \mcal{Z}\otimes \mcal{I}$ is a weak equivalence. 
\end{theorem}
\begin{proof}
We follow the nonequivariant argument, 
\cite[Theorem 4.7]{Suslin:GraySS}.
We work with the associated normalized chain complexes to the displayed simplicial abelian groups.

First we show that this map is injective on homology groups. Suppose that $\mcal{W}\in\Z/2\Cor_{k}(X\times\Delta^{n}, Y)$ is a cycle such that $\mcal{W}\otimes\mcal{I}$ is a boundary. Then there is $\mcal{V}\in \Z/2\Cor_{k}((X\times\Delta^{n+1})\wedge\G, Y\wedge\G)$ and 
$\partial_{n+1}(\mcal{V}) = \mcal{W}\otimes\mcal{I}$ and 
$\partial_{i}(\mcal{W}) = 0$ for $0\leq i\leq n$. By Lemma \ref{vwd} there is $N$ such that $\rho_{N}(\mcal{V})$ is defined. By Lemma \ref{vcp} we have that $\rho_{N}(\partial_{i}\mcal{W})$ is defined as well. Moreover by Lemma \ref{vcp} we have
\begin{align*}
& \partial_{i}(\rho_{N}(\mcal{V}))   = \rho_{N}(\partial_{i}(\mcal{V})) = 0, \;\;\; 0\leq i \leq n \\
& \partial_{n+1}(\rho_{N}(\mcal{V})) = \rho_{N}(\partial_{n+1}(\mcal{V})) = \rho_{N}(\mcal{W}\otimes \mcal{I}) = \mcal{W} .
\end{align*}

Therefore $\mcal{W}$ is itself a boundary and so the map on homology is an injection.

Now we show that the map on homology is surjective. Let $\mcal{V}\in \Z/2\Cor_{k}(X\times\Delta^{n}\wedge \G, Y\wedge \G)$ be a cycle (i.e., $\mcal{V}\in \Z/2\Cor_{k}((X\times\Delta^{n})\wedge\G, Y\wedge \G)$ 
satisfies $\partial_{i}(\mcal{V}) = 0$ for $0\leq i \leq n$). Consider the homotopy $\phi=\phi_{\mcal{V}}$ from Lemma \ref{ehp} which satisfies
$$
\phi_{0}-\phi_{1} = \mcal{V}\otimes^{(\tau)}\mcal{I} - \mcal{V}\otimes\mcal{I}.
$$
Applying $\rho_{N}$ (with respect to the second factor of $\G$) and using Lemma \ref{vcp} we have
\begin{align*}
& \rho_{N}(\mcal{V}\otimes \mcal{I}) = \mcal{V} \\
& \rho_{N}(\mcal{V}\otimes^{(\tau)}\mcal{I}) = \rho_{N}(\mcal{V})\otimes \mcal{I}.
\end{align*}
Let $\psi_{N} = \rho_{N}(\phi)$. We have
$$
(\psi_{N})_{0} - (\psi_{N})_{1} = \rho_{N}(\mcal{V})\otimes \mcal{I} - \mcal{V}
$$
and so $\partial_{i}(\psi_{N}) = \rho_{N}(\partial_{i}\phi_{\mcal{V}}) = 0$ (because $\partial_{i}\phi_{\mcal{V}} = \phi_{\partial_{i}\mcal{V}} = 0$). Thus $\psi_{N}\in 
\Cor_{k}(X\times\Delta^{\bullet}_{k}\wedge\G,Y\wedge\G)$ is a cycle. The two restrictions
\begin{equation*}
\Z/2\Cor_{k}(X\times\A^{1}\times\Delta^{\bullet}_{k}\wedge\G, Y\wedge\G)
\to \Z/2\Cor_{k}(X\times\Delta^{\bullet}_{k}\wedge\G,Y\wedge\G)
\end{equation*} 
induced by $0\in \A^{1}$ and $1\in \A^{1}$ induce the same map in homology. Therefore $(\psi_{N})_{0} - (\psi_{N})_{1} = \rho_{N}(\mcal{V})\otimes \mcal{I} - \mcal{V}$ is a boundary in $\Z/2\Cor_{k}(X\times\Delta^{\bullet}_{k}\wedge\G, Y\wedge \G)$. 
\end{proof}

Finally we have the following. Let $F$ be a homotopy invariant sheaf with equivariant transfers. 

\begin{theorem}\label{thm:can}
Let $X$ be a smooth $\Z/2$-scheme. Then
$$
H^{n}_{GNis}(X, C_*\Z_{tr,G}(Y)) = H^{n}_{GNis}(X\wedge \G, C _*\Z_{tr,G}(Y\wedge \mathbb{G})).
$$
\end{theorem}
\begin{proof}
 The argument is formally the same as in the nonequivariant case given all of the machinery developed in the previous sections. For convenience we give some details, focusing on the case of $\G_{m}^{\sigma}$. Consider the projection $\pi:X\times \G_{m}^{\sigma}\to X$. We first consider the Leray spectral sequence (which is convergent as $X$ has bounded cohomological dimension)
$$
E_{2}^{p,q} = H^{p}(X,R^{q}\pi_{*}C_{*}\Z_{tr,G}(Y\wedge \G_{m}^{\sigma})) \Longrightarrow H^{p+q}(X\times \G_{m}^{\sigma},C_{*}\Z_{tr,G}(Y\wedge\G_{m}^{\sigma})).
$$
Write $H^{q}$ for the $q$th cohomology sheaf of the complex $C_{*}\Z_{tr,G}(Y\wedge\G_{m}^{\sigma})$. To compute the complex $\mathbb{R}\pi_{*}C_{*}\Z_{tr,G}(Y\wedge\G_{m}^{\sigma})$ we use the hypercohomology spectral sequence,
$$
E_{2}^{p,q} = R^{p}\pi_{*}H^{q} \Longrightarrow H^{p+q}(\mathbb{R}\pi_{*}C_{*}\Z_{tr,G}(Y\wedge \G_{m}^{\sigma})).
$$
The stalks of $R^{p}\pi_{*}H^{q}$ are $H^{p}(S\times\G_{m}^{\sigma},H^{q})$ where $S$
is a smooth affine semilocal Henselian $G$-scheme over $k$ with a single closed orbit.
By Theorem \ref{thm:Z2hizero} we have $R^{p}\pi_{*}H^{q} = H^{p}_{GNis}(S\times\G_{m}^{\sigma},H^{q}) = 0$ for $p>0$. The spectral sequence thus degenerates and we have 
$$
H^{q}(\mathbb{R}\pi_{*}\Z_{tr,G}(Y\wedge\G_{m}^{\sigma})) = \pi_{*}H^{q}.
$$
The stalks of $\pi_{*}H^{q}$ are $H^{0}_{GNis}(S\times \G_{m}^{\sigma},H^{q})$ which, by Theorem \ref{thm:Z2cont}, 
split into the direct sum
$$
H^{0}_{GNis}(S\times \G_{m}^{\sigma},H^{q}) = H^{q}(S) \oplus H^{q}_{(-\sigma)}(S).
$$
By Proposition \ref{prop:Z2sfhcon} we have that $(\mcal{H}^{q}_{GNis})_{(-\sigma)} = (\mcal{H}^{q}_{(-\sigma)})_{GNis}$. 
Therefore we have 
$$
H^{q}(\mathbb{R}\pi_{*}C_{*}\Z_{tr,G}(Y\wedge \G_{m}^{\sigma})) = H^{q}(C_{*}\Z_{tr,G}(Y\wedge\G_{m}^{\sigma}))\oplus H^{q}(C_{*}\Z_{tr,G}(Y)).
$$
Thus $C_{*}\Z_{tr,G}(Y\wedge\G_{m}^{\sigma})\oplus C_{*}\Z_{tr,G}(Y) \to \mathbb{R}\pi_{*}C_{*}\Z_{tr,G}(Y\wedge \G_{m}^{\sigma})$ is a quasi-isomorphism and therefore
$$
H^{*}(X\times \G_{m}^{\sigma},C_{*}\Z_{tr,G}(Y\wedge\G_{m}^{\sigma})) = 
H^{*}(X,C_{*}\Z_{tr,G}(Y\wedge\G_{m}^{\sigma})) \oplus
H^{*}(X,C_{*}\Z_{tr,G}(Y)),
$$
as required.
\end{proof}

We finish by relating the complexes $C_*\Z_{tr,G}(\G^{\sigma}_{m})$ to the ones introduced in Section \ref{sec:bred}. After a change of coordinates $\P(\sigma \oplus 1)$ can be viewed as $\P^{1}$ with the action $[x:y]\mapsto [y:x]$ and $\G_{m}^{\sigma}$ becomes identified with $\P^{1} - \{[0:1],[1:0]\}$. Consider the Cartesian square in $G\Sm/k$
$$
\xymatrix{
\G_{m}^{\sigma}\times \Z/2 \ar[d] \ar[r] & \A^{1}\times \Z/2 \ar[d]^{\phi} \\
\G_{m}^{\sigma}\ar[r] & \P^1 .
}
$$
The action on $\A^{1}\times \Z/2$ is given by switching the factors and the map $\phi$ sends $(x,e)$ to $[x:1]$ and $(x,\sigma)$ to $[1:x]$. Note that $\phi^{-1}(\{[0:1],[1:0]\}) \iso \{[0:1],[1:0]\}$ is an equivariant isomorphism. In particular, the above square is an equivariant distinguished square.

Recall that we write $S^{\sigma}$ for the topological representation sphere associated to the sign representation. We have $\Z_{top}(S^{\sigma}) = \cone(\Z_{tr,G}(\Z/2) \to \Z)$, where $\Z_{top}(S^{\sigma})$ is as in Example \ref{ex:repsph}.
Using the square above, we obtain a quasi-isomorphism 
$$
(C_{*}\big(\Z_{tr,G}(\G_{m}^{\sigma})/\Z) \otimes_{tr} \Z_{top}(S^{\sigma})\big) \wkeq
C_{*}\big(\Z_{tr,G}(\P(\sigma \oplus 1))/\Z_{tr,G}(\P(\sigma))\big).
$$

For a representation $V$, define  the sheaf with equivariant transfers $\Z_{tr,G}(T^{V})$ by
$$
\Z_{tr,G}(T^{V}) := \Z_{tr,G}(\P(V\oplus 1))/\Z_{tr,G}(\P(V))
$$ 
and similarly for expressions such as $\Z_{tr,G}(X\wedge T^{V})$.
\begin{theorem}\label{thm:lastcan}
 Let $X$ be a smooth $\Z/2$-scheme and $V$ a finite dimensional representation. Then
$$
H^{n}_{GNis}(X,C_*\Z_{tr,G}(Y)) \iso H^{n}_{GNis}(X\wedge T^V, C _*\Z_{tr,G}(Y\wedge T^V)).
$$
\end{theorem}
\begin{proof}
It is enough to treat the case of a one dimensional representation. 
Using Theorem \ref{thm:gniscoh} and a standard spectral sequence argument, one sees that the displayed map of hypercohomology groups can be computed as
\begin{align*}
\Ext^{n}&(C_{*}\Z_{tr,G}(X), C_*\Z_{tr,G}(Y)) \\ &\to \Ext^{n}(C_{*}\Z_{tr,G}(X)\otimes_{tr} C_{*}\Z_{tr,G}(T^V), C_*\Z_{tr,G}(Y)\otimes_{tr} C_{*}\Z_{tr,G}(T^V)),
\end{align*}
 where $\Ext$ is computed in  $D^{-}(G\Cor_{k})$. If $V$ is a trivial representation then we have $C_{*}\Z_{tr,G}(T^V)\wkeq C_{*}(\Z_{tr,G}(\G_{m})/\Z)[1]$ and if $V$ is the sign representation then we have $C_{*}\Z_{tr,G}(T^V)\wkeq C_{*}(\Z_{tr,G}(\G_{m}^{\sigma})/\Z)\otimes_{tr}\Z_{top}(S^{\sigma})$. In either case, as both shift and $\Z_{top}(S^{\sigma})$ are invertible (see Lemma \ref{lem:invert}), the theorem follows  from Theorem \ref{thm:can}.
\end{proof}

\section*{Acknowledgements.}
It is our pleasure to thank Eric Friedlander for many helpful discussions which have greatly influenced this paper. We also thank Aravind Asok, David Gepner, Philip Herrmann, Amalendu Krishna, Kyle Ormsby, and Ben Williams for useful conversations regarding material in Section \ref{sec:enis}.
The first author thanks the University of Oslo and the IPMU for their hospitality 
during the preparation of this work. 
The third author would like to thank the MIT Mathematics Department 
for its hospitality.
\bibliographystyle{amsalpha}
\bibliography{cancellation}

\end{document}